\documentclass[10pt]{amsart}
\usepackage{etex}
\usepackage[latin1]{inputenc}
\usepackage{xspace,amssymb,amsfonts,euscript}
\usepackage{amsthm,amsmath,amscd,stmaryrd,latexsym}
\usepackage{palatino, euscript}

\usepackage{amscd, amssymb,mathtools,bbm,mathdots,arydshln,pinlabel,mathtools,hyperref,youngtab}
\usepackage[nohug]{diagrams}
\usepackage[T1]{fontenc}
\makeatletter

\usepackage{color}
\usepackage{tikz}
\usetikzlibrary{matrix,arrows}

\usepackage{comment}

\newcommand{\MA}[1]{\textcolor[rgb]{1,.5,0}{MA: #1}}

\newtheorem{theorem}[equation]{Theorem}

\newtheorem{lemma}[equation]{Lemma}
\newtheorem{proposition}[equation]{Proposition}
\newtheorem{corollary}[equation]{Corollary}
\newtheorem{conjecture}[equation]{Conjecture}

\newtheorem*{aijRelation}{Proposition \ref{prop-aijRelation}}

\theoremstyle{definition}
\newtheorem{definition}[equation]{Definition}

\theoremstyle{remark}
\newtheorem{example}[equation]{Example}
\newtheorem{remark}[equation]{Remark}

\newtheorem{construction}[equation]{Construction}
\newtheorem{notation}[equation]{Notation}

\numberwithin{equation}{section}

\newcommand{\gr}{\operatorname{gr}}

\newcommand{\inv}{^{-1}}

\newcommand{\Z}{\mathbb{Z}}
\newcommand{\N}{\mathbb{N}}
\newcommand{\Q}{\mathbb{Q}}
\newcommand{\C}{\mathbb{C}}
\newcommand{\e}{\varepsilon}
\renewcommand{\d}{\delta}
\renewcommand{\a}{\alpha}\renewcommand{\b}{\beta}

\newcommand{\End}{\operatorname{End}}
\newcommand{\Sym}{\operatorname{Sym}}
\newcommand{\Endg}{\underline{\operatorname{End}}}
\newcommand{\Homg}{\underline{\operatorname{Hom}}}

\newcommand{\Tor}{\operatorname{Tor}}
\newcommand{\FT}{\operatorname{FT}}

\newcommand{\Id}{\operatorname{Id}}

\newcommand{\ip}[1]{\langle #1 \rangle}
\renewcommand{\matrix}[1]{\begin{bmatrix}#1\end{bmatrix}}
\renewcommand{\aa}{\mathbf{a}}
\newcommand{\bb}{\mathbf{b}}
\newcommand{\xx}{\mathbf{x}}
\newcommand{\yy}{\mathbf{y}}

\newcommand{\wt}{\operatorname{wt}}

\newcommand{\Ext}{\operatorname{Ext}}
\renewcommand{\deg}{\operatorname{deg}}

\newcommand{\Tot}{\operatorname{Tot}}

\newcommand{\Cone}{\operatorname{Cone}}

\renewcommand{\sl}{\mathfrak{sl}}

\newcommand{\colim}{\operatorname{colim}}

\newcommand{\SBim}{\mathbb{S}\mathbf{Bim}}
\newcommand{\pr}[1]{P_{1^{#1}}}
\newcommand{\one}{\mathbbm{1}}

\newcommand{\K}{\mathcal{K}}
\newcommand{\HT}{\operatorname{HT}}
\newcommand{\Ch}{\operatorname{Ch}}
\newcommand{\hocolim}{\operatorname{hocolim}}

\newcommand{\ww}{\underline{w}}

\newcommand{\HH}{\operatorname{HH}}
\newcommand{\HHH}{\operatorname{HHH}}

\newcommand{\del}{\partial}

\newcommand{\holim}{\operatorname{holim}}
\newcommand{\IC}{\mathcal{I}}
\newcommand{\JC}{\mathcal{J}}

\newcommand{\KC}{\mathcal{K}}
\newcommand{\HC}{\mathcal{H}}
\renewcommand{\H}{\mathcal{H}}

\newcommand{\KBC}[1]{\KC^-(\SBim_{#1})}
\newcommand{\homfly}{HOMFLYPT }

\newcommand{\RbimR}{R\text{-mod-}R}

\begin{document}

\title[Categorified Young symmetrizers and torus links II]{Categorified Young symmetrizers and stable homology of torus links II}
\author[Michael Abel and Matthew Hogancamp]{Michael Abel and Matthew Hogancamp}
\address{Department of Mathematics, Duke University, Durham, NC 27708}
\email{maabel\char 64 math.duke.edu}
\address{Department of Mathematics, Indiana University, Bloomington, IN 47401}
\email{mhoganca\char 64 indiana.edu}


\begin{abstract}We construct complexes $\pr n$ of Soergel bimodules which categorify the Young idempotents corresponding to one-column partitions.  A beautiful recent conjecture \cite{GoRa15} of Gorsky-Rasmussen relates the Hochschild homology of categorified Young idempotents with the flag Hilbert scheme.  We prove this conjecture for $\pr n$ and its twisted variants. We also show that this homology is also a certain limit of Khovanov-Rozansky homologies of torus links. Along the way we obtain several combinatorial results which could be of independent interest. 
\end{abstract}

\maketitle


\section{Introduction}
\label{sec-intro}

The computation of Khovanov-Rozansky homology of torus links is a challenging problem which is made even more interesting by the conjectures in \cite{ORS12,GORS12} relating these homology groups with certain Hilbert schemes and Cherednik algebras.  It seems a more tractable problem to compute the stable limit of the homology of the $(n,m)$ torus link as $m\to\infty$.  The stable limit of homologies of Khovanov-Rozansky homology has a longer history, being conjectured to exist by \cite{DGR06} and proven to exist (in the $\sl_2$ case, that is, the original Khovanov homology) in \cite{Stosic07}.  Rozansky showed that in the $\sl_2$ case, the stable limit can be computed from a categorified Jones-Wenzl projector \cite{Roz10a}.  Recently, the second author \cite{Hog15} proved existence of this stable limit in the triply-graded homology, and proved a stable version of conjectures in \cite{ORS12,GORS12}.

In this paper we investigate a second kind of stable limit of Khovanov-Rozansky homology of torus links.  In the $\sl_2$ setting, the existence of two different stable limits was due to Rozansky \cite{Roz10b}; this second limit computes Hochschild homology of Khovanov's ring $H^n$ \cite{Kh02}.  In the triply graded case, we will see that the second stable limit computes a certain kind of Hochschild homology associated to polynomial rings.  This stable limit is also explicitly related to Hilbert schemes by a recent conjecture of Gorsky-Rasmussen \cite{GoRa15} (see Conjecture \ref{conj-GR}).  For  a comparison of stable homologies see \S \ref{subsec-introCompare}.

Fix an integer $n \geq 1$ and let $R = \Q[x_1,\ldots,x_n]$ be graded by setting $\deg(x_k) = 2$. Note that $W := S_n$ acts on $R$ by permuting variables. Let $R^W \subset R$ denote the ring of symmetric polynomials.  Observe that $f \otimes g \mapsto fg$ defines an algebra map $R \otimes_{R^W} R \to R$, so that $R$ is a graded $R\otimes_{R^W} R$-module.

\begin{definition}\label{def-introP}
Let $R \to \pr{n}^\vee$ denote an injective resolution of $R$, thought of as a graded $R \otimes_{R^W} R$-module.
\end{definition}

We note for future reference that $R\otimes_{R^W} R$ is self-injective.  The complex $\pr{n}^\vee$ is graded infinite, being supported in all non-negative homological degrees. This complex categorifies the one-column Young symmetrizer in a sense explained below. 

Let $\Ch(\RbimR)$ denote the category of complexes of finitely generated, graded $(R,R)$-bimodules.  Our first main theorem says that $\pr{n}^\vee$ is related to a certain braid group action on $\Ch(\RbimR)$. Associated to each $n$-strand braid, Rouquier \cite{Rou06} defines a complex $F(\beta)\in \Ch(\RbimR)$ such that $F(\b)\otimes_R F(\b')\simeq F(\b \b')$.  We have:

\begin{theorem}\label{thm-introLimit1}
	Let $\FT \in \Ch(\RbimR)$ denote the Rouquier complex associated to the (positive) full twist braid. Then $\{\FT^{\otimes k}\}_{k=0}^\infty$, after the appropriate grading shifts, can be made into a direct system with homotopy colimit $\pr{n}^\vee$.
\end{theorem}

\begin{remark}\label{rmk-introDuality}
Dually, one may consider a complex $\pr{n}$ which is defined to be a projective resolution $\pr{n} \to R$ viewed as a graded $R\otimes_{R^W}R$-module.  Analogously, if we let $\FT\inv$ denote the Rouquier complex of the negative full twist braid then $\{(\FT\inv)^{\otimes k}\}_{k=0}^\infty$ can made into an inverse system with homotopy limit $\pr n$, after the appropriate grading shifts. We will prove in \S \ref{subsec-explicitCx} that the duality functor $(\cdot)^\vee: \Ch(\RbimR) \to \Ch(\RbimR)$ sends $\pr n$ to $\pr n^\vee$.

Our notation prefers complexes which are bounded above, in order to be consistent with \cite{Hog15}.  However, in this paper the complex $\pr{n}^\vee$ is often nicer than $\pr{n}$, since $\pr n^\vee$ is an algebra in the homotopy category $\K^+(\RbimR)$, while $\pr n$ is a coalgebra.  This is discussed in Remark \ref{rmk-introAlgebras}.  The algebra structure makes writing our following main results easier. See Example \ref{ex-nondual} for a sample computation in the non-dual, $n=2$ case.
\end{remark}

Let us recall the construction of Khovanov and Rozansky's triply-graded link homology given in \cite{Kh07}.  Hochschild cohomology of bimodules defines a functor $\HHH$ from $\Ch(\RbimR)$ to the category of triply-graded vector spaces, and $\HHH(F(\beta))$ is a well-defined invariant of the oriented link $\hat{\beta}$, up to isomorphism and overall shift of triply-graded vector spaces.  This invariant is called Khovanov-Rozansky homology, hereafter referred to as \emph{KR homology}.
\begin{remark}\label{rmk-introHHH}
Actually, in \cite{Kh07} the results are stated in terms of Hochschild homology.  In this paper, as in \cite{Hog15}, we prefer Hochschild cohomology to Hochschild homology; there is not much difference since for polynomial rings the two are isomorphic up to regrading.  Our convention will ensure that $\HHH$ applied to our projector is a graded commutative algebra, rather than an algebra up to regrading.
\end{remark}
Theorem \ref{thm-introLimit1} allows us to relate KR homology of torus links with $\HHH(\pr{n}^\vee)$:

\begin{theorem}\label{thm-introFullTwist}
The KR homology of the $(n,nk)$ torus links stablizes as $k\to\infty$.  The stable limit is isomorphic to $\HHH(\pr{n}^\vee)$.
\end{theorem}
In order to state our results for the $(n,nk+m)$ torus links, we first introduce an action of the symmetric group. Let $w \in S_n$, then for any complex $C \in \Ch(\RbimR)$, let $w(C)$ denote the complex obtained from $C$ by twisting the right $R$-action by $w\inv$, that is, $f \cdot c \cdot g = fcw\inv(g)$.  We refer to $w(\pr{n}^\vee)$ as a \emph{twisted projector}.  The connection between this action of the symmetric group and the action of the braid group by Rouquier complexes is established in \S \ref{subsec-rouquier}.  An elementary property of Hochschild cohomology implies that $\HHH(w(\pr{n}^\vee))$ depends only on the conjugacy class (cycle type) of $w$.  Using this, we have:

\begin{theorem}\label{thm-introStableExistence}
Fix integers $1 \leq m \leq n$, and let $w \in S_n$ be an $n$-cycle. The triply-graded KR homology of the $(n,nk+m)$-torus links have a stable limit as $k \to \infty$, given by $\HHH(w^m(\pr{n}^\vee))$. This limit depends only on the number of components $\gcd(m,n)$, up to isomorphism.
\end{theorem}

 Thus, the problem of computing stable homology of torus links reduces to the computation of $\HHH(w(\pr{n}^\vee))$ for the appropriate $w\in S_n$.  To compute this we first reduce to the computation of the Hochschild degree zero part $\HHH^0(w(\pr{n}^\vee))$.  Here and below, we use the same notation for both exterior and polynomial algebras, and will distinguish them by declaring the variables as odd or even, respectively.

\begin{proposition}\label{prop-introHHHsimp}
Let $\xi_k$ denote an odd variable of degree $q^{-2k}a$, where we use $a$ to denote Hochschild degree.  If $C\in \Ch(\RbimR)$ is a complex whose chain bimodules are direct sums of copies of $R\otimes_{R^W} R$ with shifts, then $\HHH(C) = \Lambda[\xi_1,\ldots,\xi_n]\otimes_\Q \HHH^0(C)$, where $\HHH^0(C)$ is, up to a regrading (see Remark \ref{rmk-introHHH}), the functor which identifies the left and right $R$-actions then takes homology.  That is, $\HHH^0(C)=H(R\otimes_{R\otimes_{R^W} R} C)$.
\end{proposition}
This is proven in \S \ref{subsec-computations}.   The following beautiful conjecture  was communicated to us by E.~Gorsky and J.~Rasmussen \cite{GoRa15}:

\begin{conjecture}\label{conj-GR}
For each $1\leq i<j\leq n$, let $v_{ij}$ denote an even variable of degree $\deg(v_{ij})=q^{2i-2j-2}t^2$.  Let $E$ denote the ring $E = R[v_{ij}]/J$ where $J$ is the ideal generated by the entries of the commutator of the following formal matrices:
	\[
	X = \matrix{
		x_1	&	1		&	0		&	\cdots	&	0& 0\\
		0		&	x_2	&	1		&  \cdots & 0& 0 \\
		0		&	0		&	x_3		&  \cdots & 0& 0 \\
		\vdots		&	\vdots		&	\vdots		&  \ddots & \vdots & \vdots  \\
		0		&	0		&	0		&  \cdots & x_{n-1} &1 \\
		0		&	0		&	0		&  \cdots & 0 & x_n
	}
	,\ \
	V = \matrix{
		0	&	v_{12}		&	v_{13}		&	\cdots	& v_{1,n-1} &	v_{1n}\\
		0		&	0	&	v_{23}		&  \cdots & v_{2,n-1} & v_{2n} \\
		0		&	0		&	0		&  \cdots & v_{3,n-1} & v_{3n} \\
		\vdots		&	\vdots		&	\vdots		&  \ddots &\vdots  & \vdots  \\
		0		&	0		&	0		&  \cdots & 0 & v_{n-1,n} \\
		0		&	0		&	0		&  \cdots & 0&  0
	}
	\]
Then up to an overall grading shift, $\HHH^0(w(\pr{n}^\vee))$ is isomorphic to the quotient of $E$ in which we identify $x_{w(i)}$ with $x_i$ for all $1\leq i\leq n$.
\end{conjecture}


More generally, to each Young tableau $T$, Gorsky and Rasmussen \cite{GoRa15} associate an algebra $E_T$ coming from the flag Hilbert scheme on $\C^2$. The \emph{flag Hilbert Scheme} $\text{FHilb}^n(\C^2)$ of $n$ points on $\C^2$ is the moduli space of flags of ideals $$\text{FHilb}^n(\C^2) =\{\C[x,y] = I_0 \supset I_1 \supset \cdots \supset I_n | \dim(I_k/I_{k+1})=1\}.$$
 $\text{FHilb}^n(\C^2,\ell)$ is defined as the subscheme of $\text{FHilb}^n(\C^2)$ such that all $I_k$ are set-theoretically supported on the line $x=0$.

There is a natural $\C^* \times \C^*$-action on $\text{FHilb}^n(\C^2,\ell)$. The fixed points $z_T$ of this action are parameterized by standard Young tableaux of size $n$. In \cite{GoRa15}, Gorsky and Rasmussen explicitly construct the local algebras $E_T$ of $z_T$. They conjecture that in fact, $E_T$ is isomorphic to $\End(P_T)$ where $P_T$ is the categorified Young symmetrizer associated to the Young tableau $T$.

One of the main contributions of this paper is an exact combinatorial description of the the projector $\pr{n}^\vee$ (see \S \ref{subsec-explicitCx} and also \S \ref{subsec-introComb} of this introduction).  Using this description we can prove the above conjecture:

\begin{theorem}\label{thm-introTwistedE}
Conjecture \ref{conj-GR} is true.  Furthermore, there are isomorphisms of bigraded rings
\[
\HHH^0(\pr{n}^\vee)\cong E\cong \Endg(\pr{n}^\vee)\cong R[u_2,\ldots,u_n]/I
\]
where $I$ is the ideal generated by the elements $\sum_{i=1}^j u_i a_{i,j}(\xx,w(\xx))$ for $1\leq j\leq n$, $a_{ij}(\xx,\yy)$ are certain polynomials explained below, and  $\Endg(\pr{n}^\vee)$ here denotes the bigraded ring of bihomogeneous chain endomorphisms of $\pr{n}^\vee$, modulo homotopies. The degrees are $\deg(u_k)=t^2q^{-2k}$.
\end{theorem}
For more of an explanation of the relation between Hochschild cohomology of $\pr{n}^\vee$ and its ring of endomorphisms, see \S \ref{subsec-EndPn}.

The special polynomials $a_{ij}(\xx,\yy)$ are defined in Definition \ref{def-aij}, and are essentially the double Schubert polynomials $\mathfrak{S}_w(\xx,\yy)$ (See \cite{Mac91} for more information).  More precisely, if $w_{i,n}$ denotes the cycle $(n, n- i+1,n-i+2, \ldots , n-1) \in S_n$ and $w_0 \in S_n$ is the longest word, then $a_{in} = (-1)^i\mathfrak{S}_{w_{i,n}}(w_0(\xx),\yy).$  In Proposition \ref{prop-theRels} we explicitly give the relation between $a_{ij}$ and $v_{ij}$.  We state explicitly the following important special case of torus \emph{knots}, in which the stable homology is a free superpolynomial algebra:

\begin{corollary}\label{cor-introCycle}
If $w\in S_n$ is an $n$-cycle, then $$\HHH(w(\pr{n}^\vee)) \cong R[u_2,\ldots,u_n,\xi_1,\ldots,\xi_n]/I,$$ where $I=I_w$ is the ideal generated by $x_i-x_{i-1}$ $(1\leq i\leq n-1)$. In particular, this homology is isomorphic to a stable limit of KR homologies torus knots.  The degrees are $\deg(u_k)=t^2q^{-2k}$ and $\deg(\xi_k)=aq^{-2k}$.
\end{corollary}

In \cite{AbHog15}, we use $\pr{n}$ to give a categorification of the $\Lambda^n$-colored \homfly polynomial. This polynomial was first categorified by Webster and Williamson using a geometric framework \cite{WW09}.  These two theories are different however, as the homology of the unknot in the theories are quite different.  

We now explore the different pieces of this construction in more detail, including a combinatorial description of $\pr n^\vee$. Along the way we compare our results to the results of \cite{Hog15}.

\subsection{Hecke algebras, Young symmetrizers, and categorification}
\label{subsec-introYoung}

The Hecke algebra of $S_n$ (or briefly Hecke algebra) is a $q$-deformation of the symmetric group algebra $\C[S_n]$.  The classical Young symmetrizers are certain idempotent elements $p_T \in \C[S_n]$, indexed by Young tableaux of size $n$. In terms of representation theory, if $T$ is of shape $\lambda$, multiplication by $p_T$ is a projection of the regular representation $\C[S_n]$ onto the irreducible representation $V_\lambda$. In $\H_n$ there are natural $q$-analogues to the Young symmetrizers serving the same role \cite{Gyo86,AistonMorton98}, and we will also denote these by $p_T$. The Young symmetrizers also play an important role in quantum topology, where ``cabling and inserting $p_T$'' defines the $\lambda$-colored HOMFLYPT polynomial. One also has the central idempotents $p_\lambda = \sum p_T$, where the sum is over tableaux with shape $\lambda$.

\begin{remark}
For each integer $N\geq 1$, one can consider the quotient of $\HC_n$ by the ideal generated by the $p_{\lambda}$, where $\lambda$ has more than $N$ parts.  We call this the $\sl_N$ quotient, since the resulting algebra is isomorphic to $\End_{U_q(\sl_N)}(V^{\otimes n})$, where $V$ is $q$-deformation of the standard $N$-dimensional representation of $\sl_N$.  In case $N=2$, this is the Temperley-Lieb algebra.  In this way, Young symmetrizers act as projection operators on $V^{\otimes n}$.  The one-row idempotent $p_{(n)}$ acts as projection onto the $q$-symmetric power $\Sym^n(V)$, while the one-column idempotent $p_{1^n}$ acts as projection onto the $q$-exterior power $\Lambda^n(V)$, which is zero unless $N\geq n$.  Note that $p_{1^n}$ is zero in the $\sl_N$ quotient unless $N\geq n$.  Thus, our projectors relate to Rozansky's \cite{Roz10b} only in the cases $1\leq n\leq 2=N$.
\end{remark}

In this paper, we categorify the Young symmetrizer for the one-column partition $p_{1^n}$. This $p_{1^n}$ can be described as the unique multiple of $b_{w_0}$ which is idempotent.  Here, $b_{w_0}$ is the Kazhdan-Lusztig basis element corresponding to the long-\-est word $w_0 \in S_n$ and can be described in terms of the standard basis of $\H_n$ by $b_{w_0} = \sum_{w \in S_n} q^{\ell(w_0)-\ell(w)}T_w.$  Equivalently, $p_{1^n}$ is the unique multiple of $b_{w_0}$ such that $1-p_{1^n}$ annihilates $b_{w_0}$ from the left and right.

The category $\SBim_n$ of Soergel bimodules is a certain category of graded $(R,R)$-bimodules (explained more in \S \ref{subsec-soergel}) which categorifies the Hecke algebra $\H_n$.  Roughly what this means is that to each bounded complex $C$ of Soergel bimodules, there is a well-defined Euler characteristic $\chi(C)\in\H_n$ which satisfies a number of properties.  For instance:
\begin{enumerate}
\item $\chi(C\oplus D)=\chi(C)+\chi(D)$
\item $\chi(C\otimes_R D) = \chi(C)\chi(D)$
\item $\chi(C\ip{i}(j))=(-1)^iq^j \chi(C)$, where $\ip{i}$ and $(j)$ are homological and $q$-degree degree shifts, respectively.
\item If $C$ and $D$ are chain homotopy equivalent $C\simeq D$, then $\chi(C)=\chi(D)$
\item $\H_n$ is the $\Q(q)$-span of the Euler characteristics.
\end{enumerate}
We will say that a complex categorifies its Euler characteristic, and that a homotopy equivalence $C\simeq D$ categorifies the identity $\chi(C)=\chi(D)$, and so on.  Most of the complexes we consider are infinite (but bounded in one direction); nonetheless, their Euler characteristics will be well-defined as formal power series.  We will not make this precise, as the notion of Euler characteristic serves only to motivate our use of the word ``categorification.''

The indecomposable Soergel bimodules $B_w$ are indexed by permutations $w\in S_n$, and they categorify the Kazhdan-Lusztig basis $\{b_w\}\subset \H_n$.  As a special case, the bimodule $B_{w_0} := R \otimes_{R^W} R(-\ell)$ categorifes the element $b_{w_0}$.  Here $\ell = \frac{1}{2}n(n-1)$ is the length of the longest element of $S_n$, and $(-\ell)$ is the grading shift which places $1\otimes 1$ in degree $-\ell$.   We will let $\K^-(\SBim_n)$ denote the homotopy category of chain complexes over $\SBim_n$ which are bounded from above. 

\begin{theorem}\label{intro-Pn}
	There exists a bounded from above chain complex $\pr{n} \in \K^-(\SBim_n)$ and a map $\e: \pr{n} \to R$ such that
	\begin{enumerate}
		\item[(P1)] $\pr{n}$ is homotopy equivalent to a chain complex whose chain bimodules are direct sums of shifted copies of $B_{w_0}$.
		\item[(P2)] $\Cone(\e) \otimes_R B_{w_0} \simeq B_{w_0} \otimes_R \Cone(\e) \simeq 0$.
		\item[(P3)] $\pr{n} \otimes_R \pr{n} \simeq \pr{n}$.
	\end{enumerate}
	If another pair $(P',\e')$ satisfies conditions (P1) and (P2), then $P'$ is canonically homotopy equivalent to $\pr{n}$; the canonical homotopy equivalence $\Phi: P' \to \pr{n}$ is characterized up to homotopy by $\e \simeq \Phi \circ \e'$.
\end{theorem}

The axioms (P1) and (P2) ensure that $\e\otimes \Id_{\pr{n}}$ and $\Id_{\pr{n}}\otimes \e$ are homotopy equivalences, and thus imply (P3).  In the language of \cite{Hog15}, such an object is called a \emph{counital} idempotent, and a number of important properties follow immediately from the general theory of such objects (for instance, uniqueness).  We prove in \S \ref{subsec-axioms} that the projective resolution definition from Definition \ref{def-introP} satisfies (P1) and (P2); this gives the existence part of the proof of the above theorem. 

\subsection{A comparison of stable homologies}
\label{subsec-introCompare}
Define $\FT\in \K^-(\SBim_n)$ to be the Rouquier complex associated to the full twist.  Suppose we are given a morphism $\a:R\rightarrow G(\FT)$, where $G$ is a grading shift.  Tensoring $\a$ with $G^k(\FT^{\otimes k})$ gives a map $\a_k: G^k(\FT^{\otimes k}) \rightarrow G^{k+1}(\FT^{\otimes k+1})$.  In this way, maps $R\rightarrow G(\FT)$ yield directed systems $\{G^k(\FT^{\otimes k}) \}_{k=0}^\infty$.  One might refer to limits (or, more precisely, homotopy colimits) of such directed systems as  ``infinite full twists,'' since they are limits of powers of $\FT$ with shifts.  Motivated by work of Rozansky, one may expect that infinite full twists should be intimately related with categorified Young symmetrizers.  This is indeed the case.  

There is a canonical map $R\rightarrow \FT(n(n-1))\ip{-n(n-1)}$ which is the inclusion of the degree zero chain group.  In \cite{Hog15}, the second author showed that the resulting directed system has a homotopy colimit given by $P_n$, and that this complex categorifies the one-row Young symmetrizer.

In this paper we construct a map $R\rightarrow \FT(-n(n-1))$ so that the resulting homotopy colimit categorifies the one-column Young symmetrizer.  Actually this limit is the dual of $\pr{n}$; to obtain $\pr{n}$ one should instead take homotopy limit of an \emph{inverse} system involving $\FT\inv$ (see below).  We will first discuss the $n=2$ case.  Let $F = F(\sigma_1)$ denote the Rouquier complex associated to the (positive) elementary braid generator on two strands, so that $\FT=F^{\otimes 2}$.  A straightforward calculation gives the following, omitting grading shifts for clarity. 
$$F^{\otimes k} \ \simeq \ B \xrightarrow{\phi_\pm} B \xrightarrow{\phi_\mp} \cdots \xrightarrow{\phi_-} B \to R .$$
Here $F^{\otimes k}$ is a complex of length $k+1$, the sign alternates between $+$ and $-$ or $\pm = (-1)^{k-1}$, $B \in \SBim_2$ is the nontrivial indecomposable Soergel bimodule, and $\phi_\pm$ are some endomorphisms of $B$. If we shift so that $R$ appears in homological degree zero, then these complexes approach a well defined limit $P_2$ (precisely, a homotopy colimit).  If we instead shift $F^{\otimes k}$ so that the left-most $B$ appears in homological degree zero, then this sequence does not approach a well-defined limit because of the alternating differential $\phi_\pm$. However, there are two ``convergent subsequences'' which approach the limits (homotopy colimits) $\pr{2}^\vee$ and $s(\pr{2}^\vee)$. Here $\vee$ is the duality functor $\vee: \K^-(\SBim_2) \to \K^+(\SBim_2)$, and $s\in S_2$ is the nontrivial permutation.

Similarly, for negative powers of $F$, one has 
$$F^{\otimes -k} \ \simeq \ R\to B \xrightarrow{\phi_-} B \xrightarrow{\phi_+}  \cdots \xrightarrow{\phi_\pm} B.$$
Various homotopy limits of these complexes will produce $P_2^\vee$, $\pr{2}$, and $s(\pr{2}^\vee)$.

To describe the general situation, first let $X=\sigma_{n-1} \cdots\sigma_2\sigma_1$ denote the positive braid lift of the $n$-cycle $w=(n,n-1,\ldots,2,1)\in S_n$, so that the closure of $X^{\pm k}$ is the $(n, \pm k)$ torus link.  The following table summarizes the basic information regarding the various stable limits of Khovanov-Rozansky homology:

\begin{center}
\begin{tabular}{|c|c|c|c|}
\hline
Projector: & (Co)limit of: & Boundedness: & Algebra structure \\
\hline
$P_n$ & $\hocolim_k F(X^k)$ & above & algebra \\
$P_{1^n}^\vee$ & $\hocolim_k F(X^{nk})$ & below  & algebra \\
$w^m(P_{1^n}^\vee)$ & $\hocolim_k F(X^{nk+m})$ & below & N/A \\
$P_n^\vee$ & $\holim_k F(X^{-k})$ & below & coalgebra \\
$P_{1^n}$ & $\holim_k F(X^{-nk})$ & above & coalgebra \\
$w^m(P_{1^n}^\vee)$ & $\holim_k F(X^{-nk+m})$ & above & N/A \\
\hline
\end{tabular}
\end{center}
\begin{remark}\label{rmk-introAlgebras}
In the right-most column above, the algebra structure is implied by the general theory of unital and counital idempotents \cite{Hog15}.  The computations of $\HHH$ tend to be nicer to write down in the case of algebras; this is the reason why Theorem \ref{thm-introTwistedE} is stated in terms of the positive torus links $T_{n,nk+m}$, rather than $T_{n,-nk-m}$.
\end{remark}

Recently, the second author in  \cite{Hog15} computed the (bounded above) stable KR homology of the $(n,k)$ torus links:

\begin{theorem} \label{thm-introStblSym}
	There is an isomorphism of triply graded algebras $$\HHH(P_n) \simeq \Q[U_1,\ldots,U_n,\Xi_1,\ldots,\Xi_n]$$ where the $U_k$ are even indeterminates of tridegree $t^{2-2k}q^{2k}$ and the $\Xi_k$ are odd indeterminates of tridegree $t^{2-2k}q^{2k-4}a$. This homology is isomorphic to a colimit of KR homologies of the torus links $(n,k)$ as $k\to \infty$.
\end{theorem}
\begin{remark}
This algebra is isomorphic to the algebra for $\HHH(w(\pr{n}^\vee))$ in Corollary \ref{cor-introCycle} after regrading.   More specifically, introduce new variables $q_1=q^2$, $t_1=t^2q^{-2}$, and $a_1=aq^{-2}$.  Then in Theorem \ref{thm-introStblSym}, the even variables in have degrees $q_1 t_1^{1-k}$ and the odd variables have degrees $a_1 t_1^{1-k}$.  In Corollary \ref{cor-introCycle}, the even variables have degrees $t_1q_1^{1-k}$ and the odd variables have degrees $a_1q^{1-k}$.  Clearly swapping $t_1$ with $q_1$ exchanges these degrees.  This demonstrates a mirror symmetry which is conjectured to exist more generally in \cite{GuSt12}.
\end{remark}

\subsection{Combinatorial description of $\pr{n}$}
\label{subsec-introComb}

Recall from the definition that $\pr{n}$ is a resolution of $R$ by free $R \otimes_{R^W} R$-modules. In \S \ref{subsec-explicitCx} we are able to construct an \emph{explicit} such resolution, which is essential for the computation for $\HHH(\pr{n}^\vee)$. We briefly discuss this construction in this section. 

Put $\xx = \{x_1,\ldots ,x_n\}$, $\yy = \{y_1,\ldots,y_n\}$ and let $I \subset \Q[\xx,\yy]$ denote the ideal generated by $e_k(\xx) - e_k(\yy)$, for $1 \leq k \leq n$, where $e_k$ denotes the $k$th elementary symmetric function. It will be more convenient to consider $(R,R)$-bimodules as $\Q[\xx,\yy]$-modules for this section. Note that $B_{w_0} \simeq \Q[\xx,\yy]/I$ up to a grading shift of $-\ell$, where $\ell = \frac{n(n-1)}{2}$ is the length of $w_0$.  

To motivate the construction, we first remark that the decategorified Young symmetrizer $p_{1^n}$ equals $\frac{1}{[n]!}b_{w_0}$, where $[n]!=[n][n-1]\cdots [2]$, and $[k]=\frac{q^k-q^{-k}}{q-q\inv}$ is the quantum integer.  This is a $q$-analogue of a standard identity in the symmetric group.  In other words:
\begin{equation}\label{eq-introMultiplicity0}
p_{1^n} = q^{-\ell}\frac{(1-q^{-2})(1-q^{-2})\cdots (1-q^{-2})}{(1-q^{-2})(1-q^{-6})\cdots (1-q^{-2n})} b_{w_0}
\end{equation}
To categorify this, we will explicitly construct a complex $M_n$ in which $B_{w_0}$ appears with graded multiplicity
\begin{equation}\label{eq-introMultiplicity}
q^{-\ell}\frac{(1+tq^{-2})(1+tq^{-2})\cdots (1+tq^{-2})}{(1-t^{2}q^{-2})(1-t^{2}q^{-3})\cdots (1-t^{2}q^{-2n})}
\end{equation}

Factors of the form $(1-t^{2}q^{-2k})\inv$ correspond to the action of an even variable $u_k$ of degree $t^{2}q^{-2k}$, and factors of the form $(1+tq^{-2})$ are captured by odd variables $\theta_i$ of degree $tq^{-2}$.  Thus, to start with, we consider the bigraded algebra
\[
A_n:=\Q[\xx,\yy,u_1,u_2,\ldots,u_n,\theta_1,\theta_2,\ldots\theta_n]
\]
with variables $u_k,\theta_k$ as above.  Then we form the quotient $M_n=A_n/IA_n$, and we shift the bigradings so that $\bar{1}\in M_n$ lies in $q$-degree $-2\ell$.  By construction, the bigraded object $M_n$ will a direct sum of infinitely many copies of $B_{w_0}$, with graded multiplicity (\ref{eq-introMultiplicity}).  

Now we put differentials on $A_n$ and $M_n$ so that $A_n$ is a differential bigraded algebra, and $M_n$ is a differential bigraded $A_n$-module.  To do this, in Definition \ref{def-aij} we introduce polynomials $a_{ij}(\xx,\yy)$, ($1\leq i\leq j\leq n$) such that $\sum_{j=i}^n a_{ij}(\xx,\yy)(y_j-x_j) = 0$ modulo $I$.  Then we let $d_A$ be the unique $\Q[\xx,\yy]$-linear differential on $A_n$ given by
\begin{equation}\label{def-introDA}
\begin{aligned}
d_A(u_k) &= 0 \\
d_A(\theta_k) &= \sum_{i=1}^ka_{ik}(\xx,\yy)u_i\\
\end{aligned}
\end{equation}
together with the graded Leibniz rule.   Note that $M_n$ is generated as an $A_n$-module by a single element $\overline{1}\in M_n$, where $\overline{1}$ denotes the image of $1\in A_n$ under the quotient $A_n\rightarrow M_n$.  Thus, the differential $d_M$ is uniquely determined by $d_M(\overline{1})$ together with the graded Leibniz rule for dg modules: $d_M(am)=d_A(a)m\pm a d_M(m)$, with the sign given by the homological degree of $a$.  With this in mind, we let $d_M$ be the unique differential on $M_n$ such that
\[
d_M(\bar{1}) = \sum_{i=1}^n(y_i-x_i)\theta_i\bar{1}
\]
together with the graded Leibniz rule.  The proof that $d_M^2=0$ is straightforward, and relies on the defining property of the $a_{ij}(\xx,\yy)$.

\begin{theorem}\label{thm-introPnAsDGM}
	$M_n \simeq \pr{n}^\vee$.
\end{theorem}
To prove this we construct explict quasi-isomorphisms $\Q[\xx]\cong M_1\rightarrow M_2\rightarrow \cdots \rightarrow M_n$, from which we conclude that $M_n$ is an injective resolution of $\Q[\xx]$, regarded as a graded $\Q[\xx,\yy]/I$-module.  Then the dual of this complex is precisely $\pr{n}$, by Remark \ref{rmk-introDuality}.

The presence of odd variables $\theta_i$ and the structure of the differentials $d_M$ and $d_A$ give the complex various 2-periodicities. This is a common phenomenon in commutative algebra \cite{Eis80}, and can be encoded with an object called a matrix factorization. In \S \ref{subsec-matrixFact} we use this construction to give a proof of Theorem \ref{thm-introTwistedE}.

\begin{remark}
There are reduced versions $\~A_n\simeq A_n$ and $\~M_n\simeq M_n$, obtained by setting the variables $u_1,\theta_1$ equal to zero.  This has the effect of canceling the $(1-q^2)$ factors in (\ref{eq-introMultiplicity0}) Below, we consider the reduced case.
\end{remark}

We now end this section by studying the case $n=3$. Let $s,t\in S_3$ be the simple transpositions, so that $w_0=sts$ is the longest word, and set $B := B_{sts}$.  For the rest of this discussion let $p_{ij}=y_i-x_j$. Unpacking the definitions, including those from Definition \ref{def-aij}, we see that the differential on $\~M_3$ is determined by
\begin{equation}
\begin{aligned}
d_M(\bar{1}) &= p_{22}\theta_2\bar{1}+p_{33}\theta_3\bar{1}\\
d_M(\theta_2\bar{1}) &= p_{12}u_2\bar{1} - p_{33}\theta_2\theta_3\bar{1}\\
d_M(\theta_3\bar{1}) &= (p_{23}+p_{12})u_2\bar{1} + p_{23}p_{13}u_3\bar{1} + p_{33}\theta_2\theta_3\bar{1}\\
d_M(\theta_2\theta_3\bar{1}) &= p_{12}u_2\theta_3\bar{1}-p_{23}p_{13}u_3\theta_2\bar{1}-(p_{12}+p_{23})u_2\theta_2\bar{1}\\
\end{aligned}
\end{equation}
Recall that we have placed $\bar{1}$ in $q$-degree $-6$.  Since $B=(\Q[\xx,\yy]/I)(-3)$, taking into account the degrees of $\theta_i$ and $u_i$, we see that $\~M_3\cong \pr{3}^\vee$ is the total complex of the below ``perturbed'' double complex:

\begin{equation}
\begin{minipage}{122mm}
\begin{center}
\includegraphics[width=122mm]{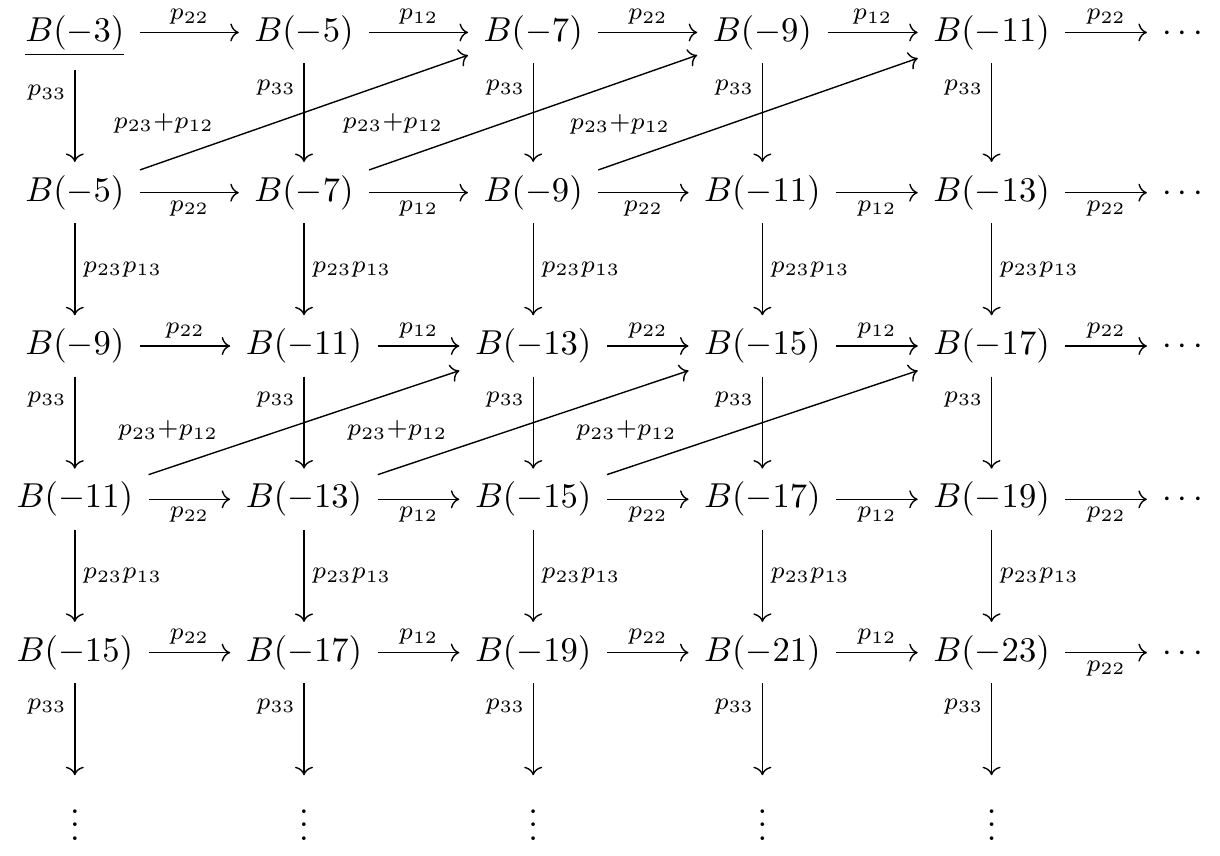}
\end{center}
\end{minipage}
\end{equation}

The double-peridiocity which is evident in the above diagram is realized by the action of the even variables $u_2,u_3$.  The ``fundamental region''
\begin{diagram}
	B(-3) & \rTo^{p_{22}} & B(-5)\\
	\dTo<{p_{33}} & & \dTo>{p_{33}}\\
	B(-5) & \rTo_{p_{22}} & B(-7)\\
\end{diagram}
is spanned by the products of the odd variables.  More precisely, this square is the image in $\tilde{M}_3$ of the exterior algebra $\Q[\xx,\yy,\theta_2,\theta_3]$, which can also be described as the Koszul complex of the sequence $p_{22},p_{33}$ acting on $B$. To write down the complex $\pr{3}$, we reverse the arrows and change all grading shifts $(-k)$ to $(k)$.  The resulting diagram is given explicitly in \S \ref{subsec-cases2and3}.  

\begin{remark}\label{rmk-twistedProj}
In this paper we also consider twisted projectors $w(\pr{n})$, which can be described as follows.  Let $w\in S_n$ be given.  Recall that $\pr{n}$ is a complex whose chain bimodules are shifted copies of $\Q[\xx,\yy]/I$.  The differential on $\pr{n}$ can be represented by matrices whose entries are polynomials $q(\xx,\yy)$.  We define $w(\pr{n})$ to be the complex given by the same underlying bimodule, but replacing $q(\xx,\yy)$ with $q(\xx,w(\yy))$.
\end{remark}

\subsection{Outline of the paper}
Section \S \ref{sec-projector} begins with a recollection of some relevant facts regarding Soergel bimodules.  We then give a set of axioms which characterize $\pr{n}$, and we construct the dual complex $\pr{n}^\vee$ as an explicit injective resolution. In \S \ref{sec-prelims} we show that $\pr{n}^\vee$ is a homotopy colimit of Rouquier complexes associated to powers of the full-twist. In \S \ref{sec-computations} we compute $\HHH(F\otimes \pr{n}^\vee)$ for arbitrary Rouquier complexes $F$, and prove Conjecture \ref{conj-GR}. In \S \ref{sec-combinatorics} we prove some combinatorial results which are likely known to experts but difficult to find in the literature.  Their proofs are elementary in any case, so we include them here.

\subsection{Acknowledgments}  The authors would like to thank E.~Gorsky and J.~Rasmussen for sharing their beautiful work.  Our \S 4 was especially influenced by their Conjecture \ref{conj-GR}. The authors would like to also thank J.~Allman for many helpful conversations.

\section{A categorified Young symmetrizer}
\label{sec-projector}
In this section we introduce the category $\SBim_n$ of Soergel bimodules and describe an idempotent complex $\pr{n}\in\K^-(\SBim_n)$ which categorifies a Young symmetrizer.

\begin{notation}\label{notation-Kom}
If $\mathcal{A}$ is an additive category, we let $\K(\mathcal{A})$ denote the homotopy category of complexes over $\mathcal{A}$ (with differentials of degree $+1$).  We use superscripts $b,+,-$ to denote the full subcategories of complexes which are bounded, respectively bounded from above, respectively bounded from below.  The homological grading shift of complexes is denoted by $\ip{k}$, so that $(C\ip{k})_i=C_{i-k}$.  By convention, the differential on $C\ip{k}$ is $(-1)^k d_C$.
\end{notation}

\subsection{The Soergel category}
\label{subsec-soergel}
Here we set up some notation which will be used throughout this paper.  Fix once and for all an integer $n\geq 1$, and put $W:=S_n$.  Let $w_0\in W$ denote the longest word (that is, $w_0(i)=n+1-i$) and $\ell=\frac{1}{2}n(n-1)$ its length.  Embed $S_{n-1}$ into $S_n$ in the standard way, as permutations of $\{1,\ldots,n\}$ which fix $n$.  Denote the longest word of $S_{n-1}\subset S_n$ by $w_1$.  Its length is $\ell-n+1$.

Set $R:=\Q[x_1,\ldots,x_n]$.  Regard $R$ as a graded ring, via $\deg(x_i)=2$ for all $i$.  Note that $W$ acts on $R$ by permuting variables.  For a simple transposition $s\in W$, let $R^s\subset R$ denote the subalgebra consisting of polynomials $f$ with $s(f)=f$.  Define a graded $(R,R)$-bimodule $B_s:=R\otimes_{R^s} R(-1)$, where $(k)$ denotes the functor which shifts grading up by $k$.  That is $(M(k))_i=M_{i-k}$.  Let $\SBim_n$ denote the smallest full subcategory of all graded $(R,R)$-bimodules containing $B_s$ and closed under direct sum, direct summands, grading shift, and tensor product $\otimes_R$.    Objects of $\SBim_n$ are called \emph{Soergel bimodules}.  We will only briefly recall some relevant facts about Soergel bimodules, and we refer the reader to \cite{EKh10} for more details.

The isomorphism classes of indecomposable Soergel bimodules are indexed by $w\in S_n$.  In this paper we are interested only in a few special cases of these bimodules, which we now describe.   For each $I\subset \{1,\ldots,n-1\}$, let $W_I\subset W$ denote the corresponding parabolic subgroup, that is, the subgroup generated by the simple transpositions $(i,i+1)$ with $i\in I$.  This subgroup has a unique element of maximum length, which we denote by $w_I\subset W_I$.  Let $R^I\subset R$ denote the subalgebra consisting of polynomials such that $w(f)=f$ for all $w\in W_I$.  It is well known that the indecomposable Soergel bimodule $B_{w_I}$ admits the following description:
\begin{equation}\label{eq-longestBimodule}
B_{w_I} = R\otimes_{R^I} R(-\ell(w_I))
\end{equation}
The following bimodules are the most important special cases:
\begin{enumerate}
\item $B_s = R\otimes_{R^s} R(-1)$.  This corresponds to the case $I = \{i\}$.  The longest word $w_I\subset W_I$ is the simple transposition $s=(i,i+1)$ which has length 1.
\item $B_{w_0} = R\otimes_{R^W}R(-\ell)$, where $R^W\subset R$ is the subalgebra of symmetric polynomials.  This corresponds to the case $I = \{1,\ldots,n-1\}$, so that $W_I=W$ with longest word $w_I=w_0$.
\item $B_{w_1}=R\otimes_{R^{S_{n-1}}} R(-\ell+n-1)$.  This corresponds to the case $I=\{1,\ldots,n-2\}$, so that $W_I = S_{n-1}\subset S_n = W$.
\end{enumerate}

We will typically denote tensor product over $R$ simply by juxtaposition: $AB:=A\otimes_R B$.  Tensor product over $\Q$ will be denoted by $\sqcup$:
\begin{definition}\label{def-sqcup}
Let $\sqcup:\SBim_i\times \SBim_j\rightarrow \SBim_{i+j}$ denote the bilinear functor $M\sqcup N=M\otimes_\Q N$. Suppose $R_i =\Q[x_1,\ldots,x_i]$ and $R_j = \Q[x_1,\ldots,x_j].$ Note that if $M$ is a graded $(R_i,R_i)$-bimodule  and $N$ is a graded $(R_j,R_j)$-bimodule, then $M\otimes_\Q N$ is naturally a graded $(R_{i+j},R_{i+j})$-bimodule, since $R_i\otimes_\Q R_j\cong R_{i+j}$.  Thus, our definition of $\sqcup$ makes sense.
\end{definition}
For any $M\in \SBim_k$ with $1\leq k\leq n$, we have $M\sqcup R_{n-k} \cong M[x_{k+1},\ldots,x_n]$.  The notational conventions ensure that $B_w\sqcup R_{n-k} \cong B_w \in \SBim_n$ for all $w\in S_k\subset S_n$.  Thus, we often abuse notation and regard an object of $\SBim_k$ as an object of $\SBim_n$ when there is no possibility of confusion.

\begin{notation}\label{abuse-induction}
For any $M\in \SBim_k$, we often will denote $M\sqcup R_{n-k}\in \SBim_n$ simply by $M$.
\end{notation}

The bimodule $B_{w_0}$ tends to absorb Soergel bimodules:
\begin{proposition}\label{prop-BwOnB0}
$B_wB_{w_0}$ and $B_{w_0}B_w$ are direct sums of shifted copies of $B_{w_0}$, for each $w\in S_n$.  In the special case where $w$ is the longest word of $S_k\subset S_n$, we have
\[
B_w B_{w_0} \cong [k]! B_{w_0} \cong B_{w_0}B_w
\]
\end{proposition}
Here $[k]!=[k][k-1]\cdots [2][1]$, and $[j]=\frac{q^j-q^{-j}}{q-q\inv}$.  We are employing the convention that for any $f(q)\in \N[q,q\inv]$, $f(q)M$ denotes the corresponding direct sum of shifted copies of $M$:

$$\left(\sum_{a \in \Z} c_aq^a \right) \cdot M := \bigoplus_{a\in\Z} M(a)^{\oplus c_a}. $$
\begin{proof}
This is proven in \cite{EW13} and is also a corollary of Theorem \ref{thm-frobBases}. 
\end{proof}

We conclude with a notational convention that we will use to shorten many expressions later:
\begin{notation}\label{notation-XandY}
We have independent sets of variables $\xx:=\{x_1,\ldots,x_n\}$ and $\yy=\{y_1,\ldots,y_n\}$.  We regard an $(R,R)$-bimodule as a $\Q[\xx,\yy]$-module, where $x_k$ acts by left multiplication by $x_k$, and $y_k$ acts via right multiplication by $x_k$.  For $1\leq k\leq n$, let $I_k\subset \Q[\xx,\yy]$ denote the ideal generated by elements
\[
e_i(\yy)-e_i(\xx)  \ \ \ (1\leq i\leq k) \ \ \ \ \ \ \ \ \ \text{together with} \ \ \ \ \ \ \ \ y_i-x_i  \ \ \ (k+1\leq i\leq n)
\]
 where $e_k$ denotes the elementary symmetric function.  Note that $B_{w_0}\cong \Q[\xx,\yy]/I_n$ and $B_{w_1}\cong \Q[\xx,\yy]/I_{n-1}$, up to grading shifts.  Taking the quotient by $y_n-x_n$ defines a canonical map $\Q[\xx,\yy]/I_n\twoheadrightarrow \Q[\xx,\yy]/I_{n-1}$.  Dually, there is a canonical map $\Q[\xx,\yy]/I_{n-1}\hookrightarrow \Q[\xx,\yy]/I_{n}$ which sends $1\mapsto \prod_{i=1}^{n-1}(y_i-x_n)$.  This is a well-defined map of $\Q[\xx,\yy]$-modules by Proposition \ref{prop-canonicalMaps}.
\end{notation}

\subsection{The ideal $\IC\subset \K^-(\SBim_n)$}
\label{subsec-mainLemma}
Recall that $p_{1^n}\in\HC_n$ is the unique multiple of $b_{w_0}$ such that $1-p_{1^n}$ annihilates $b_{w_0}$.  In \S \ref{subsec-axioms} we will give a categorical analogue of this characterization.  But first we study the categorical analogue of being a multiple of $b_{w_0}$.

\begin{definition}\label{def-I}
Let $\IC\subset \K^-(\SBim_n)$ denote the full subcategory consisting of complexes whose chain bimodules are isomorphic to direct sums of $B_{w_0}$ with shifts.  
\end{definition}

\begin{proposition}\label{prop-IisAtwoSidedIdeal}
The subcategory $\mathcal{I}\subset \K^-(\SBim_n)$ is a two-sided tensor ideal.  That is to say, if $I\in \mathcal{I}$ and $A\in\K^-(\SBim_n)$ are arbitrary, then $IA,AI\in\mathcal{I}$.  
\end{proposition}
\begin{proof}
Follows trivially from Proposition \ref{prop-BwOnB0}.
\end{proof}

\begin{definition}\label{def-Iperp}
Let $\IC^\perp\subset \K^-(\SBim_n)$ and ${}^\perp\IC\subset \K^-(\SBim_n)$ denote the full subcategories consisting of complexes $C$ such that $CB_{w_0}\simeq 0$, respectively $B_{w_0}C\simeq 0$.
\end{definition}

\begin{remark}\label{rmk-Iperp}
An easy limiting argument shows that in fact $C\in \IC^\perp$ implies that $CD\simeq 0$ for all $D\in \IC$, and similarly for ${}^\perp\IC$.
\end{remark}

Even though kernels and images don't exist in $\SBim_n$, it still makes sense to talk about the homology of $C\in\K(\SBim_n)$, regarded as a complex of $(R,R)$-bimodules.  Our next goal in this section is to show that if $C\in\K^-(\SBim_n)$ is acyclic (has zero homology), then $C\in \IC^{\perp}\cap{}^\perp\IC$.

\begin{proposition}\label{prop-exactnessOfTensors}
The functors $B_w\otimes (-)$ and $(-)\otimes B_w$ are exact for all $w\in S_n$.  In particular, if $C$ is a complex of Soergel bimodules, then
\[
H(B_w\otimes_R C) \cong B_w\otimes_R H(C) \ \ \ \ \ \ \ \ \ \ \ \ \ H(C\otimes_R B_w)\cong H(C)\otimes_R B_w
\]
for all $w\in S_n$.
\end{proposition}
\begin{proof}
Notice that $B_s$ is free as a right (or left) $R$-module.  Each $B_w$ is isomorphic to a direct summand of a tensor product of $B_s$'s.  Hence each $B_w$ is projective as a right (or left) $R$-module.  It follows that tensoring with $B_w$ is exact.
\end{proof}

\begin{lemma}\label{lemma-endB0}
The Soergel bimodule $B_{w_0}$ is free of rank 1 as a graded module over its graded endomorphism ring $\Endg(B_{w_0})\cong R\otimes_{R^W} R$.
\end{lemma}
\begin{proof}
Obvious.
\end{proof}

\begin{lemma}\label{lemma-acyclicImpliesContractible}
Suppose $C\in\mathcal{I}$ is arbitrary.   Then $C$ is contractible if and only if $C$ is acyclic.
\end{lemma}
\begin{proof}
Certainly $C$ being contractible implies that $C$ is acyclic.  To prove the converse, suppose $C\in\mathcal{I}$ is acyclic.  Let $C_i$ denote the $i$-th chain bimodule.  By definition of $\mathcal{I}$, up to equivalence we may assume that each $C_i$ is a direct sum of shifted copies of of $B_{w_0}$.  Then the $R\otimes R$-action on $C_i$ factors through $R\otimes_{R^W} R =\Endg(B_{w_0})$.  So we may regard $C$ as a complex of $R\otimes_{R^W} R$-modules.  On the other hand, $B_{w_0} = R\otimes_{R^W} R(-\ell)$ is free as a $R\otimes_{R^W}R$-module, hence the acyclic complex $C$ is contractible as a complex of $R\otimes_{R^W}R$-modules by standard arguments.  The null-homotopy $h=(h_i:C_i\rightarrow C_{i-1})_i$ commutes with the $R\otimes_{R^W} R$-action, hence $h$ also commutes with the $R\otimes R$-action.  That is to say, $C$ is contractible as a complex of $(R,R)$-bimodules.
\end{proof}

We have two immediate corollaries.  Recall that a chain map $f:C\rightarrow D$ is said to be a \emph{quasi-isomorphism} if $f$ is an isomorphism in homology.

\begin{corollary}\label{cor-quasiisoIsEquiv}
Two complexes $A,B\in \mathcal{I}$ are homotopy equivalent if and only if they are quasi-isomorphic.
\end{corollary}
\begin{proof}
Suppose $A,B\in \mathcal{I}$, and let $f:A\rightarrow B$ be a chain map.  A well known fact states that $f$ is a homotopy equivalence if and only if the mapping cone $\Cone(f)$ is contractible.  By Lemma \ref{lemma-acyclicImpliesContractible}, this holds if and only if $\Cone(f)$ is acyclic.  Another well known fact about mapping cones states that $\Cone(f)$ is acyclic if and only if $f$ is a quasi-isomorphism.
\end{proof}

\begin{corollary}\label{cor-Iperp}
Let $Z\in\K^-(\SBim_n)$ be a acyclic.  Then $Z\in\mathcal{I}^\perp\cap {}^\perp\mathcal{I}$.
\end{corollary}
\begin{proof}
Let $Z\in\K^-(\SBim_n)$ be acyclic, and let $C\in\mathcal{I}_n$ be arbitrary.  Proposition \ref{prop-exactnessOfTensors} implies that $ZC$ and $CZ$ are acyclic, and then Lemma \ref{lemma-acyclicImpliesContractible} implies that $ZC$ and $CZ$ are contractible.
\end{proof}

\subsection{Axiomatics}
\label{subsec-axioms}

Recall that there is an algebra map $R\otimes_{R^W} R\rightarrow R$.  The main object in this paper is the following free resolution:
\begin{definition}\label{def-P}
Let $\pr{n}\buildrel \e\over \rightarrow R$ denote a resolution of $R$ by free graded $R\otimes_{R^W} R$-modules.  
\end{definition}

\begin{theorem}\label{thm-Pn}
We have
\begin{enumerate}
\item[(P1)] $\pr{n}$ is in $\mathcal{I}$.
\item[(P2)] $\Cone(\pr{n}\buildrel \e\over \rightarrow R)$ is in $\mathcal{I}^\perp\cap {}^\perp\mathcal{I}$.
\end{enumerate}
Furthermore, the pair $(\pr{n},\e)$ is uniquely characterized by these properties up to canonical equivalence.  By this, we mean that if $(Q,\nu)$ is another pair satisfying (P1) and (P2) then there is a unique (up to homotopy) map $\phi:Q\rightarrow \pr{n}$ such that $\nu\circ \phi\simeq \e$; this map is a homotopy equivalence.
\end{theorem}

\begin{proof}
By construction, the chain bimodules of $\pr{n}$ are direct sums of shifted copies of $R\otimes_{R^W}R = B_{w_0}(\ell)$, so that (P1) clearly holds.  Since $\e:\pr{n}\rightarrow R$ is by construction a quasi-isomorphism, $\Cone(\e)$ is acyclic.  Then axiom (P2) holds by Corollary \ref{cor-Iperp}.

Suppose $Q\in \IC$ is some complex and $\nu:Q\rightarrow R$ is some chain map such that $\Cone(\nu)C\simeq 0 \simeq C\Cone(\nu)$ for all $C\in \IC$.  Then $\Cone(\nu)\pr{n}\simeq 0$, since $\pr{n}\in\IC$.  This implies that $\nu\otimes\Id_{\pr{n}}$ is homotopy equivalence $Q\pr{n}\rightarrow \pr{n}$.  A similar argument shows that $\Id_Q\otimes \e$ is a homotopy equivalence $Q\pr{n}\rightarrow Q$.  Thus, $Q\simeq \pr{n}$.  Via Corollary 4.23 in \cite{Hog15} one obtains a canonical equivalence.
\end{proof}

\begin{remark}
This is the categorical analogue of the characterization of $p_{1^n}\in\HC_n$.  The complex $\pr{n}$ is a ``multiple'' of $B_{w_0}$, in the sense that its chain bimodules are sums of $B_{w_0}$ with shifts.  The complex $\Cone(\pr{n}\buildrel\e\over\rightarrow R)$ should be thought of as playing the role of the ``difference'' of $R$ and $\pr{n}$, and axiom (P2) states that this complex annihilates $B_{w_0}$.
\end{remark}

The axioms imply that $\Cone(\e)$ annihilates $\pr n$ from the right and left.  Equivalently, $\e:\pr{n}\rightarrow R$ becomes an equivalence after tensoring with $\pr{n}$ on the right or left.  Thus, $(\Cone(\e),\pr{n})$ forms a pair of complementary idempotents in $\KBC n$, in the sense of \cite{Hog15}.  Many properties of $\pr{n}$ can be deduced immediately from the axioms together with some basic theory of categorical idempotents developed in \cite{Hog15} (also \cite{BD14}).  For instance:

\begin{proposition}\label{prop-PnIsCentral}
The idempotents $\pr n\in \KBC n$ and $\Cone(\e)$ are central.  That is, $A\pr n\simeq \pr n A$ and $A\Cone(\e)\simeq \Cone(\e)A$ for any complex $A\in\KBC n$, and this isomorphism is natural in $A$.
\end{proposition}
\begin{proof}
Follows from the fact that $\IC \subset \K^-(\SBim_n)$ is a two-sided tensor ideal and Theorem 4.10 in \cite{Hog15}.
\end{proof}

\subsection{The case $n=2$ and $n=3$}
\label{subsec-cases2and3}
Note that $\pr 1=R$.  In this subsection we describe $\pr 2$ and $\pr 3$.  We generalize these to an explicit combinatorial description of all of the projectors $\pr n$ in \S \ref{subsec-explicitCx}.

The indecomposable bimodules in $\SBim_2$ are $R$ and $B_s$, where $s\in S_2$ is the nontrivial permutation.  Note that $s$ is the longest word, and has length 1.    The projective resolution $\pr 2 \buildrel\e\over \rightarrow R$ is given by
\begin{equation}\label{eqn-explicitP2}
\begin{diagram}
\cdots&  \rTo{y_1-x_2}& B_s(7) &\rTo{y_2-x_2} & B_s(5) &\rTo{y_1-x_2} & B_s(3) & \rTo{y_2-x_2}& \underline{B_s(1)} & \rDashto{1} & R\\
\end{diagram}
\end{equation}

We emphasize that $\pr 2$ is the complex consisting of terms to the left of the dashed arrow. 

Now we give a construction of $\pr{3}$. Let $p_{ij} = y_i - x_j$, then $\pr{3}$ is the total complex of the following perturbed double complex supported in a single quadrant above the dashed arrows.  In $S_3$ there are two simple transpositions, $s,t\in S_3$, where $s$ swaps 1 and 2, and $t$ swaps 2 and 3.  The longest word is $w_0=sts$, and its length is 3.

\begin{equation}\label{eqn-explicitP3}
\begin{minipage}{122mm}
\begin{center}
\includegraphics[width=122mm]{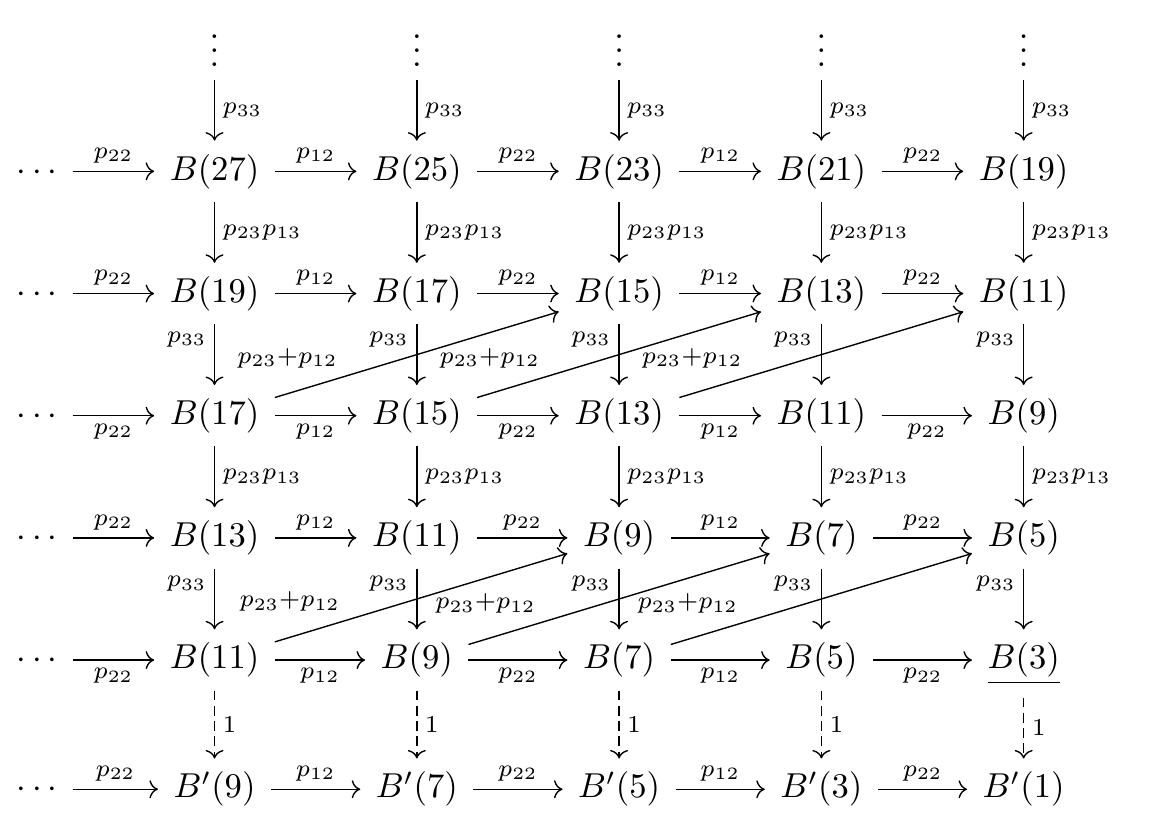}
\end{center}
\end{minipage}
\end{equation}
where $B=B_{sts}$ is the bimodule associated to the longest word in $S_3$, and $B'=B_s \in \KBC 3$.

\begin{remark}\label{rmk-signs}
For now, ignore the part of the diagram below the dashed lines.  To form the total complex $\pr 3$, first order the rows by $0,1,2,\ldots $ starting with the row just above the dashed lines.  Then to each horizontal arrow in an odd numbered row place a minus sign.  Also number the columns from right to left, starting at zero.  The total complex has $-k$-th chain bimodule equal to the direct sum, over $i\in\{0,1,\ldots,k\}$, of terms on the $i$-th row and the $k-i$-th column.  The differential is the (signed) sum of arrows indicated by the diagram.
\end{remark}

Now, the part of the diagram below the dashed arrows is simply $\pr{2}$, regarded now as an object of $\KBC 3$.  The dashed arrows visibly define a chain map $\pr{3}\rightarrow \pr{2}$.  Note the following properties of $\pr{3}$:

\begin{itemize}
	\item In the horizontal direction, the complex is $2$-periodic and each row looks like $\pr{2}$ except replacing all copies of $B_s$ with $B_{sts}$. The extra diagonal arrow is needed to compensate for the fact $p_{22}p_{21} \neq 0$ in $B_{sts}$.
	\item In the vertical direction, the complex is $2$-periodic.  It will be shown that each column of the augmented complex is acyclic, hence the dashed arrows represent a quasi-isomorphism $\pr{3}\rightarrow \pr{2}$.
\end{itemize}

As mentioned in the introduction, we will construct $\pr n$ as an explicit complex in which $B_{w_0}$ appears with graded multiplicity
\begin{equation}\label{multiplicity}
\text{multiplicity of $B_{w_0}$ in $\pr{n}$} = q^{\ell}\prod_{i=1}^n \frac{1+t\inv q^2}{1-t^{-2}q^{2i}}
\end{equation}
where $t$ represents homological degree.  The denominator is realized by an action of the bigraded polynomial ring $\Q[u_1,\ldots,u_n]$ in which $\deg(u_i)=t^{2}q^{-2i}$.  This action generalizes the multiple periodicity evident in our diagrams (\ref{eqn-explicitP2}) and (\ref{eqn-explicitP3}) for $\pr 2$ and $\pr 3$.  The ``fundamental domain'' for this action is a certain cube-like complex $B_{w_0}\otimes \Lambda[\theta_2,\ldots,\theta_n]$, which contributes the numerator in the equation (\ref{multiplicity}).  The 2-periodicity of our complex is an instance of a common phenomenon in commutative algebra \cite{Eis80}, and can be encoded with an object called a matrix factorization.  We return to this idea in \S \ref{subsec-matrixFact}.

\subsection{An iterated projective resolution construction of $\pr{n}$}
\label{subsec-iteratedProjRes}
Below, let $R^{(k)} = R^{S_k}\subset R$ ($1\leq k\leq n$) denote the subalgebra of polynomials which are symmetric in the variables $x_1,\ldots,x_k$.  Note that $R^{(n)}=R^W$ and $R^{(1)}=R$.  By (\ref{eq-longestBimodule}), $R\otimes_{R^{(k)}} R$ is equal to $B_{w}(\ell(w))$, where $w$ is the longest word of $S_k\subset S_n$.   By definition, $\pr{k}$ is a resolution of $R$ by free graded $R\otimes_{R^{(k)}} R$-modules.

\begin{definition}\label{def-Jk}
Let $\mathcal{J}_k\subset \KBC n$ denote the full subcategory consisting of complexes whose chain bimodules are direct sums of shifted copies of $R\otimes_{R^{(k)}} R$. \end{definition}

\begin{definition}\label{def-Jres}
For any complex $C\in \KBC n$, a \emph{$\JC_k$-resolution} of $C$ will mean a complex $D\in \JC_k$ and a quasi-isomorphism $\phi:D \rightarrow C$.  If $C\in \mathcal{J}_{k}$, then $C$ can be regarded as a complex of graded $R\otimes_{R^{(m)}} R$-modules, for any $m\geq k$.  In this case, a $\JC_m$-resolution of of $C$ is the same as a resolution of $C$ by \emph{free} graded $R\otimes_{R^{(m)}} R$-modules.
\end{definition}

Our construction of $\pr{n}$ is inductive, and is based on the following principle:
\begin{enumerate}
\item Assume that we have constructed a $\JC_{k-1}$-resolution $\e_{k-1}:\pr{k-1} \rightarrow R$.
\item Find an explicit $\JC_k$-resolution $C\rightarrow R\otimes_{R^{(k-1)}} R$.
\item General theory of projective resolutions states that we may replace each copy of $R\otimes_{R^{(k-1)}} R$ in $\pr{k-1}$ by $C$ up to quasi-isomorphism.  The result is a $\JC_k$-resolution $\e_k:\pr{k}\rightarrow  R$.
\end{enumerate}

In this paper we are interested in explicit formulae, so will spend some time on (2) above.  We first address the case when $k$ = 2. We have the following complex
\begin{equation} \label{eqn-exactSeqB1}
\begin{diagram}
C_1 & = & R(4) &\rTo{y_1-x_2} & B_s(3) & \rTo{y_2-x_2} & B_s(1) &\rTo{1} & R\\
\end{diagram}.
\end{equation}
As we will see, this complex is acyclic.  By ``stringing together'' infinitely many copies of this complex, we recover the expression (\ref{eqn-explicitP2}) for $\pr{2}$.  More precisely, there exists a map $\gamma :C_1 \to C_1(4)\langle -3 \rangle$ defined as follows:

\begin{equation} \label{eqn-gamma2strand}
\begin{diagram}[small]
&&&&&&R(4) &\rTo & B_s(3) & \rTo & B_s(1) &\rTo & R\\
&&&&&&\dTo^{\Id}\\
R(8) &\rTo  & B_s(7) & \rTo & B_s(5) &\rTo & R(4)\\
\end{diagram}
\end{equation}

By a straightforward short exact sequence argument, $\Cone(\gamma)$ is acyclic. $\Cone(\gamma)$ has a contractible direct summand $$R(4) \xrightarrow{1} R(4).$$ Therefore, by Gaussian elimination, $\Cone(\gamma)$ is homotopy equivalent to the complex $C_2$ pictured below:
\begin{equation*} 
\begin{diagram}[small]
R(8) &\rTo{y_1-x_2} & B_s(7) & \rTo{y_2-x_2} & B_s(5) &\rTo{y_1-x_2} & B_s(3) & \rTo{y_2-x_2} & B_s(1) &\rTo{1} & R\\
\end{diagram}
\end{equation*}

Inductively we can define $C_{\ell+1} = \Cone(\gamma_\ell: C_\ell \to C_1(4\ell)\langle -1-2\ell \rangle)$, where $\gamma_\ell$ is defined in a way similar to $\gamma$.  That is, $C_\ell$ has a copy of $R$ in homological degree $-1 - 2\ell$, and the map $\gamma_\ell$ sends this copy of $R$ by the identity map onto the copy of $R$ in $C_1(4\ell)\langle -1-2\ell \rangle$ in homological degree $-1 - 2\ell$. This gives a contractible summand of the form $$R(4\ell) \xrightarrow{1} R(4\ell).$$ in $\Cone(\gamma_{\ell})$. There is a natural map $\pi_{\ell}:C_{\ell} \twoheadrightarrow C_{\ell-1}$ coming from the mapping cone construction. The colimit of the inverse system $\{C_\ell,\pi_\ell\}$ is an acyclic chain complex since homology commutes with taking colimits.  Moreover, after stripping off the contractible summands in each step, we have a complex which has a single copy of $B_s$ in every homological degree except in homological degree 0 where we have a single copy of $R$. Therefore, $C_\infty = \colim(\{C_\ell,\pi_\ell\})$ is an $R\otimes_{R^s} R$-free resolution of $R$; in fact we recover exactly the expression of $\pr 2$ in diagram (\ref{eqn-explicitP2}). 

We now describe in a similar fashion construct a $\JC_n$-resolution of $B_{w_1}$, where $w_1\subset S_{n-1}\subset S_n$ is the longest word.  We first ask: how does the acyclic complex (\ref{eqn-exactSeqB1}) generalize to arbitrary $n$? 

\begin{proposition}\label{prop-canonicalMaps}
There is an acyclic complex of the form
\begin{equation}\label{eq-Z}
Z \ \ := \ \ 
\Big(
\begin{diagram}
B_{w_1}(2n) & \rTo & B_{w_0}(n+1) &\rTo^{1\otimes x_n - x_n\otimes 1} & B_{w_0}(n-1) &\rTo & B_{w_1}
\end{diagram}
\Big)
\end{equation}
Here, the first and third components of the differential are the bimodule maps
\begin{enumerate}
\item $B_{w_0}(n-1)\rightarrow B_{w_1}$ sending $1\otimes 1\mapsto 1\otimes 1$,
\item $B_{w_1}(n-1)\rightarrow B_{w_0}$ sending $1\otimes 1\mapsto \prod_{i=1}^{n-1}(x_n\otimes 1 - 1\otimes x_i)$.
\end{enumerate}
\end{proposition}
Recall that $1\otimes 1\in B_{w_0}$ has degree $\frac{1}{2}n(n-1)$ and $1\otimes 1\in B_{w_1}$ has degree $\frac{1}{2}(n-1)(n-2)$.  The difference of these is $n-1$, which explains the grading shifts above.
\begin{proof}[Proof of Proposition \ref{prop-canonicalMaps}]
It is known that $R^I$ is a Frobenius algebra over $R^W$, and a pair of dual bases is provided by $\{(-1)^{k}x_n^k\}_{k=0}^{n-1}$ and $\{e_{n-1-k}(x_1,\ldots,x_{n-1})\}_{k=0}^{n-1}$, where $e_k$ denotes the elementary symmetric polynomial.  For a proof of this fact, see Theorem \ref{thm-frobBases} in \S \ref{subsec-frobExtension}. 

In general, suppose $k$ is a commutative ring and $S$ is a Frobenius $k$-algebra with trace $\partial:S\rightarrow k$.  If $\{a_i\}$ and $\{b_i\}$ are dual bases, in the sense that $\partial(a_ib_j)=\d_{ij}$, then there is an $(S,S)$-bilinear map $\Delta:S\rightarrow S\otimes_k S$ sending $1\mapsto \sum_{i=1}^r a_i\otimes b_i$.  This map is independent of the given pair of dual bases.  In our present situation, this gives rise to a map $\Delta:R^I\rightarrow R^I\otimes_{R^W} R^I$ defined by $\Delta(1) = \sum_{i=0}^{n-1}(-1)^i x_n^i\otimes e_{n-i-1}(x_1,\ldots,x_{n-1})$.  

Below, we will adopt Notation \ref{notation-XandY}.  We identify $R^I$ with a subspace of $\Q[\xx]$, and we identify $R^I\otimes_{R^W} R^I$ with a subquotient of $\Q[\xx,\yy]$.   Let us abbreviate $e_k:=e_k(y_1,\ldots,y_{n-1})$ and $\bar e_k:= e_k(y_1,\ldots,y_n)$.  The above discussion gives us a map $R^I\rightarrow R^I\otimes_{R^W} R^I$ such that
\[
\Delta(1)=\sum_{i=0}^{n-1}(-1)^ix_n^i e_{n-1-i} = \prod_{i=1}^{n-1}(y_i-x_n)
\]
Consider the following sequence of $(R^I,R^I)$-bimodules and bimodule maps
\begin{equation}\label{eq-4termAcyclicDude}
\begin{diagram}
R^I &\rTo^{d_0} & R^I\otimes_{R^W}R^I & \rTo^{d_1} & R^I\otimes_{R^W}R^I & \rTo^{d_2} & R^I
\end{diagram}
\end{equation}
where $d_0=\Delta$, $d_1$ is multiplication by $y_n-x_n$, and $d_2$ is the map which sends $1\mapsto 1$ (hence sends $f(\yy)-f(\xx)\mapsto 0$ for all polynomials $f$).   We claim that the above is a chain complex.   Clearly $d_2\circ d_1=0$.  The identity $d_1\circ d_0=0$ follows since $\prod_{i=1}^n(x_n-y_i)$ is symmetric in $y_1,\ldots,y_n$, hence equals $\prod_{i=1}^n(x_n-x_i)=0$ in $R^I\otimes_{R^W} R^I$.

Now we show that the complex is acyclic.  By Theorem \ref{thm-frobBases}, we have $R^I\cong \bigoplus_{i=0}^{n-1}(-x_n)^i R_W$, hence 
\[
R^I\otimes_{R^W} R^I \cong \bigoplus_{i=0}^{n-1} (-x_n)^i\otimes R^I
\]
as right $R^I$-modules.  Thus, we may identify $R^I\otimes_{R^W}R^I$ with $(R^I)^{\oplus n-1}$ as right $R^I$-modules.  With respect to this decomposition, the components $d_0,d_1,d_2$ are represented by the matrices
\[
\matrix{e_{n-1}\\ e_{n-2} \\ \vdots \\ e_1 \\ 1},
\ \ \ \ \ \ \ 
\matrix{
y_n		&		0 			&		 \cdots 	& 		0 			& 		0  	&	 -\bar e_{n-1} \\
1 			&		y_n		&		\cdots 		&		0			&		0 		& -\bar e_{n-2} \\
0			&		1			&		\cdots		&		0			&		0 		& -\bar e_{n-3} \\
\vdots	&		\vdots &      \ddots 		&		\vdots	& 	\vdots & \vdots \\
0			&		0		   &			\cdots &			1 		& y_n 		& -\bar e_2 \\
0			&		0		   &			\cdots &			0		& 	1 			&  y_n-\bar e_1},
\ \ \ \ \ \ \ \ \  
\matrix{1& -y_n & \cdots & (-1)^{n-1}y_n^{n-1}}
\]
The form of the last column of the middle matrix comes from expressing $x_n^n$ in terms of $\bar e_i$ and $x_n^j$ for $1\leq j<n$, using the relation \ref{eq-xnRelation}.  After Gaussian elimination, this complex is clearly acyclic, in fact contractible as a complex of left $R^I$-modules.

Now, to obtain the complex (\ref{eq-Z}) from the complex (\ref{eq-4termAcyclicDude}), apply the functor $R\otimes_{R^I}(-)\otimes_{R^I}R$.  Note that $R$ is free as a right or left $R^I$-module, hence this functor sends acyclic to complexes to acyclic complexes.  This completes the proof.
\end{proof}

We will refer to the maps $\phi:B_{w_1}(n-1)\rightarrow B_{w_0}$ and $\psi:B_{w_0}(n-1)\rightarrow B_{w_1}$ as the \emph{canonical maps} in the sequel.  In case $n=2$ we recover the usual canonical maps (the ``dots'' in the Elias-Khovanov diagram category \cite{EKh10}) $B_s(1)\rightarrow R$ and $R(1)\rightarrow B_s$. 

This $Z$ forms the building block of a $\JC_n$-resolution of $B_{w_1}$.  That is, we consider the following semi-infinite chain complex
\[
\begin{diagram}[small]
&& &&&& A & \rTo & B &\rTo^{y_n-x_n} & B & \rTo & \underline{A}\\
&& &&&& & \rdTo^{-1} & &&&&\\
&& A & \rTo & B &\rTo^{y_n-x_n} & B &\rTo & {A}&&&&\\
&& & \rdTo^{-1} & &&&& &&&&\\
\cdots &\rTo^{y_n-x_n} &  B &\rTo & A &&&& &&&&
\end{diagram}
\]
where we are omitting the degree shifts, we have abbreviated $B=B_{w_0}$, $A=B_{w_1}$, the unlabeled arrows are the canonical maps.  This complex can be described as a colimit as in the discussion at the beginning of this section, and is acyclic.  After Gaussian elimination we obtain:

\begin{proposition}\label{prop-resolutions}
Let $w_0\in S_n$ and $w_1\in S_{n-1}\subset S_n$ denote the longest words.  Set $B:=B_{w_0}$ and $A:=B_{w_1}$.  Let 
$D_1$ denote the semi-infinite complex
\[
\begin{diagram}[small]
 \cdots & \rTo^{a_{nn}} & B(3n+1) & \rTo^{z} & B(3n-1) & \rTo^{a_{nn}} & B(n+1) & \rTo^{z}& B(n-1) & \rTo  &  {A} &\rTo&0
\end{diagram},
\]
and let $D_2$ denote the semi-infinite complex
\[
\begin{diagram}[small]
0 &\rTo & A & \rTo & B(1-n) & \rTo^{z} & B(-1-n) & \rTo^{a_{nn}}& B(1-3n) & \rTo^{z} &  \cdots
\end{diagram}.
\]
where the unlabeled arrows are the canonical maps, $a_{nn}=\prod_{i=1}^{n-1}(y_i-x_n)$, and $z=y_n-x_n$. Then $D_1$ and $D_2$ are acyclic.  In particular $D_1$ represents a $\JC_n$-resolution $C\buildrel \e\over\rightarrow B_{w_1}$.  \qed
\end{proposition}

The polynomial $a_{nn}(\xx,\yy)$ is a special case of Definition \ref{def-aij}.  We want to now iterate the above to construct a resolution of $R$ by free graded $R\otimes_{R^W}R$-modules, and thus an explicit construction of $\pr{n}$.  We illustrate this with an example:

\begin{example}[Construction of $\pr{3}$]
	Let $n = 3$. The discussion at beginning of this section gives $\pr{2}$ as a $\JC_2$-resolution of $R$.  We now replace each copy of $B_s$ with its $\JC_3$-resolution given by Proposition \ref{prop-resolutions}. After doing this we get the following sequence of complexes:
	
	\begin{equation*}
\begin{minipage}{122mm}
\begin{center}
\includegraphics[width=122mm]{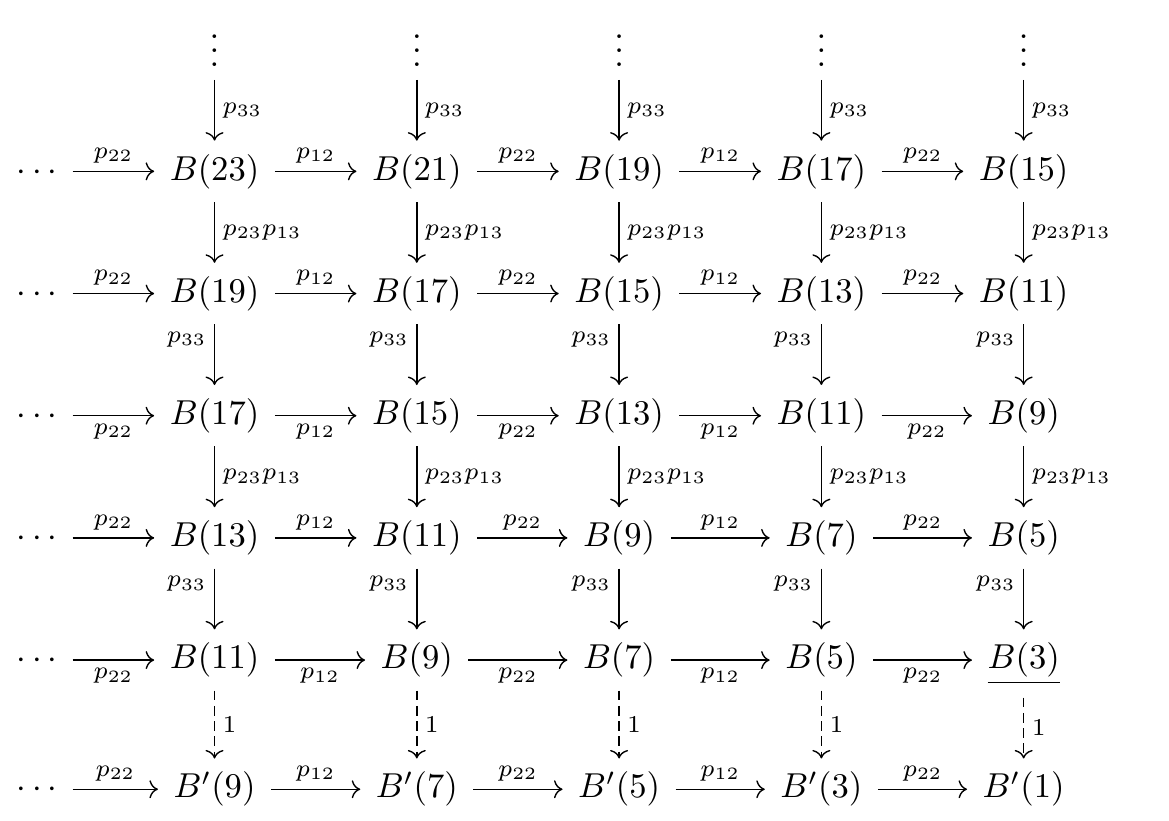}
\end{center}
\end{minipage}
\end{equation*}
This a not a double complex, since the horizontal differential does not square to zero.  Thus, in forming the ``total complex'' (precisely: \emph{convolution}) it is necessary to add extra components to the differential, which in this case are certain diagonal arrows.  Doing so gives us the expression for $\pr 3$ in (\ref{eqn-explicitP3}).
\end{example}

The existence of the diagonal arrows in the above example is implied by a standard fact from homological algebra:

\begin{lemma}[Functorial projective resolutions] \label{lem-projResCom}
	Let $X = (X, d)$ be a chain complex of objects in an abelian category $\mathcal{C}$. Suppose $P_i$ is a projective resolution of $X_i$, then we can replace each $X_i$ with $P_i$ to form a sequence
	$$ P(X) = \cdots \xrightarrow{} P_i \xrightarrow{\tilde{d}_i} P_{i-1} \xrightarrow{} \cdots$$
where $\tilde{d}_i$ is the unique lift (up to homotopy) of $d_i$. Let $f_i: P_i \to X_i$ be the augmentation map of the projective resolution $P_i$ of $X_i$. There is a unique convolution $\Tot(P(X))$ and  the maps $f_i$ determine a quasi-isomorphsm $f : \Tot(P(X)) \to X$.  The assignment $X\mapsto \Tot(P(X))$ defines a functor $\K^-(\mathcal{C})\rightarrow \K^-(\mathcal{C})$.
\end{lemma}

	This is proven in many standard homological algebra texts such as Weibel \cite{Wei94}.
\begin{corollary}\label{cor-functorialResolution}
For each $X\in \mathcal{J}_k$ there is a $\JC_{k+1}$-resolution of $X$, canonical up to homotopy, denoted $P(X)\buildrel \e_X\over \rightarrow X$.  The assignment $X\mapsto P(X)$ defines a functor $\JC_k\rightarrow \JC_{k+1}$.\qed
\end{corollary}

This result together with Proposition \ref{prop-resolutions} gives our construction of $\pr n$.
\begin{construction}\label{const-iteratedProjRes}
Consider $\pr{n-1}$ as a $\JC_{n-1}$-resolution of $R$, and let $\Phi:\JC_{n-1}\rightarrow \JC_n$ denote the functor from Corollary \ref{cor-functorialResolution}.  Then $\pr{n} = \Phi(\pr{n-1})$.  More precisely, each chain module $X_i$ of $\pr{n-1}$ is a direct sum of shifted copies of $B_{w_1}$.  Replace each $X_i$ with the corresponding direct sum of shifted copies of
\[
\begin{diagram}[small]
 \cdots & \rTo^{y_n-x_n} & B_{w_0}(3n-1) & \rTo^{a_{nn}} & B_{w_0}(n+1) & \rTo^{y_n-x_n}& B_{w_0}(n-1)
\end{diagram},
\]
The result is a sequence $$P = \cdots \xrightarrow{} P_i \xrightarrow{} P_{i-1} \xrightarrow{} \cdots.$$ Then $\Tot(P) \simeq \pr{n}$ as objects of $\KBC n$.
\end{construction}

It is clear that in the $\pr{n}$ just constructed, the bimodule $B_{w_n}$ appears with graded multiplicity $\prod_{i=2}^n\frac{q\inv+t\inv q}{q^{-i}-t^{-2}q^{i}}$.  Substituting $t=-1$ yields $\frac{1}{[n]!}$ (see also \S \ref{subsec-introComb}).

\subsection{A combinatorial construction of the projector}
\label{subsec-explicitCx}
Now that we have a rough idea of the structure of $\pr{n}$, we are ready to describe $\pr{n}$ explicitly.  The diagonal arrows (see Diagram \ref{eqn-explicitP3}) are described by an explicit family of polynomials $a_{ij}(\xx,\yy)$, which we define first.

\begin{definition}\label{def-aij}
Let $a_{ij}(\xx,\yy)$ ($2\leq i\leq j$) denote the polynomials defined by
\[
a_{ij}(\xx,\yy):=\sum_{\gamma} \prod_{k=1}^{i-1}(y_{\gamma_k} - x_{\gamma_k+i-k})
\]

where the sum is over decreasing sequences $\gamma = (\gamma_1,\ldots,\gamma_{i-1})$ with $j-1\geq \gamma_1>\cdots >\gamma_{i-1}\geq 1$.  Note that the degree $a_{ij}$ is $i-1$.  By convention, we set $a_{1,j}:=1$ for all $j\geq 1$.  We also find it useful to abbreviate $p_{ij}:=y_i-x_j$.
\end{definition}

\begin{example}\label{ex-a1-4}
For $1< i\leq j\leq 4$, the $a_{ij}$ are given below:
\begin{enumerate}
\item $a_{22} = p_{12}$
\item $a_{23}=p_{12}+p_{23}$.
\item $a_{24}=p_{12}+p_{23}+p_{34}$.
\item $a_{33}=p_{13}p_{23}$.
\item $a_{34}=p_{13}p_{23}+p_{13}p_{34}+p_{24}p_{34}$
\item $a_{44}=p_{14}p_{24}p_{34}$.
\end{enumerate}
\end{example}
Note that $a_{22}$, $a_{23}$, and $a_{33}$ appear prominently in our diagrams which define $\pr{2}$ and $\pr{3}$ (Diagrams (\ref{eqn-explicitP2}) and (\ref{eqn-explicitP3})).  Also, one of the differentials in the complexes from Proposition \ref{prop-resolutions} appears as a special case of the $a_{ij}$:

\begin{example}\label{ex-ann}
We have $a_{jj}=p_{1j}p_{2j}\cdots p_{j-1,j}$.
\end{example}

Recall the ideal $I_n\in \Q[x_1,\ldots,x_n,y_1,\ldots,y_n]$ generated by $e_k(\xx)-e_k(\yy)$ ($1\leq k\leq n$).  In \S \ref{subsec-explicitRels} we prove the following:

\begin{proposition}\label{prop-aijRelation}
The $a_{ij}(\xx,\yy)$ satisfy $\sum_{j=k}^n a_{i,j}(\xx,\yy)(x_j-y_j) =0$ modulo $I_n$.  In particular $\sum_{j=k}^n a_{i,j}(\xx,\yy)(x_j-y_j) $ acts by zero on $B_{w_0}$.
\end{proposition}

Using these polynomials we will define a complex (the equation $d^2=0$ will follow from Proposition \ref{prop-aijRelation}), and then we will show that this complex is homotopy equivalent to $\pr{n}^\vee$, which is dual to $\pr{n}$ (see Definition \ref{def-duality} below).  The reason that $\pr{n}^\vee$ is easier to describe is that it is naturally a unital algebra in $\K^+(\SBim_n)$,  whereas $\pr{n}$ is a coalgebra in $\K^-(\SBim_n)$.  The language of differential bigraded (dg) algebras and modules will help make the discussion brief, and will make manifest the periodicity in our complexes (evident in our expressions for $\pr 2$ and $\pr 3$).

\begin{definition}\label{def-dga}
Let $u_k$ denote an even formal indeterminate of bidegree $t^{2}q^{-2k}$, and let $\theta_k$ denote an odd formal indeterminate of bidegree $t^{1}q^{-2}$.  Let $A_n$ denote the super-polynomial algebra
\[
A_n:=\Q[\xx,\yy,u_1,\ldots,u_n,\theta_1,\ldots,\theta_n]
\]
with $\Q[\xx,\yy,u_1,\ldots,u_n]$-linear differential determined by $d_A(1)=0$ and
\begin{eqnarray}
d_A(\theta_k)&=& \sum_{i=1}^k a_{i,k}(\xx,\yy) u_i \label{eq-thetaDiff},
\end{eqnarray}
for each $k=1,\ldots,n$, together with the graded Leibniz rule $d(ab)=d(a)b+(-1)^{|a|}ad(b)$.
\end{definition}
Note that the differential $d_A$ is defined independently of $n$.  That is, the obvious inclusion $A_{n}\rightarrow A_{n+1}$  is a chain map for all $n\geq 1$.

\begin{definition}\label{def-dgMod}
Recall the ideal $I_n\in \Q[\xx,\yy]$ from Notation \ref{notation-XandY}.  We define the $A_n$-module $M_n:=A_n/I_nA_n$.  Denote the image of $a\in A_n$ under the quotient map $A_n \twoheadrightarrow M_n$ by $\overline{a}$.  Shift $M_n$ so that $\overline{1}\in M_n$ lies in degree $q^{-2\ell}$ where $\ell=\binom{n}{2}$ as usual.   Define a differential $d_M$ on $M_n$ by the rules
\begin{enumerate}
\item $d_M(\overline{1})=\sum_{k=1}^n(y_k-x_k)\overline{\theta}_k$
\item $d_M(a m) = d_A(a)m+ (-1)^{|a|}a d_M(m)$
\end{enumerate}
for all $a\in A$, $m\in M$.  Here, $|a|$ denotes the homological degree of a homogeneous element $a\in A$.
\end{definition}

We check that the above differential is well-defined, hence makes $M_n$ into a dg $A_n$-module.  Since $M$ is generated by $\overline{1}\in M_n$ as a left $A_n$-module, $d_M$ is unique (if it exists).  On the other hand, the formula
\[
d_M(\overline{a})=d(a\overline{1})={d_A(a)}\overline{1} +(-1)^{|a|}\sum_{k=1}^n(y_k-x_k)a\overline{\theta}_k
\]
shows that $d_M$ is a well-defined $\Q[\xx,\yy,u_1,\ldots,u_n]$-linear map $M_n\rightarrow M_n$.  We claim that $d_M^2=0$.  Given the Leibniz rule (2), it suffices to show that $d_M^2(\bar{1})=0$. 
 Applying ${d_M}$ to the the equation $d_M(\bar{1}) = \sum_{i=1}^n (y_i-x_i)\theta_i\overline{1}$ and using the Leibniz rule gives
\begin{eqnarray*}
d_M^2(1)
&=&\sum_{i=1}^n (y_i-x_i) d_A(\theta_i)\overline{1} - \sum_{i=1}^n (y_i-x_i) \theta_i d_M(\overline{1})\\
&=&\sum_{i=1}^n\sum_{j=1}^i (y_i-x_i) a_{ji}(\xx,\yy)u_j \overline{1} - \sum_{i=1}^n\sum_{j=1}^n (y_i-x_i)(y_j-x_j)\theta_i \theta_j\\
&=& \sum_{j=1}^n u_j \sum_{i=j}^n (y_i-x_i) a_{ji}(\xx,\yy) + 0\\
&=&0
\end{eqnarray*}
In the third line we used that the odd variables $\theta_i$ anti-commute and square to zero.  In the fourth line we appeal to Proposition \ref{prop-aijRelation}.

\begin{remark}\label{rmk-cubeCx}
We may visualize the complex $M_n$ in the following way.  Set $B_n:=\Q[\xx,\yy]/I_n$.  Then $M_n\cong B_n[u_1,\ldots,u_n,\theta_1,\ldots,\theta_n]$ is a direct sum of infinitely many copies of $B_n$, indexed by the positive orthant $\Z_{\geq 0}^{n}$.  Multiplication by $u_k$ acts by translating a distance of two units along the $k$-th axis.  The summands corresponding to the ``unit cube'' $\{0,1\}^{n}$  contribute a copy of the exterior algebra $B_n[\theta_1,\ldots,\theta_n]\subset A_n$, which then generates all of $A_n$ under the action of the $u_k$.  The differential $d_M$ makes the the unit cube $B_n[\theta_1,\ldots,\theta_n]$ into a Koszul complex associated to the sequence $p_{11},p_{22},\ldots,p_{nn}\in B_n$, and the $\Q[u_1,\ldots,u_n]$-equivariance ensures that the translates of the unit cube all have the same internal differential. 
	
There are components between different copies of the unit cube which are encoded by the algebra differential $d_A$; these are defined in terms of the polynomials $a_{ij}(\xx,\yy)$.
\end{remark}

\begin{remark}
There are reduced versions of $A_n$ and $M_n$ defined as follows.  Set $ \tilde{\theta}_k:=\theta_k-\theta_1$.  Then set $ \tilde{A}_n = \Q[\xx,\yy,u_2,\ldots,u_n, \tilde{\theta}_2,\ldots, \tilde{\theta}_n]$ with differential determined by $d(\tilde{\theta}_k)=\sum_{i=2}^k u_i a_{ik}(\xx,\yy)$.  The reader may check that $A_n \cong B\otimes_\Q  \tilde{A}_n$, where $B=\Q[u_1,\theta_1]$ with differential $d(\theta_1)=u_1$.  Clearly $B\simeq \Q$, so that $A_n\simeq  \tilde{A_n}$.  Similarly, the reduced version of $M_n$ is $\tilde{M_n} =  \tilde{A_n}/ I_n \tilde{A}_n$ with differential determined by $d_M(\bar 1) = \sum_{i=2}p_{ii} \tilde{\theta}_i$.  The reader may check that $M_n\cong B\otimes  \tilde{M}_n\simeq \tilde{M}_n$.  
\end{remark}

By $\Q[\xx,\yy,u_1,\ldots,u_n]$-equivariance, $M_n$ is determined by the $d_M$ applied to a product of $\theta_k$'s.  It is instructive to work these out in the cases $n=2,3$; The $n=3$ case is done in \S \ref{subsec-introComb} of the introduction.

After extending scalars, we will regard $A_k$ and $M_k$ as complexes of $\Q[\xx,\yy]$-modules, for all $1\leq k\leq n$.  For example $A_1 = \Q[\xx,\yy,u_1,\theta_1]$ with differential $d_A(\theta_1)=u_1$, and $M_1= \Q[\xx,u_1,\theta_1]$ with differential $d_M(\bar{1})=p_{11}\bar{\theta_1}$.  Clearly $A_1\simeq \Q[\xx,\yy]$ and $M_1\simeq \Q[\xx]$.  We wish to relate $M_{n-1}$ with $M_{n}$.  The following is obvious, since the formula for $d_A(\theta_k)$ does not involve $n$:

\begin{proposition}\label{prop-mapOfdga}
There is a unique map of dg $\Q[\xx,\yy]$-algebras $A_{n-1}\rightarrow A_n$ sending $\theta_k\rightarrow \theta_k$ for $1\leq k\leq n-1$.\qed
\end{proposition}
Thus, by restriction, we regard $M_n$ as a dg $A_{n-1}$-module.  We want to construct a canonical map $M_{n-1}\rightarrow M_n$ of dg $A_{n-1}$-modules.  First, note that from Proposition \ref{prop-canonicalMaps}, there is a well defined $\Q[\xx,\yy]$-module map
\begin{equation}\label{eq-mapOfCoeffs}
\Q[\xx,\yy]/I_{n-1} \rightarrow \Q[\xx,\yy]/I_n \ \ \ \ \ \ \ \ \ \ \ \ \ \overline{1}\mapsto \prod_{i=1}^{n-1}(y_i-x_n)\overline{1}
\end{equation}

\begin{proposition}\label{prop-mapOfdgm}
There is a unique map $\phi:M_{n-1}\rightarrow M_n$ of dg $A_{n-1}$-modules sending $\overline{1}\mapsto \prod_{i=1}^{n-1}(y_i-x_n)\overline{1}\in M_n$.  
\end{proposition}
We remark that $\phi$ preserves the bidegrees.  This is the reason for choosing $\deg_q(\overline{1}_n)=n-n^2$, where $\overline{1}_n = \bar{1} \in M_n$.  
\begin{proof}
Uniqueness is clear, since $M_{n-1}$ is generated by $\overline{1}$ as an $A_{n-1}$-module.  The map (\ref{eq-mapOfCoeffs}) extends to a map $\phi:M_{n-1}\rightarrow M_n$ of $A_{n-1}$-modules.  We must check that $\phi$ commutes with the differentials $d_{M_n}$ and $d_{M_{n-1}}$.  The computation reduces easily to the computation $\phi(d_{M_{n-1}}(\overline{1}))=d_{M_n}(\phi(\overline{1}))$, which can be checked directly:
\begin{eqnarray*}
\phi(d_{M_{n-1}}(\overline{1}))
&=&\phi\bigg(\sum_{i=1}^{n-1} (y_i-x_i)\phi(\theta_i\overline{1})\bigg)\\
&=&\sum_{i=1}^{n-1} \theta_i (y_i-x_i)\prod_{j=1}^{n-1}(y_j-x_n)\overline{1}\\
&=&\sum_{i=1}^{n} \theta_i (y_i-x_i)\prod_{j=1}^{n-1}(y_j-x_n)\overline{1}\\
&=&d_{M_n}(\phi(\overline{1}))
\end{eqnarray*}
In the third line we used the fact that $\prod_{i=1}^n (y_i-x_n)$ is zero in $\Q[\xx,\yy]/I_n$, by Proposition \ref{prop-explicitRelations} (see also Example \ref{example-theRelation}).
\end{proof}

\begin{lemma}\label{lemma-filtration}
The map $\phi:M_{n-1}\rightarrow M_n$ is a quasi-isomorphism.
\end{lemma}
\begin{proof}
We prove this by constructing a filtration on $\Cone(\phi)$ whose subquotients are the acyclic complex $D$ from Proposition \ref{prop-resolutions}.

First, note that there is an obvious isomorphism
\[
A_n\cong \Q[\xx,\yy,u_n,\theta_n]\otimes_\Q \Q[u_1,\ldots,u_{n-1},\theta_1,\ldots,\theta_{n-1}]
\]
as algebras (ignoring the differentials, and also using the usual sign rule for the tensor product of two superalgebras).   Taking the quotient by the ideal generated by $e_k(\yy)-e_k(\xx)$ gives an isomorphism
\begin{equation}\label{eq-Mdecomp}
M_n\cong (\Q[\xx,\yy,u_n,\theta_n]/I_n) \otimes_\Q \Q[u_1,\ldots,u_{n-1},\theta_1,\ldots,\theta_{n-1}]
\end{equation}
Similarly, we have
\[
M_{n-1}\cong (\Q[\xx,\yy]/I_{n-1})\otimes \Q[u_1,\ldots,u_{n-1},\theta_1,\ldots,\theta_{n-1}]
\]
We filter each of the above complexes by homological degree on the second tensor factor.  We remark that any homogeneous element $c_1\otimes c_2$ of any of the above complexes satisfies $\deg_h(c_1\otimes c_2) = \deg_h(c_1)+\deg_h(c_2)$.  Since $\deg_h(c_1)\geq 0$, this means that the filtration degree of any homogeneous element is bounded above by its homological degree.  We will use this fact at the end of this proof.

The differentials respect the filtrations, as does the map $\phi$.  With respect to this filtration we have
\begin{enumerate}
\item $d_A(1\otimes a_2)=0+(\text{higher})$ since we are filtering by homological degree on the second factor.
\item $d_A(\theta_n\otimes 1) = a_{nn}(\xx,\yy)u_n\otimes 1 + (\text{higher})$.
\item $d_M(\overline{1}\otimes 1) = (y_n-x_n)\overline{\theta}_n\otimes 1 + (\text{higher})$.
\end{enumerate}
It is now straightforward to check that the induced filtration on $\Cone(\phi)$ has associated graded
\[
\gr \Cone(\phi) \cong D\otimes_\Q \Q[u_1,\ldots,u_{n-1},\theta_1,\ldots,\theta_{n-1}]
\]
where $D$ is the complex
\[
\begin{diagram}
 \Q[\xx,\yy]/I_{n-1} & \rTo^{a_{nn}} & \Q[\xx,\yy]/I_n & \rTo^{y_m-x_m} & \Q[\xx,\yy]/I_n & \rTo^{a_{nn}}& \Q[\xx,\yy]/I_n & \rTo^{y_m-x_m} & \cdots
\end{diagram}
\]
Since $D$ is acyclic, it follows that $\gr \Cone(\phi)$ is acyclic.  The following argument shows that $\Cone(\phi)$ is acyclic, hence $\phi$ is a quasi-isomorphism.

We have a filtration $\Cone(\phi)=F^0\supset F^1\supset \cdots$ whose associated graded is acyclic.  This implies that any class $z\in H(F^{k})$ is homologous to a class in $H(F^{k+1})$, as follows from the long exact sequence in homology associated to the short exact sequence
\[
0\rightarrow F^{k+1}\rightarrow F^k\rightarrow F^k/F^{k+1}\rightarrow 0.
\]
In particular, any class $z\in H(\Cone(\phi))$ is homologous to a class in $H(F^k)$ for all $k\geq 0$.  On the other hand,  by the remark above, the filtration degree of any nonzero homogenous element $c\in \Cone(\phi)$ is bounded above by its homological degree.  It follows that any class $z\in H(\Cone(\phi))$ is null-homologous.  This completes the proof that $\Cone(\phi)$ is acyclic.
\end{proof}

Now we show $M_n$ is dual to $\pr{n}$. 
\begin{definition}\label{def-duality}
If $M$ is a graded $(R,R)$-bimodule, set $M^\vee=\Homg_R(M,R)$, where this latter object is the graded space of right $R$-module maps $M\rightarrow R$.  Since $R$ is commutative, this is a graded $(R,R)$-bimodule via $(\a \cdot f\cdot \b)(m) = f(\b m \a)$ for all $\a,\b\in R$, $m\in M$, $f\in M^\vee$.
\end{definition}
It is well known that if $w_I\in S_n$ is the longest word of a parabolic subgroup $W_I\subset S_n$ (see \S \ref{subsec-soergel}), then $B_{w_I}^\vee \cong B_{w_I}$.  In particular, $B_s^\vee\cong B_s$, so the duality functor sends Soergel bimodules to Soergel bimodules.  Extending to complexes gives a contravariant functor $\K^-(\SBim_n)\leftrightarrow \K^+(\SBim_n)$.
\begin{proposition}\label{prop-Pdual}
Suppose $C\in \K^+(\SBim_n)$ is a complex whose chain bimodules are direct sums of copies of $B_{w_0}$ with grading shifts, and suppose there is a quasi-isomorphism $\nu:R\rightarrow C$.  Then $C^\vee\simeq \pr{n}$, with structure map $\nu^\vee: C^\vee\rightarrow R$.
\end{proposition}
\begin{proof}
Since Soergel bimodules are free as right $R$-modules, $(-)^\vee$ is an exact functor.  In particular, $H(C)\cong 0$ if and only if $H(C^\vee)\cong 0$, and $\nu:R\rightarrow C$ is a quasi-isomorphism if and only if $\nu^\vee$ is a quasi-isomorphism.  From the remarks preceding the proposition, the chain bimodules of $C^\vee$ are direct sums of shifted copies of $B_{w_0}^\vee\cong B_{w_0}$, so that $C^\vee\in \IC$.  Thus, $(C^\vee,\nu^\vee)$ satisfies the axioms from Theorem \ref{thm-Pn}, so $C^\vee \simeq \pr{n}$ by uniqueness.
\end{proof}

\begin{theorem}\label{thm-MisP}
$M_n\simeq \pr{n}^\vee$.  
\end{theorem}
\begin{proof}
Composing the quasi-isomorphisms from Lemma \ref{lemma-filtration} gives a quasi isomorphism $\e: \Q[\xx]=M_1\rightarrow M_n$.  Note that $\e$ is a map of dg $A_1=\Q[\xx,\yy]$-modules, that is, complexes of $(R,R)$-bimodules.  Proposition \ref{prop-Pdual} now says that $M_n^\vee\simeq \pr{n}$, which completes the proof.
\end{proof}

\section{The projector as an infinite full twist}
\label{sec-prelims}

We introduce the Rouquier complexes, which define an action of the braid group on the homotopy category of complexes of Soergel bimodules.  Upon restricting to $\IC\in \K^-(\SBim_n)$ from Definition \ref{def-I}, this action factors through the symmetric group.  We then use this fact to show that the projector $\pr{n}$ is a homotopy colimit of Rouquier complexes associated to powers of the full twist. 

\subsection{Action of the symmetric group}
\label{subsec-symGrp}
In this section we construct a categorical action of the symmetric group on the ideal $\mathcal{I}\subset \K^-(\SBim_n)$ from Definition \ref{def-I}.

\begin{definition}\label{def-standardBimod}
Given $w\in W$, let $R_w$ denote the $(R,R)$-bimodule which equals $R$ as a left module, and whose right $R$-action is twisted by $w$.  That is, $f\cdot h \cdot g = f h w(g)$ for all $f,g\in R$, $h\in R_w$.  The bimodules $R_w$ are not Soergel bimodules, but are often referred to as \emph{standard bimodules}.
\end{definition}
Note the important difference between superscipts and subscripts: $R^s$ is the algebra of $s$-invariant polynomials, whereas $R_s$ is the $(R,R)$-bimodule whose right $R$-action is twisted by $s$.

\begin{proposition}\label{prop-B0absorbsStandards}
We have $B_{w_0}R_w\cong B_{w_0}$ for every $w\in S_n$.  If $f=f(\xx,\yy)$, thought of as an endomorphism of $B_{w_0}$, then the following square comments:
\[
\begin{diagram}
B_{w_0}R_w & \rTo^{f(\xx,\yy)\otimes \Id_{R_w}} & B_{w_0}R_w \\
\uTo^{\cong} && \uTo^{\cong}\\
B_{w_0} & \rTo^{f(\xx,w\inv(\yy))} & B_{w_0}
\end{diagram}
\]
There is a similar description of the functor $R_w\otimes (-)$.
\end{proposition}
\begin{proof}
Ignore grading shifts throughout this proof, and identify $B_{w_0}$ with the quotient $\Q[\xx,\yy]/I_n$.  By abuse of notation, we will denote elements of $B_{w_0}$ by $g(\xx,\yy)$.  Note that $B_{w_0}R_w = \Q[\xx,\yy]/I_n$, but with $\Q[\xx,\yy]$-module structure $g(\xx,\yy)\cdot 1 = g(\xx,w(\yy))$.  With these identifications, it is trivial to check that $g(\xx,\yy)\mapsto g(\xx,w(\yy))$ defines an isomorphism $B_{w_0}\rightarrow B_{w_0}R_w$.  The statement regarding commutativity of the square is equally straightforward.
\end{proof}

Proposition \ref{prop-B0absorbsStandards} implies that the functor $(-)\otimes R_w$ gives a right action of the symmetric group on $\mathcal{I}$.   We will give a convenient description of this action.

\begin{definition}
Recall that any $C\in \IC$ is homotopy equivalent to a complex whose chain bimodules are direct sums of shifted copies of $B_{w_0}$.  Note that each component of the differential can be regarded a matrix with entries in $\Q[\xx,\yy]$.  Let $w(C)$ denote the complex with the same underlying bimodule, but with differential obtained from $d_C$ by applying $f(\xx,\yy)\mapsto f(\xx,w(\yy))$.
\end{definition}
The following is clear:

\begin{proposition}
For each $C\in \IC$ and each $w\in S_n$, we have $w(C)\cong CR_{w\inv}$.
\end{proposition}
\begin{definition}
We will refer to the complexes $w(\pr{n})$ as \emph{twisted projectors}.
\end{definition}

\begin{remark}\label{rmk-twistPn}
The explicit formulae for $\pr{n}$ also apply to $w(\pr{n})$, by applying the transformation $f(\xx,\yy)\mapsto f(\xx,w(\yy))$ everywhere.
\end{remark}

\subsection{Action of the braid group}
\label{subsec-rouquier}
In this section we construct an action of the braid group on $\K^-(\SBim_n)$, originally due to Rouquier (see also \cite{Rou06}).  We show that when restricted to $\IC\subset \K^-(\SBim_n)$, this action factors through the symmetric group.  This fact will be used in our argument that $\pr{n}$ is a limit of Rouquier complexes.

The following is standard in the theory of Soergel bimodules \cite{EKh10}:
\begin{proposition}\label{prop-SES}
We have short exact sequences
\begin{equation}\label{eq-SES1}
0\rightarrow  R_s(1) \rightarrow B_s \rightarrow R(-1)\rightarrow 0
\end{equation}
\begin{equation}\label{eq-SES2}
0\rightarrow  R(1) \rightarrow B_s \rightarrow R_s(-1)\rightarrow 0
\end{equation}\qed
\end{proposition}

The complexes $F_s^\pm$ defined below should be thought of as certain kinds of deformations of $R_s$.  They are quasi-isomorphic but not chain homotopy equivalent to $R_s$.

\begin{definition}\label{def-rouquier}
Define the following complexes of Soergel bimodules:
\begin{eqnarray*}
F_s &:=& \Big(\begin{diagram}[small]\underline{B_s} & \rTo^{1} & {R}(-1) \end{diagram}\Big)\\
F_s\inv &:=& \Big(\begin{diagram}[small] {R}(1) & \rTo^{x_1-y_2} & \underline{B_s} \end{diagram}\Big)
\end{eqnarray*}
where we have underlined the degree zero chain bimodules, and we are using the conventions of Notation \ref{notation-XandY}.  Suppose $b$ is a braid word, i.e.~a word in $(s,\pm)$, where $s$ is a simple transposition.  We will denote the corresponding product $F_{s_1}^\pm \cdots F_{s_r}^\pm$ simply by $F_{b}$.  We call each $F_b$ a \emph{Rouquier complex}.  If $\underline{w}$ is a word in the simple transpositions we will denote by $F_{\underline{w}}^+, F_{\underline{w}}^-$ the Rouquier complex associated to the positive (resp.~negative) braid lift of $\underline{w}$.
\end{definition}
If two braid words $b$ and $b'$ represent the same braids, then $F_b\simeq F_{b'}$, and the homotopy equivalence is canonical up to homotopy \cite{EliasKrasner2010}

\begin{definition}\label{def-rouquierQuasiisos}
Let $\phi_s: R_s(1)\rightarrow F_s$ denote the chain map induced by the first map in the short exact sequence (\ref{eq-SES1}).  For any positive braid word $\b$ with length $r$, let $\phi_\b:R_s(r)\rightarrow F(\b)$ denote the corresponding tensor product of maps $\phi_s$.  Similarly, let $\psi_s:F_s\inv\rightarrow R_s(-1)$ denote the chain map induced by the second map in the short exact sequence (\ref{eq-SES2}), and let $\psi_\b:F(\b)\rightarrow R_s(-r)$ denote the coresponding tensor product of maps $\psi_s$, where $\b$ is a negative braid with length $r$.
\end{definition}

\begin{proposition}\label{prop-RouquierQuasiIso}
The maps $\phi_{\b}:R_w(r)\rightarrow F(\b)$ and $\psi_{\b}:F(\b)\rightarrow R_w(-r)$ are quasi-isomorphisms.
\end{proposition}
\begin{proof}[Sketch of proof]
First, note that exactness of the sequences (\ref{eq-SES1}) and (\ref{eq-SES2}) implies that $\phi_s$ and $\psi_s$ are quasi-isomorphisms.  A standard homological algebra fact states that the derived tensor product of two quasi-isomorphisms is a quasi-isomorphism.  Soergel bimodules are free as left or right $R$-modules, so for complexes of Soergel bimodules, derived tensor product coincides with ordinary tensor product.

\end{proof}
 As a corollary we obtain the following:
\begin{theorem}\label{thm-catSymAction}
The categorical braid group action on $\mathcal{I}$ factors through the symmetric group.  More precisely, let $\ww=(s_1^{\nu_1},\ldots, s_r^{\nu_r})$ denote a braid word ($\nu_i=\pm 1$), let $w\in S_n$ denote the permutation represented by $\ww$, and let $e=\nu_1+\cdots +\nu_r$ denote the braid exponent.  Then $C\mapsto F_{\ww} C \simeq R_wC(e)$ and $C F_{\ww}\simeq C R_w(e)$ for all $C\in \IC$.    This is an isomorphism of functors $\IC\rightarrow \IC$.
\end{theorem}
\begin{proof}
It suffices to show that tensoring with $R_s(\pm 1)$ and $F_s^{\pm 1}$ give isomorphic functors $\mathcal{I}\rightarrow \mathcal{I}$.  From Proposition \ref{prop-RouquierQuasiIso}, $\phi_s:R_s(1)\rightarrow F_s$ and $\psi_s:F_s\inv\rightarrow R(-1)$ are quasi-isomorphisms.  By Corollary \ref{cor-quasiisoIsEquiv}, these quasi-isomorphisms become homotopy equivalences after tensoring with any object of $\mathcal{I}$.  Thus, $\phi_s$ and $\psi_s$ define natural isomorphisms of endofunctors of $\mathcal{I}$.
\end{proof}

 \subsection{Constructing the projector via infinite full-twists}
 \label{subsec-infiniteTwists}
 In this section we show that $\pr n ^\vee $ is a homotopy limit of Rouquier complexes.  As a result, $\HHH(\pr{n}^\vee)$ is a limit of Khovanov-Rozansky homologies of torus links.  We will compute the homology $\HHH(w(\pr{n}^\vee))$ for all $w$ in subsequent sections.
 
\begin{remark}
We focus on the dual projector $\pr n ^\vee$ for the rest of this section. This discussion can easily be altered to work for $\pr n$ as well. We make this choice because we focus on $\pr n ^\vee$ in the sequel.
\end{remark}

 \begin{definition}\label{def-FT}
 	Fix an integer $n\geq 1$ and let $\HT_n = F_{w_0}$ denote Rouquier complex associated to the positive lift of the longest word.  This is the half-twist on $n$ strands, which has $\ell = \frac{1}{2}n(n-1)$ crossings.  Let $\FT_n := \HT_n\HT_n$ denote the Rouquier complex associated to the positive full twist on $n$-strands.  Let $\phi:R(2\ell)\rightarrow \FT_n$ denote the quasi-isomorphism from Proposition \ref{prop-RouquierQuasiIso}. The index $n$ is fixed throughout this section, so we will omit the subscript throughout.\end{definition}
 
 The following is a special case of Theorem \ref{thm-catSymAction}:
 \begin{corollary}\label{cor-FTfixedSet}
 	For any $C\in \IC$, the maps $\Id_C\otimes \phi: C(2\ell)\rightarrow C\FT$ and $\phi\otimes \Id_C:C(2\ell)\rightarrow \FT C$ are homotopy equivalences.  In particular $\FT \pr{n}\simeq \pr{n}(2\ell)$.\qed
 \end{corollary}
 
 \begin{example}
 	The full twist on two strands is homotopy equivalent to the complex
 	\[
 	\FT_2\ \  \simeq \ \ \begin{diagram} (\underline{B_s}(1) &\rTo^{y_1-x_2} & B_s(-1) &\rTo & {R(-2)} )\end{diagram}
 	\]
 	The chain map $\phi:R(2)\rightarrow \FT_2$ has mapping cone
 	\[
 	\Cone(\phi) \ \ \simeq \ \ \begin{diagram} (\underline{R(2}) &\rTo^{y_1-x_2} & B_s(1) &\rTo^{y_2-x_2} & {B_s(-1)} & \rTo^{1} & R(-2) )\end{diagram}
 	\]
 	This complex is acylic by Proposition \ref{prop-canonicalMaps}, hence kills $B_s$ from the left and right by Corollary \ref{cor-Iperp}.  It is a useful exercise to prove this directly.
 \end{example}

\begin{lemma}\label{lemma-powersOfHT}
Let $k\geq 1$ be a given integer.  There is a complex $X\simeq \HT$ such that the chain bimodules of $X^{\otimes k}$ in homological degrees $<k$ are direct sums of $B_{w_0}$ with shifts.
\end{lemma}
\begin{proof}
First, note that $\HT$ is supported in non-negative homological degrees, since it is the Rouquier complex associated to a positive braid.  In the case $k=1$, we appeal to the fact (see Theorem 6.9 in \cite{EW14}) that if $F(\b)$ is the Rouquier complex associated to the positive braid lift of a \emph{reduced} expression for $w\in S_n$, then stripping off contractible summands yields the minimal complex $(F(\b))_{\text{min}}\simeq F(\b)$ whose degree zero chain bimodule is $B_w$.  Setting $w=w_0$ proves the $k=1$ case of the proposition.

Let $X=\HT_{\text{min}}$, so that the degree zero chain bimodule of $X$ is $X_0=B_{w_0}$. Then the degree $j$ chain bimodule of $X^{\otimes k}$ is
\[
(X^{\otimes k})_j = \bigoplus_{i_1+\cdots +i_k = j} \ X_{i_1}\cdots X_{i_k}
\]
If $j<k$, then this forces at least one of the indices $i_m$ to satisfy $i_m=0$, hence at least one of the factors above equals $B_{w_0}$.  By Proposition \ref{prop-BwOnB0}, this forces $(C^{\otimes k})_j$ to be a direct sum of $B_{w_0}$ with shifts.
\end{proof}
 
Now basic idea of stablization is this: by the above, the homological degree $<2k$ chain bimodules of $X^{\otimes 2k}\simeq \FT^{\otimes k}$ form an object of $\IC$.  By Corollary \ref{cor-FTfixedSet}, $\FT$ fixes objects of $\IC$ up to shift, so $\FT^{\otimes k+1}$ and $\FT^{\otimes k}$ are homotopy equivalent in homological degrees $<2k$, up to a grading shift.  This is the stablization we seek.  We now make these arguments precise.  We first define our directed system, and the notion of homotopy colimit.

 \begin{definition}\label{def-theSystem}
$f_0=\phi(-2\ell):R\rightarrow \FT(-2\ell)$, where $\phi$ is as in Definition \ref{def-FT}.  For $k \geq 1$ set $f_k:=f_0\otimes \Id^{\otimes k}:\FT^{\otimes k}(-2k\ell)\rightarrow \FT^{\otimes k+1}(-2(k+1)\ell)$.
 \end{definition}
 
 \begin{definition}\label{def-holim}
 	Let $\{A_k, f_k:A_{k}\rightarrow A_{k+1}\}_{k=0}^\infty$ be a direct system of chain complexes.  The \emph{homotopy colimit} of $\{A_k, f_k\}$ is by definition the mapping cone
 	\[
 	\hocolim A_k := \Cone\bigg(\bigoplus_{k=0}^\infty A_k \buildrel \Id-S\over \longrightarrow \bigoplus_{k=0}^\infty A_k\bigg)
 	\]
 	Where $S:\bigoplus_{k=0}^\infty A_k\rightarrow \bigoplus_{k=0}^\infty A_k$ is the chain map with components given by the $f_k$.
 \end{definition}
 
 \begin{theorem}\label{thm-PasInfiniteBraid}
We have $\pr{n}^\vee\simeq \hocolim A^{\otimes k}$, where $\{A^{\otimes k}, f_k\}$ is the directed system from Definition \ref{def-theSystem}.
 \end{theorem} 

The proof utilizes several lemmas.  First, set $P:= \pr n^\vee$.  Let $\eta=\e^\vee:R\rightarrow \pr{n}^\vee$ be the unit map, and set $C:=\Cone(\eta)\ip{1}$.  Note that $P$ and $C$ are supported in non-negative homological degrees.  The reader may anticipate that this fact will be used in an essential way in the proof, since otherwise we might work with $\pr{n}$ instead of $\pr{n}^\vee$.  We will utilize the existence of an exact triangle
 \begin{equation}\label{eq-resOfId}
 P\ip{1}\rightarrow C \rightarrow \one  \rightarrow P
 \end{equation}
 in $\K^+(\SBim_n)$, where $\one = R$ is the trivial bimodule, and $\ip{1}$ denotes the upward shift in homological degree.  Also, the properties of $P$ (Theorem \ref{thm-Pn}) ensure that
 \begin{equation}\label{eq-orthoIdempts}
 PC\simeq CP\simeq 0,
 \end{equation}
 In the language of \cite{Hog15}, this means that $(P,C)$ is a pair of complementary idempotents in $\K^+(\SBim_n)$.  An object $A\in \K^+(\SBim_n)$ is fixed by $P$ if and only if it is annihilated by the complementary idempotent $C$, and vice versa.  In particlar, we have:

 \begin{lemma}\label{lemma-CkillsB0}
 $B_{w_0}C\simeq 0\simeq CB_{w_0}$.
 \end{lemma}
 \begin{proof}
 The proof of Theorem \ref{thm-Pn} (after applying the duality functor) says that $PB_{w_0}\simeq B_{w_0}$, hence $CB_{w_0}\simeq C(PB_{w_0})\cong (CP)B_{w_0}\simeq 0$ by (\ref{eq-orthoIdempts}).  A similar argument shows that $B_{w_0}C\simeq 0$.
 \end{proof}

\begin{lemma}\label{lemma-cauchy}
$\HT^{\otimes k}C$ is homotopy equivalent to a complex which is supported in homological degrees $\geq k$.
\end{lemma} 
\begin{proof}
Follows easily from Lemma \ref{lemma-CkillsB0} and Lemma \ref{lemma-powersOfHT}.
\end{proof}
 
 \begin{lemma}\label{lemma-trivialHocolims}
 	Suppose $\{N_k,f_k:N_k\rightarrow N_{k+1}\}$ is a directed system of complexes, and $N_k$ is homotopy equivalent to a complex which is supported in homological degrees $\leq c_k$, where $c_k\to -\infty$ as $k\to \infty$.  Then $\hocolim N_k$ is contractible.
 \end{lemma}
 \begin{proof}
 	This is a fairly standard argument.  See, for instance, Proposition 3.3 in \cite{Hog15}.
 \end{proof}

\begin{proof}[Proof of Theorem \ref{thm-PasInfiniteBraid}]
We first claim that $\Cone(f_k)$ is homotopy equivalent to a complex which is supported in homological degrees $\geq 2k-1$, hence $f_k$ should be thought of as a homotopy equivalence in degrees $<2k-1$.  First, observe that by Corollary \ref{cor-FTfixedSet}, $\Cone(f_0)P\simeq 0$.  Tensoring the exact triangle (\ref{eq-resOfId}) on the left with $\Cone(f_0)$ gives an exact triangle
 	\[
 	0\rightarrow \Cone(f_0)C\rightarrow \Cone(f_0)\rightarrow 0
 	\]
 A standard fact about mapping cones now implies that $\Cone(f_0)C\simeq \Cone(f_0)$.  For $k\geq 0$, we have $f_k = f_0\otimes \Id^{\otimes k}$, so
 \[
 \Cone(f_k)\cong \Cone(f_0)\otimes \FT^{\otimes k}(-2k\ell) \simeq \Cone(f_0)C\FT^{\otimes k}(-2k\ell)
 \]
$\Cone(f_0)$ is supported in degrees $\geq -1$ (because of the shift in degrees in forming the mapping cone), and $C\FT^{\otimes k}$ is supported in degrees $\geq 2k$ by Lemma \ref{lemma-cauchy}, so $\Cone(f_k)$ is supported in degrees $\geq 2k-1$.  Thus, the directed system $\{\FT^{\otimes k}(-2k\ell), f_k\}$ is \emph{Cauchy} in the sense of Rozansky \cite{Roz10a}, and the homotopy colimit $L:=\hocolim \FT^{\otimes k}(-2k\ell)$ is defined in $\K^+(\SBim_n)$.
 	
 We now claim that $LC\simeq 0$ and $LP\simeq P$.  Suppose we have shown this.  Tensor the exact triangle (\ref{eq-resOfId}) on the left with $L$.  Since $LC\simeq 0$, the result is an exact triangle $LP\ip{1}\rightarrow 0\rightarrow L\rightarrow LP$, which implies that $LP\simeq L$.  On the other hand $LP\simeq P$, so $P\simeq L$.  
 	
 To complete the proof therefore, we must show that $LC\simeq 0$ and $LP\simeq P$.  By definition, $L$ is the homotopy colimit of complexes $\FT^{\otimes k}(-2k\ell)$, and so $LC$ is the homotopy colimit of complexes $\FT^{\otimes k}C(-2k\ell)$, which is contractible by Lemma \ref{lemma-trivialHocolims} and Lemma \ref{lemma-cauchy}.  That is to say, $LC\simeq 0$.
 	
 Finally, we will show that $LP\simeq P$.  There is a map of directed systems
 	\begin{equation}\label{eq-mapOfsystems}
 	\begin{diagram}
 	P &\rTo^{P\phi} & P\FT &\rTo^{P \phi \FT} & P\FT^{\otimes 2} &\rTo^{P\phi \FT^{\otimes 2}} & P\FT^{\otimes 3} &\rTo^{P \phi \FT^{\otimes 3}} & \cdots\\
 	\uTo^{P} && \uTo^{P\phi} && \uTo^{P\phi^{\otimes 2}} && \uTo^{P\phi^{\otimes 3}} &&\\
 	P &\rTo^{P} & P &\rTo^{P} & P &\rTo^{P} & P &\rTo^{P} & \cdots
 	\end{diagram}
 	\end{equation}
 	where are abusing notation by denoting $\Id_X$ by $X$, for any chain complex $X$, and we are omitting all explicit grading shifts.  The homotopy colimit of the bottom row is simply $P$, and the homotopy colimit of the top row is $PL$.  By Corollary \ref{cor-FTfixedSet}, each vertical arrow is homotopy equivalence.  It is a standard fact (and an easy exercise) to deduce from this that the induced map $P\rightarrow PL$ is a homotopy equivalence.  This completes the proof.
 \end{proof}

\section{Computations of stable homology}
\label{sec-computations}
In this section we introduce the functor $\HHH$, whose input is a complex of Soergel bimodules and whose output is a triply graded vector space.   We introduce matrix factorizations as a way of expressing $\pr{n}$ (or $\pr{n}^\vee$) via compact formulae.  Using this we compute $\HHH(w(\pr{n}^\vee))$ for arbitrary permutations $w\in S_n$, proving a conjecture of Gorsky-Rasmussen.   In case $w\in S_n$ is an $n$-cycle, $\HHH(w^m(\pr{n}^\vee))$ is a stable limit of Khovanov-Rozansky homology (or briefly KR homology in the sequel) of the $(n,kn+m)$-torus links as $k\to \infty$.

\subsection{Hochschild cohomology}
\label{subsec-HH}
We refer the reader to \cite{Kh07} and references therein for more details on Hochschild (co)homology and the connection to link homology.  We review only a few essential facts here.

Let $R^e = R\otimes_\Q R$.  We will regard graded $(R,R)$-bimodules as  graded left $R^e$-modules. The \emph{Hochschild cohomology} of $M$ is defined by $$\HH^k(R;M) = \bigoplus_{i,j}\Ext_{R^e}^j(R(i),M).$$  We let $\HH(R;M)$ denote  the bigraded space $\HH(R;M)=\bigoplus_k \HH^k(R;M)$.  When $R$ is understood, we omit it from the notation and write $\HH(M)$.  The degree of a homogeneous element $z\in \HH(M)$ will be written multiplicatively, as $\deg(z)=q^i a^j$, here $a$ indicates the Hochschild degree. By definition  $\HH^0(M)$ is the graded of graded bimodule maps $R\rightarrow M$.  This will be denoted $\Homg(R,M)$.

We may regard $\HH$ as a functor from the category of $(R,R)$-bimodules to the category of $R$-modules, where all objects are understood to be graded in our case.  The $R$-action on $\HH(M)$ is induced by the $R$-action on the first argument of $\Ext_{R^e}(R,M)$.

Since $\HH$ is a linear functor, we can extend to complexes.  If $C$ is a complex of graded $(R,R)$-bimodules,  then $\HH(C)$ is the complex such that $\HH(C)_k=\HH(C_k)$, with differential $d_{\HH(C)}=\HH(d_C)$.  Note that $\HH(C)$ is triply graded; the differential has tridegree $(0,0,1)$.

\begin{definition}\label{def-HHH}
If $C$ is a complex of graded $(R,R)$-bimodules, let $\HH(C)$ denote the complex obtained by applying $\HH$ to the bimodules of $C$.  Set $\HHH(C):=H(\HH(C))$.  This is a triply graded vector space.  By convention, $z\in \HHH_{i,k}^j(C)$ means that $z$ has $q$-degree $i$, and Hochschild degree $j$, and homological degree $k$.  In this case we write $\deg(z)=q^ia^jt^k$.  
\end{definition}

\begin{remark}\label{rmk-HomflyHomology}
It is possible to show (see \cite{Kh07} and \S 3 of \cite{Hog15}) that if $\b$ is an $n$-strand braid and $F(\b)$ is the Rouquier complex associated to $\b$, then $\HHH(F(\b))$ depends only on the closure of $\b$ up to isomorphism of triply graded vector spaces and an overall grading shift.  The grading indeterminacy can be fixed by a renormalization; the resulting link invariant is KR homology.  Thus, $\HHH(\pr{n}^\vee)$ computes the stable limit of KR homologies of the $(n,nk)$-torus links, and similarly for $\HHH(w(\pr{n}^\vee))$ where $w\in S_n$.
\end{remark}

It is a standard property of Hochschild cohomology that $\HH(M\otimes_R N)\cong \HH(N\otimes_R M)$, and this isomorphism is natural in $M$, $N$. Naturality of this isomorphism means we can extend to complexes as well, and we obtain:

\begin{proposition}\label{prop-HHHisTrace}
Suppose $C$ and $D$ are bounded above complexes of $(R,R)$-bimodules.  Then $\HH(C\otimes_R D)\cong \HH(D\otimes_R C)$.\qed
\end{proposition}

We now turn our attention to the computation of $\HHH(w(\pr{n}^\vee))$, for $w\in S_n$.  The following simplifies this task:
\begin{corollary}\label{cor-twistedHHH}
Let $w\in S_n$ be given.  Then $\HHH(w(\pr{n}^\vee))$ depends only on the conjugacy class (cycle type) of $w$, up to isomorphism.
\end{corollary}
\begin{proof}
Let $R_w$ be the standard bimodule, so that $w\inv(\pr{n}^\vee):=\pr{n}^\vee R_{w}$.  Since $\pr{n}^\vee$ is central, we have
\[
\pr{n} R_w \cong \pr{n}^\vee R_{v\inv}R_{vw} \cong  R_{v\inv} \pr{n}^\vee R_{vw}
\]
Applying $\HHH(-)$ gives
\[
\HHH(\pr{n}^\vee R_{w}) \cong \HHH(R_{v\inv}\pr{n}^\vee R_{vw})\cong \HHH(\pr{n^\vee} R_{vw}R_{v\inv}) \cong \HHH(\pr{n}^\vee R_{vwv\inv})
\]
This completes the proof.
\end{proof}

Our next simplification says that $\HHH(C)$ is especially simple whenever $C\in \IC$, that is, if the chain bimodules of $C$ are direct sums of $B_{w_0}$ with shifts.  Recall that $\Endg(B_{w_0}) = R\otimes_{R^W} R\cong \Q[\xx,\yy]/I_n$ (Notation \ref{notation-XandY}), so that we can regard the differential $d_k:C_k\rightarrow C_{k+1}$ as a matrix with entries from $\Q[\xx,\yy]/I_n$.

\begin{proposition}\label{prop-HHIn}
If $C \in \IC$, then  $\HH(C) \simeq \Lambda[\xi_1,\ldots,\xi_n] \otimes_\Q \HH^0(C)$, where the $\xi_i$ are odd variables of degree $q^{-2i}a$.  Furthermore, $\HH^0$ is the functor that sends $B_{w_0}\mapsto R(\ell)$ and $\Endg(B_{w_0})\owns f(\xx,\yy)\mapsto  f(\xx,\xx)$.
\end{proposition}
The rest of this section is concerned with proving this proposition.   As before, let $R = \Q[\xx]$, $R^e = \Q[\xx,\yy]$, $\ell = \frac{1}{2}n(n-1)$.  Write $B:=B_{w_0} = \Q[\xx,\yy]/I_n$, as usual.  Let $K$ denote the Koszul resolution of $R$ by free graded $R^e$-modules.  



Below, $C$ will denote a complex whose chain bimodules are direct sums of $B$ with grading shifts.  We use $C(i,j)$ to denote a grading shift in both $q$ and homological degree: $$(C(i,j))_k = C_{k-j}(i).$$

\begin{lemma}\label{lem-HHBn}
	The Hochschild cohomology of $B$ is isomorphic to a shift of a bigraded superpolynomial ring with odd variables $\xi_i$ and even variables $x_i$: $$\HH(B) \cong \Q[x_1,\ldots,x_n,\xi_1,\ldots,\xi_n](\ell),$$ where $\xi_i$ has degree $q^{-2i}a$ and $x_i$ has degree $q^2$.
\end{lemma}
\begin{proof}
For this proof it is actually much more convenient to work with Hochschild homology $\HH_\bullet$ rather than cohomology $\HH^\bullet$.  Hochschild homology is defined by $\HH_{-k}(M) = \Tor^k_{R^e}(R,M)$.  For later convenience, our grading is such that $\HH_{\bullet}$ is supported in \emph{negative} Hochschild degrees.  It is well known that for our graded polynomial ring $R$,
\begin{equation}\label{eq:HHandHH}
\HH^\bullet(R;M)=q^{-2n}a^n\HH_\bullet(R;M).
\end{equation}  Throughout this proof, we use variables $q$ and $a$ to denote the grading shifts in those degrees.  Now, recall $R^e = R \otimes R \simeq \Q[\xx,\yy]$ and $q^{\ell}B\cong \Q[\xx,\yy]/I_n$ as in Notation \ref{notation-XandY}.  It is convenient to compute $\HH_\bullet(q^\ell B)$ using the Koszul resolution:
	\[
	K = \bigotimes_{i = 1}^n (\Q[\xx,\yy](2i) \xrightarrow{e_i(\xx)-e_i(\yy)} \underline{\Q[\xx,\yy]})
	\]
	where the tensor product is over $\Q[\xx,\yy]$ and the underlined term is in Hochschild degree 0.  Then $\HH_\bullet (q^{\ell}B)$ is computed by identifying the left and right actions on $K$---that is, setting $y_i=x_i$ for $1\leq i\leq n$---and taking homology.  Identifying the left and right actions yields
	\[
	R\otimes_{R^e} K = \bigotimes_{i = 1}^n ( R(2i) \xrightarrow{0} \underline{R}),
	\] 
	 If we let $\xi_i^\vee$ be an odd variable of Hochschild degree $-1$ and $q$-degree $2i$ then we can rewrite this as $$ R\otimes_{R^e} K\cong \Q[\xx,\xi_1^\vee,\ldots,\xi_n^\vee]$$
	with zero differential.  Taking homology gives the computation of $\HH_\bullet(q^{\ell}B)=q^{\ell}\HH_\bullet(B)$.  Given Equation (\ref{eq:HHandHH}), we obtain $\HH^\bullet(B) = q^{-2n-\ell}a^n \Q[\xx,\xi_1^\vee,\ldots,\xi_n^\vee]$.  It is a trivial exercise to check that this is isomorphic to $q^{\ell}\Q[\xx,\xi_1,\ldots,\xi_n]$ as we claimed.  The isomorphism is essentially the Hodge star operator, which for instance sends $\xi_1^\vee\cdots\xi_n^\vee\mapsto 1$ and $1\mapsto \xi_1\cdots \xi_n$.  
\end{proof}
%


\begin{proof}[Proof of Proposition \ref{prop-HHIn}]
Suppose $C$ is a complex whose chain bimodules are direct sums of copies of $B=B_{w_0}$ with shifts.  Lemma \ref{lem-HHBn} tells us what happens to the bimodules uppon application of the functor $\HH$.  Specifically, $\HH(B)\cong \HH^0(B)\otimes \Q[\xi_1,\ldots,\xi_n]$, and $\HH^0(B)\cong R(\ell)$.   We need only study how $\HH$ acts on elements of $\Endg(B)$.  But every endomorphism of $B$ is given by multiplication by some $f(\xx,\yy)\in R^e$.  Since $\HH$ identifies the left and right $R$-actions, the induced endomorphism $\HH(f)$ is simply multiplication by $f(\xx,\xx)$ (recall that $\HH(M)$ is naturally an $R$-module).  This completes the proof.
\end{proof}

\begin{corollary} \label{cor-HHPn}
We have $\HH(w(\pr{n}^\vee)) \cong \Lambda[\xi_1,\ldots,\xi_n] \otimes_\Q \HH_0(w(\pr{n}^\vee))$, where $\deg(\xi_i) \\ =q^{-2i}a$.\qed
\end{corollary}

\subsection{Some computations of $\HHH(w(\pr n^\vee))$}
\label{subsec-computations}

In this section we explicitly compute the Hochschild homology of twisted projectors, $\HHH(w(\pr{n}^\vee))$, in a few special cases to give motivation for the results of \S \ref{subsec-EndPn} and \S \ref{subsec-matrixFact}.   We let $(i,j)=(i)\ip{j}$ denote a shift in $q$-degree and homological degree.

\begin{example}\label{ex-P2Unknot}[$\HHH^0(\pr{2}^\vee)$]
By Corollary \ref{cor-HHPn}, the computation of $\HHH(\pr{2}^\vee)$ reduces to the computation of the homology of $\HH^0(\pr{2}^\vee)$, which we now describe.  We replace every copy of $B= B_s$ with $R(1)$ and replace each $y_i$ with $x_i$ in the dual complex of \ref{eqn-explicitP2} to obtain the complex that follows.
	\begin{equation}\label{eqn-HH0P2}
	\begin{diagram}
	\underline{R}&\rTo{0}&\underline{R}(-2)&\rTo{x_2-x_1}&\underline{R}(-4)&\rTo{0}&\underline{R}(-6)&\rTo{x_2-x_1}&\cdots,\\
	\end{diagram}
	\end{equation}
where the underlined term is in homological degree 0.  If we let $\alpha = x_2-x_1$ then we can write $R = \Q[\alpha,x_1]$ and $R/(\alpha) \simeq \Q[x_1]$. (\ref{eqn-HH0P2}) splits into $R$ and the short chain complexes
	$$R(-2-4i,1+2i) \xrightarrow{\alpha} R(-4-4i,2+2i)$$
	We can reduce each of the above chain complexes, obtaining
	
	$$H(\HH^0(\pr{2})) \cong R(-2) \oplus \bigoplus_{i \geq 0} R/(\alpha)(-4-4i,2+2i).$$
	If we let $u_2$ denote the periodicity map described in $\S \ref{subsec-explicitCx}$ then it is easy to see that $$H(\HH^0(\pr{2})) \cong \Q[x_1,\alpha,u_2]/(\alpha u_2)$$
\end{example}

\begin{example}\label{eqn-P2ClrdUnknot} [$\HHH^0(s(\pr{2}^\vee))$]
Let $s\in S_2$ be the nontrivial permutation.  We know, by Remark \ref{rmk-twistPn}, that $\HH_0(s(\pr{2}^\vee))$ is the complex
	
	\begin{equation}\label{eqn-HH0sP2}
	\begin{diagram}
	\underline{R}&\rTo{x_2-x_1}&\underline{R}(-2)&\rTo{0}&\underline{R}(-4)&\rTo{x_2-x_1}&\underline{R}(-8)&\rTo{0}&\cdots,\\
	\end{diagram}
	\end{equation}
	
	Therefore, $\HH^0(s(\pr{2}^\vee))$ is a direct sum of complexes of the form $$R(-4i,2i) \xrightarrow{x_2-x_1} R(-2-4i,2i+1),$$ for all $i \geq 0$. The homology of each of these short complexes is $R/(\alpha)(-2-4i,2i+1)$. Therefore $\HHH(s(\pr{2})) \cong \Lambda[\xi_1,\xi_2] \otimes_\Q \Q[x_1,u_2](-2,1)$.
\end{example}

\begin{example}\label{ex-nondual}[$\HHH^0(\pr{2})$]
Let $S=\Q[x_1,x_2,u_2]$, with zero differential.  Consider the dg $S$-module $M_\infty=\Q[x_1,x_2,u_2^{\pm},\theta_2]$ with differential $d(\theta_2)=(x_1-x_2)u_2$.  Define the following sub and quotient dg $S$-modules:
\[
M_+:=\Q[x_1,x_2,u_2,\theta_2]\subset M_\infty \ \ \ \ \ \ \ \ \ \ \ \ M_-:=M_\infty / u_2M_+.
\]
Note that $M_+$ and $M_-$ are generated by monomials in which $u_2$ appears with non-negative, respectively non-positive, exponent.

In Example \ref{ex-P2Unknot} we saw that $\HHH^0(\pr{2}^\vee)\cong H(M_+)$ up to a grading shift.  A similar computation shows that $\HHH^0(\pr{2})\cong H(M_-)$ up to a grading shift.  Note that $\HHH^0(\pr{2})$ is \emph{not} finitely generated over the endomorphism ring $\Q[x_1,x_2,u_2]/(x_1u_2=x_2u_2) \cong \End(\pr{2})$. A similar comparison can be done for $\HHH(\pr{n})$ and $\HHH(\pr{n}^\vee)$.
\end{example}

\begin{example}\label{eqn-P3ClrdUnknot} [$\HHH^0(w(\pr{3}^\vee))$, $w=(1,2,3)$]
Recall the complex for $\pr{3}$ from \ref{eqn-explicitP3}. We form the dual complex and twist this complex by the permutation $w:=(1,2,3)$ (see Remark \ref{rmk-twistPn}) by rewriting each $p_{ij}$ with $p_{w(i)j}$. Signs can be placed in a similar manner to what was done in Remark \ref{rmk-signs}.
	\begin{equation}\label{eqn-explicitwP3}
\begin{minipage}{122mm}
\begin{center}
\includegraphics[width=122mm]{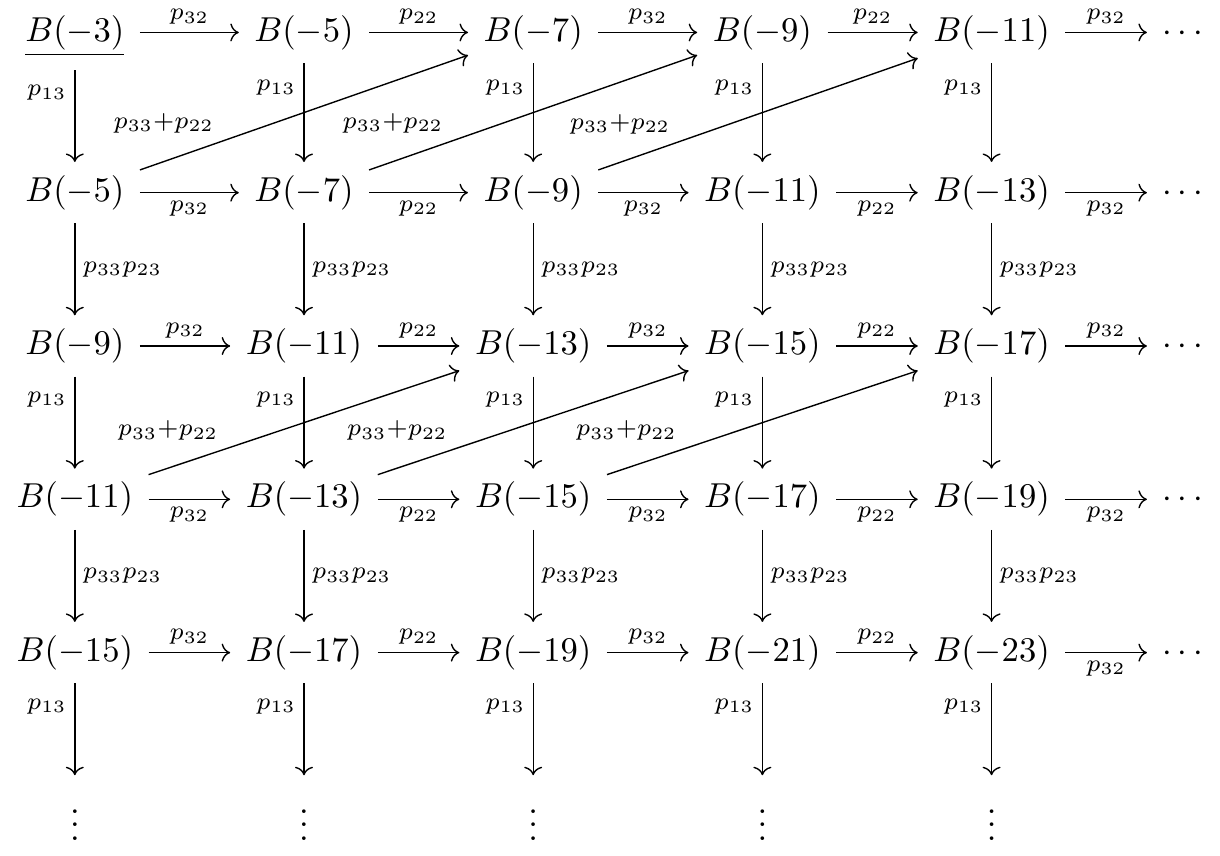}
\end{center}
\end{minipage}
\end{equation}

	By Corollary \ref{cor-HHPn}, when we apply $\HH$ to $w(\pr{3}^\vee)$, we can factor out a copy of the exterior algebra $\Lambda[\xi_1,\xi_2,\xi_3]$ and we are left with the complex $\HH^0(w(\pr{3}^\vee))$.  Recall that $\HH^0$ sends $B\mapsto R(3)$ and $f(\xx,\yy)\mapsto f(\xx,\xx)$.  Set $q_{ij} = \HH(p_{ij})= x_i-x_j$.  Note that $q_{ii}=0$, so $\HH^0(w(\pr{3}^\vee))$ is the total complex of
	
	\begin{equation}\label{eqn-HH0wP3}
\begin{minipage}{122mm}
\begin{center}
\includegraphics[width=122mm]{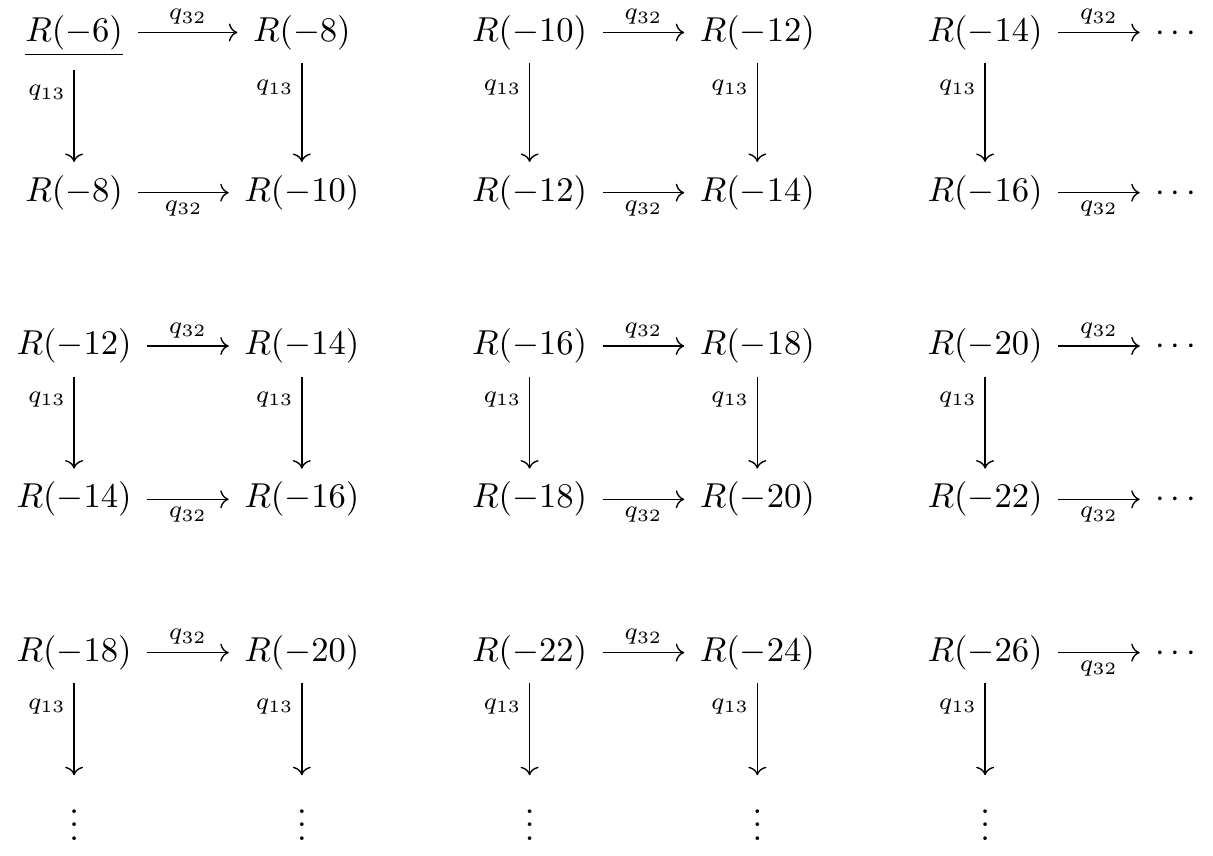}
\end{center}
\end{minipage}
\end{equation}
	
	This complex has split into a complex of the form $ \bigoplus_{i,j\geq 0} A(-4i-6j,2i+2j) ,$ where $A$ is the complex
	
	\begin{equation}\label{eqn-niceHH0P3}
\begin{minipage}{2.25in}
\begin{center}
\includegraphics[scale=1]{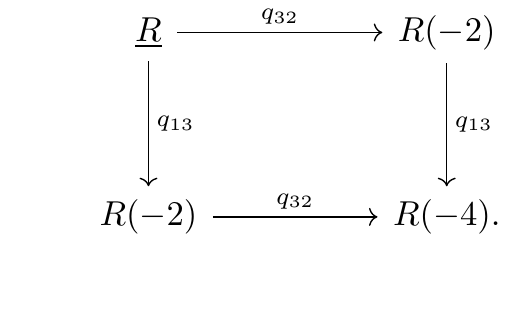}
\end{center}
\end{minipage}
\end{equation}
	
Since $R=\Q[x_1,x_2,x_3]$, we have an isomorphism $R \cong \Q[q_{13},q_{32},x_1]$.  Then $A$ is the Koszul complex associated to the regular sequence $\{q_{13},q_{32}\}$, and has homology isomorphic to $R/(q_{13},q_{32})$.  Thus, if we let $u_2$ and $u_3$ denote the perodicity maps from \S \ref{subsec-explicitCx} then we see $$\HHH(w(\pr{3}^\vee)) \simeq \Q[\xi_1,\xi_2,\xi_3,x_1,u_2,u_3](-4,2).$$	
\end{example}

We can extend this pattern to general $n$; we will see in \S \ref{subsec-matrixFact}
that $$\HHH(w(\pr{n}^\vee)) \simeq \Q[\xi_1,\ldots,\xi_n,x_1,u_2,\ldots,u_n](-2n,n),$$ where the $u_i$ are the periodicity maps from \S \ref{subsec-explicitCx}, $\xi_i$ are generators of an exterior algebra with $\deg{\xi_i}= q^{-2i}a$.

For general permutations $w\in S_n$,  $\HHH(w(\pr{n}^\vee))$ is not a free superpolynomial algebra. However, in Example \ref{ex-P2Unknot} we were able to describe $\HHH(\pr{2}^\vee)$ as a free superpolynomial algebra modulo certain relations (after grading shifts). This will be true for all $n$ and $w \in S_n$, as proven in \S \ref{subsec-matrixFact}.  

\subsection{The flag Hilbert scheme and $\End(\pr{n})$}
\label{subsec-EndPn}
The $q$-Young symmetrizers are a standard collection of primitive idempotents $p_T\in \HC_n$ in the Hecke algebra, indexed by Young tableaux.  One expects there to be complexes $P_T\in \K^-(\SBim_n)$ which categorify these elements.  If $T$ is the unique one-column tableau, then $P_T$ is our projector $\pr{n}$, and if $T$ is the unique one-row tableau, then $P_T$ is the idempotent complex constructed in \cite{Hog15}.  A construction of $P_T$ for all $T$ will be proposed by the second author and Ben Elias in \cite{ElHog15}.

In \cite{GoRa15}, E.~Gorsky and J.~Rasmussen conjecture that $\HHH(P_T)$ of the idemotents $P_T\in K^-(\SBim_n)$ can be described in terms of certain rings associated to the flag Hilbert scheme.  Our work here seems to suggest that their rings are in fact the homologies of the complexes $\Endg(P_T)$.  The relationship between $\HHH(P_T)$ and $\Endg(P_T)$ seems to be subtle in general.

In this paper we are concerned with the case where $T$ is the unique tableau whose shape is the one-column partition  $1+\cdots +1=n$, which is usually written $1^n$.   Throughout this section, set $P:=\pr{n}$ ($n$ is understood).  The purpose of this section is to prove a conjecture of Gorsky-Rasmussen in this case.  Here and below, we use the following:

\begin{definition}\label{def-homg}
if $M$ and $N$ are graded bimodules then $\Homg_{R^e}(M,N)$ will denote the graded space of all homogeneous bimodule maps.  If $C$ and $D$ are complexes of graded bimodules, then $\Homg_{R^e}(C,D)$ will denote the bigraded chain complex of all bihomogeneous $R^e$-linear maps.
\end{definition}

Our main theorem in this section is a computatio  of $H(\Homg(R,\pr{n}^\vee))$, thus a computation of $\HHH(\pr{n}^\vee)$.  The following ring is defined by Gorsky-Rasmussen in \cite{GoRa15}.

\begin{definition}\label{def-flagRing}
Let $x_1,\ldots,x_n$ denote formal (even) indeterminates of bidegree $q^2t^0$.  For all integers $1\leq i<j\leq n$, let $v_{ij}$ denote formal (even) indeterminates of bidegree $q^{2(i-j)-2}t^2$. Let $X$ and $V$ denote the matrices
\begin{equation}\label{eq:XandV}
X = \matrix{
x_1	&	1		&	0		&	\cdots	&	0& 0\\
0		&	x_2	&	1		&  \cdots & 0& 0 \\
0		&	0		&	x_3		&  \cdots & 0& 0 \\
\vdots		&	\vdots		&	\vdots		&  \ddots & \vdots & \vdots  \\
0		&	0		&	0		&  \cdots & x_{n-1} &1 \\
0		&	0		&	0		&  \cdots & 0 & x_n
}
,\ \ 
V = \matrix{
0	&	v_{12}		&	v_{13}		&	\cdots	& v_{1,n-1} &	v_{1n}\\
0		&	0	&	v_{23}		&  \cdots & v_{2,n-1} & v_{2n} \\
0		&	0		&	0		&  \cdots & v_{3,n-1} & v_{3n} \\
\vdots		&	\vdots		&	\vdots		&  \ddots &\vdots  & \vdots  \\
0		&	0		&	0		&  \cdots & 0 & v_{n-1,n} \\
0		&	0		&	0		&  \cdots & 0&  0
}
\end{equation}
Let $r_{ij}\in\Q[x_k,v_{ij}]$ denote the components of the commutator $[X,V]$.  Then set $E:=\Q[x_k,v_{ij}]/(r_{ij})$.
\end{definition}

Our main theorem in this section is:
\begin{theorem}\label{thm-untwistedHHH}
We have  isomorphisms of algebras
\[
E\cong \HHH^0(\pr{n}^\vee)\cong H(\Homg(R,\pr{n}^\vee)) \cong H(\Endg(P^\vee))\cong H(\Endg(P))
\]
As a consequence,
\[
\HHH(\pr{n}^{\vee})\cong E\otimes_{\Q} \Lambda[\xi_1,\ldots,\xi_n]
\]
where $\deg(\xi_i)=aq^{-2i}$.
\end{theorem}
The second isomorphism is by definition, the third isomorphism is by general arguments, since $\pr{n}^\vee$ is a unital idempotent, and the last isomorphism holds by an application of the duality functor.  Thus, the essential content is in the first isomorphism. The rest of this section is dedicated to the proof of this theorem.  First, we manipulate the definition of $E$ into a workable form.

\begin{proposition}\label{prop-flagRingExplicitRels}
Let us introduce the extraneous variables $v_{ii}$ for $1\leq i\leq n$.  Then $E$ is isomorphic to the quotient of $\Q[\xx,v_{ij}]_{1\leq i\leq j\leq n}$ by the relations
\begin{eqnarray}
v_{i,j} &= &v_{i-1,j-1} - (x_{i-1}-x_j)v_{i-1,j}\ \ \ \ \ \ \ \ \ \ \ (2 \leq i \leq  j\leq n) \label{eq-rel1}\\
v_{ii} & =& 0 \label{eq-rel2}
\end{eqnarray}
\end{proposition}
\begin{proof}
Straightforward exercise.
\end{proof}
The reader may be wondering why we introduce variables $v_{ii}$ and then immediately set them equal to zero.  The reason for this will become clear shortly.  First, note that the relation (\ref{eq-rel1}) allows us to write $v_{ii}$ in terms of $v_{11},v_{12},\ldots,v_{1,n}$.  The first few examples are
\begin{enumerate}
\item $v_{22} = v_{11} - q_{12} v_{12}$
\item $v_{33} = v_{11} - (q_{12}+q_{23})v_{12} + q_{13}q_{23}v_{13}$
\item $v_{44} = v_{11} - (q_{12}+q_{23}+q_{34})v_{12} + (q_{13}q_{23}+q_{13}q_{34}+q_{24}q_{34})v_{13} + q_{14}q_{24}q_{34}v_{14}$
\end{enumerate}
where $q_{ij}:=x_i-x_j$.  The reader may have noticed the similarity between these coefficients and the polynomials $a_{ij}(\xx,\yy)$ from Proposition \ref{prop-aijRelation}.  This is no accident.

\begin{proposition}\label{prop-theRels}
Using the relation (\ref{eq-rel1}) we have
\[
v_{kk} = \sum_{i=1}^k (-1)^{i-1} a_{ik}(\xx,\xx) v_{1i}
\]
\end{proposition}
\begin{proof}
Fix $k\in\{2,\ldots,n\}$, and consider the task of simplifying $v_{k,k}$.  Each application of the relation (\ref{eq-rel1}) reduces the first index by 1, and doubles the number of terms.  Thus, $k-1$-applications of this relation allows us to write $v_{kk}$ as a $\Q[\xx]$-linear combination of the variables $v_{1,i}$, with $2^{k-1}$ terms.  We may organize the terms as follows: let $\Gamma$ denote a word of length $k-1$ in the letters `R' (for ``right branch'') and `L' (for ``left branch'').  An occurence of `L' contributes a factor of 1, and an occurence `R' in the $i$-th coordinate contributes a factor of $-(x_{k-i}-x_{k-i+j+1})$ where $j$ is the number of occurences of `R' among the first $i-1$ coordinates.  For instance, the word RRLLR (corresponding to $k=6$ contributes a factor of
\[
(-1)^3(x_5-x_6)(x_4-x_6)(x_1-x_4)
\]
More precisely, let $i_1<\cdots < i_j$ denote the indices such that   $\Gamma_{i_j} = \text{`R'}$, and let $|\Gamma|=j$ denote the number of occurences of `R.'  Then define a ``Boltzmann weight'' by
\[
\wt(\Gamma) = (-1)^j(x_{k-i_1} - x_{k-i_1+1})(x_{k-i_2} - x_{k-i_2+2})\cdots (x_{k-i_j} - x_{k-i_j+j})
\] 
By construction, we have
\[
v_{kk}=\sum_{j=0}^{k-1} (-1)^{j} v_{1,j+1}\sum_{|\Gamma|=j} \wt(\Gamma) 
\]
Furthermore, comparison with the formula for $a_{ij}(\xx,\xx)$ in Definition \ref{def-aij} makes it clear that $\sum_{|\Gamma|=j} \wt(\Gamma) = (-1)^j a_{j+1,k}(\xx,\yy)$.  This completes the proof.
\end{proof}

\begin{proof}[Proof of Theorem \ref{thm-untwistedHHH}]
First, recall the differential bigraded algebra $A_n$ from Definition \ref{def-dga}: $A_n=\Q[\xx,\yy,u_1,\ldots,u_n,\theta_1,\ldots,\theta_n]$ with $\Q[\xx,\yy,u_1,\ldots,u_n]$-linear differential determined by $d_A(\theta_j)=\sum_{i=1}^k u_i a_{ij}(\xx,\yy)$, together with the graded Leibniz rule.  Then $P^\vee \simeq M_n = A_n/ I_n A_n$, with $\Q[\xx,\yy,u_1,\ldots,u_n]$-linear differential $d_M$ determined by $d_M(1)=\sum_{i=1}^n(y_i-x_i)\theta_i$ and the graded Leibniz rule $d_M(am) = d_A(a)m+(-1)^{|a|}ad_M(m)$ for all $a\in A_n$, $m\in P^\vee$.  This complex is then shifted so that $1$ lies in degree $q^{-2\ell}t^0$.

Note that $d_M(1)$ is zero modulo $(x_i-y_i)$.  By Proposition \ref{prop-HHIn}, the functor $\Homg_{\Q[\xx,\yy]}(\Q[\xx],-)$ sends $\Q[\xx,\yy]\mapsto \Q[\xx](2\ell)$ and $f(\xx,\yy)\mapsto f(\xx,\xx)$.  The grading shifts cancel, and we see that $\Homg(\Q[\xx],M)$ is isomorphic to the superpolynomial ring $\Q[\xx,u_1,\ldots,u_n,\theta_1,\ldots,\theta_n]$ with differential determined by
\[
d(u_k) = 0 \ \ \ \ \ \ \ \ \ \ \ \ d(\theta_j)=\sum_{i=1}^k u_i a_{ij}(\xx,\xx)
\]
together with the graded Leibniz rule.  This is the Koszul complex associated to the sequence
\[
(u_1), \ \  \ (u_1 + u_2 a_{22}(\xx,\xx)), \ \ \   (u_1+ u_2 a_{23}(\xx,\xx) + u_3 a_{33}(\xx,\xx)), \ \ \ \cdots \ \ \ 
\]
in $\Q[\xx,u_1,\ldots,u_n]$.  It is easy to see that this is a regular sequence, hence the homology is simply the quotient of $\Q[\xx,u_1,\ldots,u_n]$ by the ideal generated by the above elements.  By Proposition \ref{prop-flagRingExplicitRels} and Proposition \ref{prop-theRels} it follows that the homology of $\Homg_{\Q[\xx,\yy]}(\Q[\xx],M_n)$ is isomorphic to $E$, via an isomorphism which sends $u_k\mapsto (-1)^{k-1}v_{1,k}$.  This proves that $H(\Homg(R,\pr{n}^\vee))\cong E$.  The isomorphism $\Homg(R,\pr{n}^\vee)\cong H(\Endg(\pr{n}))$ follows from general theory of categorical idempotents (see \S 5 of \cite{Hog15}).  This completes the proof.
\end{proof}

\subsection{Homology of twisted projectors}
\label{subsec-matrixFact}

We first state the Gorsky-Rasmussen conjecture for twisted projectors. Let $x_i$ denote a formal indeterminate of bidegree $q^2a$ and $v_{ij}$ denote a formal indeterminate of bidegree $q^{2(i-j)-2}t^2$ as before.

\begin{conjecture}\label{conj-GRtwisted}
Recall the ring $E$ from Definition \ref{def-flagRing}.
Then for each permutation $w\in S_n$, the homology $\HHH^0(w(\pr{n}^\vee))$ is isomorphic to $E / J_w$, where $J_w$ is the ideal generated by the differences $x_{w(i)}-x_i$ for all $1\leq i\leq n$.
\end{conjecture}

The symmetric group acts on $\HHH(\pr{n}^\vee)$ by conjugation by Rouquier complexes.  More precisely, any chain map $\phi:R\rightarrow P$ gives rise to $R\simeq F_wF_w\inv\rightarrow F_wPF_w\inv$ for the Rouquier complex $F_w$ associated to the positive braid lift of $w$.  Since $P$ absorbs Rouquier complexes we have an equivalence $F_wPF_w\inv\simeq P$.  Pre-composing with the map $F\phi F\inf$ defines $w(\phi):R\rightarrow P$.  This gives rise to an action of $S_n$ on $E\cong \HHH^0(\pr{n}^\vee)$.  If $w=vwv\inv$, then $\phi\mapsto v(\phi)$ descends to an isomorphism $E/J_w \rightarrow E/J_{vwv\inv}$.  Thus, both algebras in Conjecture \ref{conj-GRtwisted} depend only on the conjugacy class of $w$ up to isomorphism.

It will be helpful to encode $\pr{n}^\vee$ in the form of a Koszul matrix factorization, the basics of which we now recall.  Let $S$ be a commutative ring and $z\in S$ a fixed element.  For us, a $z$-matrix factorization will be a $\Z/2$-graded $S$-module $M=M_0\oplus M_1$ together with an $S$-linear endomorphism $d\in \End_S(M)$ such that
\begin{enumerate}
\item $d$ is odd. That is, $d$ restricts to $S$-module maps $M_0\rightarrow M_1$ and $M_1\rightarrow M_0$.
\item $d^2=z$, where $z\in S$ is regarded as an endomorphism of $M$.
\end{enumerate}
We call $z\in S$ the \emph{potential}.
\begin{remark}
We regard $S$ as $\Z/2$-graded, where each element of $S$ is even.  A map of $\Z/2$-graded objects is called even it preserves parity, and odd if it swaps parity.
\end{remark}

A \emph{morphism} $(M,d)\rightarrow (N,d')$ between matrix factorizations with the same potential is an even $S$-module map $f:M\rightarrow N$ such that $d'\circ f = f\circ d$.  If $(M,d)$ is a $z$-matrix factorization and $(M,d')$ is a $z'$-matrix factorization, then the tensor product $M\otimes_S M'$ is naturally a $(z+z')$-matrix factorization with $d_{M\otimes M'} = d\otimes\Id_{M'} \pm \Id_M\otimes d'$.

Now, let $a,b\in S$ a pair of elements.  We have the basic \emph{Koszul matrix factorization}
\[
( a | b) \ :=  \ 
\begin{minipage}{1.5in}\begin{tikzpicture}
\tikzstyle{every node}=[font=\small]
\node (a) at (0,0) {$S$};
\node (b) at (3,0) {$S$};
\path[->,>=stealth',shorten >=1pt,auto,node distance=1.8cm,
  thick]
([yshift=3pt] a.east) edge node[above] {$a$}		([yshift=3pt] b.west)
([yshift=-2pt] b.west) edge node[below] {$b$}		([yshift=-2pt] a.east);
\end{tikzpicture}\end{minipage}
\]
We regard the left-most $S$ as sitting in even degree $0\in \Z/2$, and the right-most $S$ as sitting in odd degree $1\in \Z/2$.  This is a matrix factorization with potential $ab$.  In general given elements $a_i,b_i\in S$ with $1\leq i\leq r$, we define
\[
\left(
\begin{tabular}{c|c}
$a_1$ & $b_1$ \\
$a_2$ & $b_2$ \\
$\vdots$ & $\vdots$\\
$a_r$ & $b_r$
\end{tabular}\right)
:= \bigotimes_{i=1}^r\Big(
\begin{minipage}{1.5in}\begin{tikzpicture}
\tikzstyle{every node}=[font=\small]
\node (a) at (0,0) {$S$};
\node (b) at (3,0) {$S$};
\path[->,>=stealth',shorten >=1pt,auto,node distance=1.8cm,
  thick]
([yshift=3pt] a.east) edge node[above] {$a_i$}		([yshift=3pt] b.west)
([yshift=-2pt] b.west) edge node[below] {$b_i$}		([yshift=-2pt] a.east);
\end{tikzpicture}
\end{minipage}\Big)
\]
The tensor product here is over $S$.  This is a matrix factorization with potential $\sum_{i=1}^r a_ib_i\in S$.  If $\mathbf{a}=(a_1,\ldots,a_r)$ and $\mathbf{b}=(b_1,\ldots,b_r)$, then we also write $(\aa,\bb)$ for the corresponding Koszul complex.  Usually the ring $S$ will be implicit, but if we wish to include it in the notation, we will do so with a subscript, as in $(\aa,\bb)_S$.  Similarly, if $M$ is an $S$-module, we denote $(\aa,\bb)_M := (\aa,\bb)\otimes_S M$.  

\begin{remark} A Koszul matrix factorization is the super-position of two Koszul complexes.  There is the ``forward differential'' which is the Koszul differential associated to the  sequence $a_1,\ldots,a_r$, and there is the ``backward differential'' which is the Koszul differential associated to $b_1,\ldots,b_r$.  In particular, as an $S$-module $(\aa,\bb)$ is an exterior $S$-algebra with $r$-generators.
\end{remark}

\begin{remark}
Here is another description of a Koszul matrix factorization.  If $\mathbf{a}=(a_1,\ldots,a_r)$ and $\mathbf{b}=(b_1,\ldots,b_r)$ are given sequences of elements of $S$, then the Koszul matrix factorization $(0,\mathbf{b})$ is a differential graded algebra
\[
(0,\mathbf{b}) = \Lambda[\theta_1,\ldots,\theta_r] \ \ \text{ with differential } \ \ d_A(\theta_i)=b_i
\]
This is simply the usual Koszul complex associated to $b_1,\ldots,b_r$.  Let us call this dg algebra $A$.

Now, $M:=(\mathbf{a},\mathbf{b})$ can be thought of as a dg $A$-module generated by one element, with differential $d_M(1) = \sum_{i=1}^r a_i\theta_i$.  We leave it to the reader to check this.
\end{remark}

\begin{remark}\label{rmk-gradedMFs}
We will be primarily interested in the case where $S$ is a bigraded ring.  The degree of a bihomogeneous element $a\in S$ is written $\deg(a)=q^i t^j$, and $t$ is regarded as the homological degree.  We will require our matrix factorizations to be bigraded $S$ modules, with parity given by the mod 2 reduction of homological degree, and we require the differential to satsify $\deg(d) = q^0t^1$.  If $\deg(b)=q^it^{j}$ and $\deg(a)=q^{-i}t^{1-j}$, then we incorporate the grading shifts as follows:
\[
(a|b) = \Big(
\begin{minipage}{1.63in}\begin{tikzpicture}
\tikzstyle{every node}=[font=\small]
\node (a) at (0,0) {$S$};
\node (b) at (3,0) {$S(i)\ip{j-1}$};
\path[->,>=stealth',shorten >=1pt,auto,node distance=1.8cm,
  thick]
([yshift=3pt] a.east) edge node[above] {$a$}		([yshift=3pt] b.west)
([yshift=-2pt] b.west) edge node[below] {$b$}		([yshift=-2pt] a.east);
\end{tikzpicture}
\end{minipage}
\Big)_S
\]
so that the right-pointing arrow has degree $q^0t^0$ and the left-pointing arrow has degree $q^0t^1$.  Note that in this case
\[
(a|b) \cong q^{-i}t^{1-j} (b|a)
\]
\end{remark}

Before constructing $\pr{n}^\vee$ as a matrix factorization, we recall two basic operations on Koszul matrix factorizations which preserve the isomorphism type of the matrix factorization. The next two propositions are proven in \cite{KR08}.

\begin{proposition} \label{prop-changeBasisMF}
Let $\mathbf{a} = (a_1,\ldots,a_n)$ and $\mathbf{b} = (b_1,\ldots,b_n)$ be lists of elements in $S$. Let $(\mathbf{a},\mathbf{b})$ be the associated Koszul matrix factorization. The change of basis transformation
	\[
	\left(
	\begin{tabular}{c|c}
	$a_i$ &$b_i$ \\
	$a_j$ & $b_j$
	\end{tabular}\right) \ \ \xrightarrow \ \ \left(
	\begin{tabular}{c|c}
	$a_i$ &$b_i+\lambda b_j$ \\
	$a_j-\lambda a_i$ & $b_j$
	\end{tabular}\right),
	\]
where all other rows are fixed, yields isomorphic matrix factorizations for all $\lambda \in S$.  Similarly, if $\lambda\in S$ is invertible, then $(a_i | b_i) \cong (\lambda a_i | \lambda\inv b_i)$.	
\end{proposition}

Now assume that $S$ is a polynomial ring. We can also simplify Koszul matrix factorizations by canceling rows of the form $(0|b_i)$ or $(a_i|0)$ (called \emph{exclusion of variables}), subject to certain restrictions. Precisely, suppose $(\mathbf{a},\mathbf{b})$ is a Koszul matrix factorization with potential $z = \sum_{i=1}^n (a_ib_i)$. Let $y$ be a generator of $S$ and write $S = S'[y]$. Assume $z \in S'$ and that one of the rows of $(\mathbf{a},\mathbf{b})
$ has the form $(0, y - p)$ for $p \in S'$.
\begin{proposition}\label{prop-exclusionVarMF}
Given the above hypotheses, $(\mathbf{a},\mathbf{b})_{S}$ is homotopy equivalent to the matrix factorization of potential $\sum(a_ib_i)$ over $S'$ where we remove the row $(0,y-p)$ and substitute $p$ for $y$ everywhere in all other rows.
\end{proposition}
We have not yet introduced the notion of homotopy equivalence of matrix factorizations.  There will be no need to do so, since we are mostly interested in the case when the matrix factorization $M$ has zero potential (that is, $M$ is a chain complex), in which case the notion homotopy equivalence is standard. The basic idea of exclusion of a variable is that $(0 | y-p)$ is the complex
\[
\begin{diagram} 0 & \rTo & S'\otimes \Q[y] & \rTo^{1\otimes y-p\otimes 1}  & S'\otimes \Q[y] & \rTo & 0\end{diagram}
\]
This complex is homotopy equivalent to $S'$ via Gaussian elimination.  We will now specialize to the case of interest.

\begin{definition}\label{def-MasMF}
Fix an integer $n$.  Let $S = \Q[\xx,\yy,u_1,\ldots,u_n]$, bigraded so that $\deg(u_k)=t^2q^{-2k}$.  Define elements $b_j\in S$ by $b_j = \sum_{i=1}^j u_i a_{ij}(\xx,\yy)$, where the polynomials $a_{ij}(\xx,\yy)$ are as in Definition \ref{def-aij}.  Consider the Koszul matrix factorization
\begin{equation}\label{eq-MasMF}
 q^{-2\ell}\left(
\begin{tabular}{c|c}
$p_{11}$ &$b_1$ \\
$p_{22}$ & $b_2$ \\
$\vdots$ & $\vdots$\\
$p_{nn}$ & $b_n$
\end{tabular}\right) \ =\  q^{-2\ell}\ \bigotimes_{j=1}^n \Big(
\begin{minipage}{1.63in}\begin{tikzpicture}
\tikzstyle{every node}=[font=\small]
\node (a) at (0,0) {$S$};
\node (b) at (3,0) {$S(-2)\ip{1}$};
\path[->,>=stealth',shorten >=1pt,auto,node distance=1.8cm,
  thick]
([yshift=3pt] a.east) edge node[above] {$p_{jj}$}		([yshift=3pt] b.west)
([yshift=-2pt] b.west) edge node[below] {$b_j$}		([yshift=-2pt] a.east);
\end{tikzpicture}
\end{minipage}
\Big)_S \otimes_S \bar{S}
\end{equation}
Here $\bar{S}$ is the quotient of $S$ in which we identify $e_k(\xx)=e_k(\yy)$ for all $k\in\{1,\ldots,n\}$, and $p_{ij}(\xx,\yy)=y_i-x_j$.  Note that $\deg(p_{ii}b_i)=q^0t^1$, so that (\ref{eq-MasMF}) is graded as in Remark \ref{rmk-gradedMFs}.
\end{definition}

The following is clear from the definitions.

\begin{proposition}
The matrix factorization defined in Definition \ref{def-MasMF} is equal to $M_n$ (from \S \ref{subsec-explicitCx}) as a bigraded chain complex. Therefore the matrix factorization (\ref{eq-MasMF}) is isomorphic to $\pr{n}^\vee$ as a bigraded chain complex by Theorem \ref{thm-MisP}.
\end{proposition}

For the grading shift of $q^{-2\ell}$, recall that $B_{w_0} = q^{-\ell}\Q[\xx,\yy]/I$ and the degree zero chain bimodule of $M_n$ is $q^{-\ell}B_{w_0}=q^{-2\ell}\Q[\xx,\yy]/I$ (see also Definition \ref{def-dgMod}).

\begin{example}
For $n=2$, the complex looks like
\[
M_2 = \left(
\begin{tabular}{c|c}
$p_{11}$ & $u_1$ \\
$p_{22}$ & $u_1+p_{12}u_2$ 
\end{tabular}\right)_{\bar{S}}
\]
After a change of basis (Proposition \ref{prop-changeBasisMF}) this is isomorphic to
\[
\left(
\begin{tabular}{c|c}
$p_{11} + p_{22} $ & $u_1$ \\
$p_{22}$ & $p_{12}u_2 $
\end{tabular}\right)_{\bar{S}}
\]
Note that $p_{11}+p_{22}$ is zero in ${\bar{S}}$, so we obtain
\[
M_2 \cong \left(
\begin{tabular}{c|c}
$0 $ & $u_1$ \\
$p_{22}$ & $p_{12}u_2 $
\end{tabular}\right)_{\bar{S}} \simeq (p_{22} | p_{12})_{\bar{S}/(u_1)}
\]
by exlusion of variables (Proposition \ref{prop-exclusionVarMF}).  This is a very compact way of expressing Diagram \ref{eqn-explicitP2} (with the arrows reversed).
\end{example}

\begin{example}
For $n=3$ we have
\[
M_3 = \left(
\begin{tabular}{c|c}
$p_{11}$ & $u_1$ \\
$p_{22}$ & $u_1+p_{12}u_2$ \\
$p_{33}$ & $u_1+(p_{12}+p_{23})u_2 + p_{13}p_{23}u_3$ 
\end{tabular}\right)_{\bar{S}}
\]
After a change of basis we obtain that $M_3$ is isomorphic to
\[
\left(
\begin{tabular}{c|c}
$p_{11}+p_{22}+p_{33}$ & $u_1$ \\
$p_{22}$ & $p_{12}u_2$ \\
$p_{33}$ & $(p_{12}+p_{23})u_2 + p_{13}p_{23}u_3$ 
\end{tabular}\right)
\ \ 
\]
which by Proposition \ref{prop-changeBasisMF} is isomorphic to 
\[ 
\left(
\begin{tabular}{c|c}
$0$ & $u_1$ \\
$p_{22}$ & $p_{12}u_2$ \\
$p_{33}$ & $(p_{12}+p_{23})u_2 + p_{13}p_{23}u_3$ 
\end{tabular}\right).
\]
An application of Proposition \ref{prop-exclusionVarMF} yields
\[
M_3 \simeq 
\left(\begin{tabular}{c|c}
$p_{22}$ & $p_{12}u_2$ \\
$p_{33}$ & $(p_{12}+p_{23})u_2 + p_{13}p_{23}u_3$ 
\end{tabular}\right)_{\bar{S}/(u_1)},
\]
which is simply the complex from Diagram \ref{eqn-explicitP3}, with the arrows reversed.
\end{example}

Now, let $w\in S_n$ be given.  We wish to compute $H(\HH^0(R,w(\pr{n}^\vee)))$.  This homology depends only on the conjugacy class of $w$, so we may as well assume that $w$ has the special form (\ref{eq-specialw}).
\begin{equation}\label{eq-specialw}
w = (1,\ldots,m_1)(m_1+1,\ldots,m_2)\cdots (m_{r-1}+1,\ldots,m_r)
\end{equation}
written in cycle notation, for some integers $1\leq m_1<m_2<\cdots <m_r = n$.  Applying the permutation $w$ to $M_n$ simply permutes the variables $y_i$ in the formula (\ref{eq-MasMF}).  Taking $\HH^0(R,-)$ sends $B\mapsto R(\ell)$.  Thus $\HH^0$ dends $\Q[\xx,\yy]/I\mapsto q^{2\ell}\Q[\xx]$.  This cancels the shift of $q^{-2\ell}$ in (\ref{eq-MasMF}) then identifies $\xx$ with $\yy$.  The result is:

\begin{lemma}\label{lemma-HomRPsimp}
We have
\begin{equation}\label{eq-HomMn}
\HH^0(R,w(M_n)) \ \ \cong \ \ 
\left(
\begin{tabular}{c|c}
$x_{w(1)}-x_1$ &$\bar{b}_1$ \\
$x_{w(2)}-x_2$ & $\bar{b}_2$ \\
$\vdots$ & $\vdots$\\
$x_{w(n)}-x_n$ & $\bar{b}_n$
\end{tabular}\right)_{\Q[\xx,u_1,\ldots,u_n]}
\end{equation}
where $\bar{b}_j = \sum_{i=1}^j u_i a_{ij}(\xx,w(\xx))$.\qed
\end{lemma}

\begin{proposition}\label{prop-twistedDtheta}
If $j_1$ and $j_2$ are in the same $w$-orbit, then $$a_{i,j_1}(\xx,w(\xx))=a_{i,j_2}(\xx,w(\xx)).$$
\end{proposition}
\begin{proof}
We are assuming that $w$ has a special form, so $j_1$,$j_2$ are in the same $w$-orbit if and only if there is a $1\leq k\leq r$ such that $m_{k-1}<j_1,j_2\leq m_k$.  So without loss of generality, assume that $j_1= j_2+1$.  From the definition of $a_{ij}(\xx,\yy)$ it is easy to see that $a_{i,j+1}(\xx,\yy) - a_{ij}(\xx,\yy)$ is divisible by $y_j-x_{j+1}$, for all $1\leq i\leq j+1$, hence is zero upon specializing $y_{w(i)}=x_i$.  This completes the proof.
\end{proof}

We will now simplify the complex (\ref{eq-HomMn}) by focusing on one block at a time.  Let $w$ be as above; for each cycle $(m,m+1,\ldots,m+j)$ of $w$ we have the following tensor factor of the right-hand side of (\ref{eq-HomMn}):
\[
\left(
\begin{tabular}{c|c}
$x_{m+1}-x_{m}$ & $\sum_{i=1}^m u_ia_{i,m}(\xx,w(\xx))$ \\
$x_{m+2} - x_{m+1}$ & $\sum_{i=1}^{m+1} u_ia_{i,m+1}(\xx,w(\xx))$ \\
$\vdots$ & $\vdots$\\
$x_{m+j} - x_{m+j-1}$ & $\sum_{i=1}^{m+j-1} u_ia_{i,m+j-1}(\xx,w(\xx))$\\
$x_{m} - x_{m+j}$  & $\sum_{i=1}^{m+j}u_ia_{i,m+j}(\xx,w(\xx))$
\end{tabular}\right)
\]
Set $\a_i = x_i - x_{i+1}$.  Then $x_m-x_{m+j} = \a_m+\a_{m+1}+\cdots +\a_{m+j-1}$.  By Proposition \ref{prop-twistedDtheta}, all of the entries in the right column are equal to one another.  Rewriting gives
\[
\left(
\begin{tabular}{c|c}
$-\a_m $&$\sum_{i=1}^{m+j}u_ia_{i,m+j}(\xx,w(\xx))$ \\
$-\a_{m+1} $&$\sum_{i=1}^{m+j}u_ia_{i,m+j}(\xx,w(\xx))$ \\
$\vdots $&$ \vdots$\\
$-\a_{m+j-1} $&$\sum_{i=1}^{m+j}u_ia_{i,m+j}(\xx,w(\xx))$ \\
$\a_m + \cdots + \a_{m+j-1} $&$\sum_{i=1}^{m+j}u_ia_{i,m+j}(\xx,w(\xx))$ 
\end{tabular}\right)
\]
After a change of basis this is isomorphic to 
\[
\left(
\begin{tabular}{c|c}
$\a_m $&$ 0 $\\
$\a_{m+1} $&$ 0$\\
$\vdots $&$ \vdots$\\
$\a_{m+j-1} $&$ 0$\\
$0 $&$\sum_{i=1}^{m+j}u_ia_{i,m+j}(\xx,w(\xx))$ 
\end{tabular}\right)
\]
Applying such transformations to all of the cycles in $w$, we obtain the following:
\begin{proposition}\label{prop-twistedSimplifications}
We have
\[
\HH^0(R,w(M_n))\cong \left(
\begin{tabular}{c|c}
$\a_1 $&$ 0 $\\
$\vdots $&$ \vdots $\\
$\a_{m_1-1}  $&$ 0 $\\
$0$ & $\sum_{i=1}^{m_1}u_ia_{i,m_1}(\xx,w(\xx))$ \\
$\a_{m_1+1}  $&$ 0 $\\
$\vdots $&$ \vdots $\\
$\a_{m_2-1} $&$ 0 $\\
$0 $& $\sum_{i=1}^{m_2}u_ia_{i,m_2}(\xx,w(\xx))$ \\
$\vdots $&$ \vdots $\\
$\a_{m_{r-1}+1} $&$ 0 $\\
$\vdots $&$ \vdots $\\
$\a_{m_r-1} $&$ 0 $\\
$0 $& $\sum_{i=1}^{m_r}u_ia_{i,m_r}(\xx,w(\xx))$ \\
\end{tabular}\right)_{\Q[\xx,u_1,\ldots,u_n]}
\]\qed
\end{proposition}

We are ready to prove our main theorem:
\begin{theorem}
Let $w\in S_n$ be a permutation with $r$-cycles.  Then $\HHH^0(w(\pr{n}^\vee))$ is isomorphic to $E/(x_i-x_{w(i)}|1\leq i\leq n)$, shifted so that $1$ lies in degree $(q^{-2}t)^{n-r}$, where $E$ is as in Definition \ref{def-flagRing}.  That is to say, Conjecture \ref{conj-GRtwisted} is true.
\end{theorem}
\begin{proof}
Fix $w\in S_n$; without loss of generality, we assume that $w$ has the special form (\ref{eq-specialw}).  Define the following ideals in $\Q[\xx,u_1,\ldots,u_n]$:
\begin{enumerate}
\item $I=(\a_i)$, where $i\in\{1,\ldots,n\}\setminus\{m_1,\ldots,m_r\}$.
\item $I'=(x_i-x_{w(i)})$ where ${1\leq i\leq n}$.
\item $J$ is the ideal generated by $\sum_{i=1}^{m_k}u_i a_{i,m_k}(\xx,w(\xx))$ where $1\leq k\leq r$.
\item $J'$ is the ideal generated by $\sum_{i=1}^ju_i a_{i,j}(\xx,w(\xx))$ where $1\leq j\leq n$.
\item $J''$ is the ideal generated by $\sum_{i=1}^ju_i a_{i,j}(\xx,\xx)$ where $1\leq j\leq n$.
\end{enumerate}
It is clear that $I=I'$ and $I'+J'=I'+J''$.  Proposition \ref{prop-twistedDtheta} implies that $J=J'$, so we conclude that $I+J = I'+J''$.

The matrix factorization from Proposition \ref{prop-twistedSimplifications} is clearly isomorphic to the Koszul complex associated to the sequence $\{\a_i: i\in \{1,\ldots,n\}\setminus\{m_1,\ldots,m_r\}\}\cup \{\overline{b}_{m_1},\ldots,\overline{b}_{m_r}\}\subset R[u_1,\ldots,u_n]$, where $\overline{b}_j = \sum_{i=1}^ju_i a_{i,j}(\xx,w(\xx))$.  This sequence is easily seen to be regular, so
\[
H(\HH^0(R,w(M_n))\cong \Q[\xx,u_1,\ldots,u_n]/(I+J),
\]
which is isomorphic to $\Q[\xx,u_1,\ldots,u_n]/(I'+J'')$ by the remarks above.  On the other hand, $\Q[\xx,u_1,\ldots,u_n]/J''\cong E$ by Theorem \ref{thm-untwistedHHH}, so that $H(\HH^0(R,w(M_n))$ is isomorphic to $ E/(x_i-x_{w(i))})$, as claimed.

Finally, the overall grading shift is due to the fact (see Remark \ref{rmk-gradedMFs}) that $(\a_i| 0 )\cong q^{-2}t(0|\a_i)$.
\end{proof}

\begin{corollary}\label{cor-PSeries}
Let $w$ be a permutation with $r$ cycles. Then $\HHH(w(\pr n^\vee))$ has Poincar\'{e} series $$\mathcal{P}(q,a,t)=\frac{(q^{-2}t)^{n-r} (1-t^2q^{-2})^{r-1}\prod_{i=1}^n(1+aq^{-2i})}{(1-q^2)^r\prod_{j=2}^n(1-t^2q^{-2j})}$$ Furthermore, $\mathcal{P}(q,a,-1)$ is the unknot of the $1^n$-colored HOMFLYPT polynomial.
\end{corollary}

\begin{proof}
Once again, we fix $w \in S_n$ and assume $w$ has the special form (\ref{def-flagRing}). Also recall that we represent $\deg(z)$ by a monomial $q^ia^jt^k$. In the proof of Theorem \ref{thm-untwistedHHH}, we show that the ideal $(I + J)$ is generated by a regular sequence. A result of Stanley \cite{Stan78} states the following: If $A = \Q[z_1,\ldots,z_m]$ is a graded ring and $I = (a_1,\ldots,a_k)$ is an ideal generated the regular sequence $\{a_1,\ldots,a_k\}$, then the Poincar\'{e} series of $R/I$ is given by $$\mathcal{P}(A/I) = \frac{\prod_{i=1}^k(1 - \deg(a_i))}{\prod_{j=1}^m(1-\deg(z_j))}.$$

Therefore, it is easy to see that the denominator of the Poincar\'{e} series of $\HHH^0(w(\pr n^\vee))$ is given by $(1-q^2)^n\prod_{j=1}^n(1-t^2q^{-2j})$. As for the numerator, the regular sequence generating $I + J$ is given by $n-r$ elements of the form $\alpha_i$ and $r$ elements of the form $\bar{b}_{m_j}$. Each $\alpha_i$ has degree $q^2$, and each $\bar{b}_{m_k}$ has degree $t^2q^{-2}$.  Therefore the numerator has the form $(q^{-2}t)^{n-r}(1-q^2)^{n-r}(1-t^2q^{-2})^{r}$, where $(q^{-2}t)^{n-r}$ takes into account the overall grading shift.  The denominator has the form $(1-q^2)^n\prod_{i=1}^n(1-t^2q^{-2i})$, corresponding to the generators $x_i$, $u_i$.  The $(1-t^2q^{-2})$ factor cancels with one factor from the numerator, so $$\mathcal{P}(\HHH^0(w(\pr n^\vee)))=\frac{(q^{-2}t)^{n-r} (1-t^2q^{-2})^{r-1}}{(1-q^2)^r\prod_{j=2}^n(1-t^2q^{-2j})}.$$ 

To compute the Poincar\'{e} series for $\HHH(w(\pr n^\vee))$ we multiply by $\prod_{i=1}^n(1+aq^{-2i})$. This is because $\HHH(w(\pr n^\vee)) \cong \HHH^0(w(\pr n^\vee)) \otimes_\Q \Q[\xi_1,\ldots\xi_n]$ and that $\deg(\xi_k)= aq^{-2k}$. We leave checking the second claim to the reader.
\end{proof}

One may check that when $\mathcal{P}(\HHH(w(\pr n ^\vee)))$ is expanded into a Laurent series that all coefficients are positive, though this is not obvious upon first glance. This, of course, is expected because the coefficients are dimensions of vector spaces.

\section{Combinatorial results}
\label{sec-combinatorics}
In this section we record some combinatorial results necessary for giving an explict dg-module construction of $\pr n$ in \S \ref{sec-projector}. These results are known to experts, but we present elementary proofs of these facts for completeness.
\subsection{The Frobenius extension $R^W\subset R^I$}
\label{subsec-frobExtension}
Throughout this section, we will let $I:=\{1,\ldots,n-2\}$, so that $R^I=R^{S_{n-1}}$ is the algebra of polynomials $f\in \Q[x_1,\ldots,x_n]$ which are symmetric in $x_1,\ldots,x_{n-1}$.  In this section we recall the result that $R^W\subset R^I$ is a Frobenius exension, and we give an explicit dual pair of basis of $R^I$ as a free $R^W$-module.  Our main result in this section is the following:

\begin{theorem}\label{thm-frobBases}
$R^I$ is free of rank $n$ over $R^W$.  There are two natural ordered bases, given by
\[
\{(-1)^ix_n^i\}_{i=0}^{n-1} \ \ \ \ \ \ \  \ \ \ \ \{e_{n-1-i}(x_1,\ldots,x_{n-1})\}_{i=0}^{n-1}
\]
These bases are dual with respect to the trace pairing $(f,g):=\partial_{1,2,\ldots,n-1}(fg)$, where $\partial_{1,2,\ldots,n-1}$ denotes a divided difference operator (Definition \ref{def-dividedDiff}).  In particular $R^W\subset R^I$ is a Frobenius extension.

\end{theorem}
Iterating this, we recover as a corollary the classical result that $R$ is free over $R^W$ of rank $n!$.  If we consider $R$ to be graded with $\deg{x_i} = q^2$, and carefully keep track of gradings, then we can recover Proposition \ref{prop-BwOnB0}. Our main use for this theorem is in the proof of Proposition \ref{prop-canonicalMaps}, which itself is used in the proof of Theorem \ref{thm-MisP}.  The proof of Theorem \ref{thm-frobBases} occupies the remainder of this section.  We first introduce the trace $R^I\rightarrow R^W$, which is defined in terms of the \emph{divided difference operators}.

\begin{definition}\label{def-dividedDiff}
	Let $R:=\Q[x_1,\ldots,x_n]$, let $i\in\{1,\ldots,n-1\}$ be given, and let $s=(i,i+1)\in S_n$ denote the simple transposition which swaps $i$ and $i+1$.  The \emph{divided difference operator} $\partial_i:R\rightarrow R$ is defined by
	\[
	\partial_i(f):=\frac{f-s(f)}{x_i-x_{i+1}}
	\] 
	for each $f\in R$.  We may also write $\partial_s = \partial_i$, by abuse.
\end{definition}
Note that $f-s(f)$ is anti-symmetric in $x_i$ and $x_{i+1}$, hence is divisible by $x_i-x_{i+1}$.  So $\partial_s(f)$ is indeed a polyomial.  The following basic properties are easily checked:
\begin{proposition}\label{prop-dividedDiffProps}
	The divided difference operators satisfy:
	\begin{enumerate}
		\item $\partial_i^2=0$
		\item $\partial_i\partial_{i+1}\partial_i=\partial_{i+1}\partial_i\partial_i$
		\item if $f$ is symmetric in $x_i,x_{i+1}$, then $\partial_i(fg)=f\partial_i(g)$
		\item in general, $\partial_s(fg) = \partial_s(f)g+s(f)\partial_s(g)$.
	\end{enumerate}
\end{proposition}
\begin{proof}
	Straightforward.
\end{proof}

The following computation will be very useful:
\begin{proposition}\label{prop-traceOfXn}
	Fix an integer $n\geq 1$, and let $R:=\Q[x_1,\ldots,x_{n}]$.   Let $\partial_{1,2,\ldots,n-1}=\partial_1\circ \partial_2\circ \dots \partial_{n-1}$ denote the composition of divided difference operators.  Then
	\[
	\partial_{1,2,\ldots,n-1}(x_{n}^k) = (-1)^{n-1}h_{k+1-n}(x_1,\ldots,x_n)
	\]
	where $h_k$ denotes the $k$-th complete symmetric function.
\end{proposition}
\begin{proof}
	Induction on $n\geq 1$.  For the base case $n=1$, the composition $\partial_1\circ \cdots \circ \partial_{n-1}$ is empty, hence $\partial_{1,2,\ldots,n-1}$ is to be interpreted as the identity.  Then the base case is trivial: $x_1^k = h_k(x_1)$.  Since this case is degenerate, let's also consider the case $n=2$, which will be relevant for the inductive step as well.  Compute:
	\begin{equation}\label{eq-aDemazureFormula}
	\partial_1(x_2^k) = -\frac{x_1^k-x_2^k}{x_1-x_2} = -(x_1^{k-1}+x_1^{k-2}x_2+\cdots + x_2^{k-1})
	\end{equation}
	which equals $-h_{k-1}(x_1,x_2)$, as desired.
	
	Assume by induction that we have proven
	\[
	\partial_{1,2,\ldots,n-2}(x_{n-1}^i) =(-1)^{n} h_{i+2-n}(x_1,\ldots,x_{n-1})
	\]
	Compute:
	\begin{eqnarray*}
		\partial_{1,2,\ldots,n-1}(x_n^k)
		& = & \partial_{1,2,\ldots,n-2}(\partial_{n-1}(x_n^k))\\
		& = & -\partial_{1,2,\ldots,n-2}\Big(\sum_{i+j=k-1}x_{n-1}^ix_n^j\Big)\\
		& = & -\sum_{i+j=k-1}\partial_{1,2,\ldots,n-2}(x_{n-1}^i)x_n^j\\
		& = & -\sum_{i+j=k-1}(-1)^{n} h_{i+2-n}(x_1,\ldots,x_{n-1})x_n^j\\
		& = & (-1)^{n-1}h_{k+1-n}(x_1,\ldots,x_n)
	\end{eqnarray*}
	In the second line we used the formula (\ref{eq-aDemazureFormula}), in the fourth we used the inductive hypothesis, and in the last line we used the well-known recursion satisfied by the complete symmetric functions (see Proposition \ref{prop-symFnRecursion}).  This completes the inductive step.
\end{proof}

\begin{lemma}\label{lemma-RIgeneration}
Recall the sequence of inclusions of algebras $R^W\subset R^I\subset R$.  As an $R^W$-module, $R^I$ is generated by $\{1,x_n,\ldots,x_n^{n-1}\}$.
\end{lemma}
\begin{proof}
	Recall that $R^I\subset \Q[x_1,\ldots,x_n]$ is the algebra of polynomials which are symmetric in $x_1,\ldots,x_{n-1}$.  It is well known that $R^I$ is generated, as a unital algebra, by $x_n$ and the (say) elementary symmetric functions $e_k(x_1,\ldots,x_{n-1})$.  The recursion from Proposition \ref{prop-symFnRecursion} says that
	\[
	e_k(x_1,\ldots,x_{n-1}) = e_k(x_1,\ldots,x_n) - x_ne_{k-1}(x_1,\ldots,x_{n-1})
	\]
	Iterating, we can express each partially symmetric function $e_k(x_1,\ldots,x_{n-1})$ as
	\[
	e_k(x_1,\ldots,x_{n-1}) = \sum_{i=0}^k (-1)^{k-i} e_{i}(x_1,\ldots,x_n)x_n^{k-i}
	\]
	This shows that $R^I$ is generated by the powers of $x_n$, as an $R^W$ algebra.  Recall from Example \ref{example-theRelation} that in $R$ one has the following identity
	\[
	x_n^n=-\sum_{i=1}^n(-1)^{n-i}e_i(x_1,\ldots,x_n)x_n^{n-i}
	\]
	hence $1,x_n,\ldots,x_n^{n-1}$ suffice to generate $R^I$ as an $R^W$-module.
\end{proof}

\begin{proposition}\label{prop-theDemazureOp}
The divided difference operator $\partial_{1,2,\ldots,n-1}:R\rightarrow R$ restricts to an $R^W$-linear map $R^I\rightarrow R^W$.
\end{proposition} 
\begin{proof}
By part (3) of Proposition \ref{prop-dividedDiffProps} it is clear that $\partial_{1,2,\ldots,n-1}$ is $R^W$-linear.  Proposition \ref{prop-traceOfXn} shows that $\partial_{1,2,\ldots,n-1}x_n^k$ is an element of $R^W$, for each $k$.  Since the $x_n^k$ generate $R^I$ as an $R^W$-module (Lemma \ref{lemma-RIgeneration}), the proposition follows.
\end{proof}

We are now ready to prove the main theorem of this section:

\begin{proof}[Proof of Theorem \ref{thm-frobBases}]
In this proof we will use the brief notation $\del = \del_{1,2,\ldots,n-1}$. We will focus on proving that
\begin{equation}\label{eq-traceEX}
\partial\Big((-1)^i x_n^i e_{n-1-j}(x_1,\ldots,x_{n-1})\Big) \ = \ \d_{ij},
\end{equation}
where $\d_{ij}$ is the Kronecker delta.  Assume for a moment that we have proven this.  A standard linear algebra argument will imply that the set $\b:=\{(-1)^j x_n^j\}_{j=0}^{n-1} \subset R^I$ is $R^W$-linearly independent, as is $\b^{\ast}:=\{e_i(x_1,\ldots,x_n)\}_{i=0}^{n-1}$.  Lemma \ref{lemma-RIgeneration} says that $\b$ spans $R^I$ as an $R^W$-module.  Thus, $\b$ forms a basis, and Equation (\ref{eq-traceEX}) implies that $\b^{\ast}$ is the dual basis, up to reordering.

It remains to prove Equation (\ref{eq-traceEX}).  For this, consider the two-variable generating function
\[
F(t,u):=\sum_{i,j=0}^\infty t^iu^j (-1)^j e_i(x_1,\ldots,x_{n-1})x_n^j 
\]
Observe:
\begin{eqnarray*}
F(t,u)
&=& \frac{\prod_{i=1}^{n-1}(1+tx_i)}{1+ux_n}\\
&=& \frac{\prod_{i=1}^{n}(1+tx_i)}{(1+ux_n)(1+tx_n)}\\
&=& \Big(\sum_{i=0}^nt^ie_i(x_1,\ldots,x_n)\Big)\Big(\sum_{j=0}^\infty (-1)^jt^j x_n^j\Big)\Big(\sum_{k=0}^\infty (-1)^ku^k x_n^k\Big)\\
&=& \sum_{i,j,k}(-1)^{j+k} t^{i+j}u^k e_i(x_1,\ldots,x_n) x_n^{j+k}
\end{eqnarray*}
Now, taking $\partial$ gives
\begin{eqnarray*}
\partial(F(t,u))
&=& \sum_{i,j,k}(-1)^{j+k} t^{i+j}u^k \partial_{1,2,\ldots,n-1}(e_i(x_1,\ldots,x_n) x_n^{j+k})\\
&=& \sum_{i,j,k}(-1)^{j+k} t^{i+j}u^k e_i(x_1,\ldots,x_n) \partial_{1,2,\ldots,n-1}(x_n^{j+k})\\
&=& \sum_{i,j,k}(-1)^{j+k} t^{i+j}u^k e_i(x_1,\ldots,x_n) (-1)^{n-1}h_{k+j+1-n}(x_1,\ldots,x_n)\\
&=& \sum_{m,k}t^mu^k \sum_{a+b=m+k+1-n} (-1)^{b} e_a(x_1,\ldots,x_n) h_{b}(x_1,\ldots,x_n)\\
&=& \sum_{m,k}t^mu^k \d_{m+k,n-1}
\end{eqnarray*}
In the second line we used the fact that the divided difference operators are $R^W$-linear.  In the third line we used Proposition \ref{prop-traceOfXn} to compute $\partial(x_n^{j+k})$.  In the fourth line we introduced the new indices $m=i+j$, $a=i$, and $b=k+j+1-n$.  Finally, the last line follows from a well known identity satisfied by the elementary and complete symmetric functions.  This computation gives the desired formula.
\end{proof}

\subsection{An explicit family of relations in $R\otimes_{R^W} R$}
\label{subsec-explicitRels}
Throughout this section we ignore all grading shifts.  Let $I_k\subset \Q[\xx,\yy]$ denote the ideal from Notation \ref{notation-XandY}, so that $B_{w_0}\cong\Q[\xx,\yy]/I_n$ and $B_{w_1}\cong \Q[\xx,\yy]/I_{n-1}$.  In this section we prove Proposition \ref{prop-aijRelation}.

First, recall that the elementary symmetric functions $e_k(x_1,\ldots,x_n)\in R$ and the complete symmetric functions $h_k(x_1,\ldots,x_n)\in R$ are defined by their \emph{generating functions}:
\begin{equation}\label{eq-stdGeneratingFns}
\sum_{k=0}^n t^k e_k(x_1,\ldots,x_n) = \prod_{i=1}^n(1+tx_i) \ \ \ \ \ \ \ \ \ \ \sum_{k=0}^\infty t^kh_k(x_1,\ldots,x_n) = \frac{1}{\prod_{i=1}^n (1-tx_i)}
\end{equation}
where $t$ is a formal indeterminate.  These expressions make the following recursion relations obvious:
\begin{proposition}\label{prop-symFnRecursion}
	We have the following recursive relations, for $1\leq k\leq n-1$:
	\begin{enumerate}
		\item $e_k(x_1,\ldots,x_{n+1}) = e_k(x_1,\ldots,x_n) + x_ne_{k-1}(x_1,\ldots,x_n)$.
		\item $h_k(x_1,\ldots,x_{n+1}) = \sum_{i=0}^kx_n^i h_{k-i}(x_1,\ldots,x_n)$.
	\end{enumerate}\qed
\end{proposition}

\begin{proposition}\label{prop-theRelations}
	Fix integers  $1\leq m\leq n$.  The polynomials
	\[
	z_{m,n}(x_m,\ldots,x_n,y_1,\ldots,y_n):=\sum_{i+j=m}(-1)^je_i(y_1,\ldots,y_n)h_j(x_m,\ldots,x_n) \in \Q[\xx,\yy]
	\]
	are zero modulo $I_n$.
\end{proposition}
\begin{proof}
Fix an integer $m\in\{1,\ldots,n\}$, and consider the following formal power series in $t$:
	\[
	Z(t):=\frac{\prod_{i=1}^n(1+ty_i)}{\prod_{j=m}^n (1+tx_i)}
	\]
	The numerator is completely symmetric in the variables $y_1,\ldots,y_n$ hence each $y_k$ can be replaced by $x_k$ (modulo $I_n$).  The result is a polynomial of degree $m-1$:
	\begin{equation}\label{eq-ZmodI}
	Z(t) \equiv \prod_{i=1}^{m-1} (1+tx_i) \ \ \ \ (\text{mod }I_n).
	\end{equation}
	On the other hand, by (\ref{eq-stdGeneratingFns}) the generating function $Z(t)$ can be written as
	\[
	Z(t) = \sum_{k=0}^\infty t^k \sum_{i+j=k}(-1)^je_i(y_1,\ldots,y_n)h_j(x_m,\ldots,x_n)
	\]
	By (\ref{eq-ZmodI}, we see that the coefficient on $t^k$ vanishes for $k\geq m$, modulo $\IC_n$.  That is to say, for each $k\geq m$ we have
	\[
	\sum_{i+j=k}(-1)^je_i(y_1,\ldots,y_n)h_j(x_m,\ldots,x_n)\equiv 0  \ \ \ \ \ \ \ \ (\text{mod }I_n)
	\]
	Specializing to $k=m$ gives the identity in the statement.
\end{proof}

\begin{example}\label{example-theRelation}  Setting $m=n$ in Proposition \ref{prop-theRelations}, we obtain that
	\[
	\sum_{i+j=n} (-1)^j e_i(y_1,\ldots,y_n)x_n^j=0  \ \ \ \ \ \ \ \ \ \ \ \ \ \ \ \ \ (\text{mod }I_n)
	\]
	We can factor this as $\prod_{i=1}^n(y_i-x_n)=0$ (mod $I_n$), so we recover the relations from Proposition \ref{prop-theRelations} as a special case.  Let $I_1\subset \Q[x_1,\ldots,x_n,y_1,\ldots,y_n]$ denote the ideal generated by the differences $y_i-x_i$ for $1\leq i\leq n$.  Viewing the above relation in the quotient $R\cong \Q[x_1,\ldots,x_n,y_1,\ldots,y_n]/I_1$, we find that
	\begin{equation}\label{eq-xnRelation}
	x_n^n=-e_n(x_1,\ldots,x_n) + x_ne_{n-1}(x_1,\ldots,x_n)-\cdots +(-1)^{n} x_n^{n-1} e_1(x_1,\ldots,x_n)
	\end{equation}
	in $\Q[x_1,\ldots,x_n]$.
\end{example}

The following gives an explicit formula for the relations $z_{m,n}$ in $R\otimes_{R^W} R$.

\begin{proposition}\label{prop-explicitRelations}
Fix an integer $n\geq 1$.  The identity, 
\[
z_{m,n} = \sum_{\gamma}\prod_{i=1}^{m}(y_{\gamma_i} - x_{\gamma_i+m-i})
%
\]
,holds in $\Q[x_1,\ldots,x_n,y_1,\ldots,y_n]$, where the sum is over sequences of integers $\gamma=(\gamma_1,\ldots,\gamma_{m})$ with $1\leq \gamma_1<\cdots <\gamma_{m}\leq n$.
\end{proposition}
Before proving, we give some examples which we hope will illustrate the formula.  First:
\begin{definition}\label{def-pij}
Set $p_{ij}=y_i-x_j\in \Q[\xx,\yy]$.  In all ocurrences in this paper, the variables will satisfy $1\leq i\leq j\leq n$.
\end{definition}

\begin{example}\label{example-lotsOrelations}
At one extreme, we have $z_{1,n}=p_{11}+p_{22}+\cdots + p_{nn}$, which equals $e_1(\yy)-e_1(\xx)\in\Q[\xx,\yy]$.  At the other extreme we have $z_{n,n}=p_{1n}p_{2n}\cdots p_{nn}$.  The relations $z_{mn}$ with $2\leq m\leq n\leq 4$ are:
\begin{enumerate}
\item $z_{22}=p_{12}p_{22}$
\item $z_{23}=p_{12}p_{22} + (p_{12}+p_{23})p_{33}$
\item $z_{24}=p_{12}p_{22} + (p_{12}+p_{23})p_{33}+ (p_{12}+p_{23}+p_{34})p_{44}$
\item $z_{33} = p_{13}p_{23}p_{33}$
\item $z_{34} = p_{13}p_{23}p_{33}+(p_{13}p_{23}+p_{13}p_{34}+p_{24}p_{34})p_{44}$
\item $z_{44} = p_{14}p_{24}p_{34}p_{44}$
\end{enumerate}
\end{example}

\begin{proof}[Proof of Proposition \ref{prop-explicitRelations}]
Observe:
\begin{eqnarray*}
\prod_{1\leq i\leq \ell \atop m\leq j\leq \ell}\frac{1+ty_i}{1+tx_j} - \prod_{1\leq i\leq \ell-1 \atop m\leq j\leq \ell-1}\frac{1+ty_i}{1+tx_j}
&=& \prod_{1\leq i\leq \ell-1 \atop m\leq j\leq \ell-1}\frac{1+ty_i}{1+tx_j}\bigg(\frac{1+ty_{\ell}}{1+tx_{\ell}}-1\bigg)\\
&=& \prod_{1\leq i\leq \ell-1 \atop m\leq j\leq \ell-1}\frac{1+ty_i}{1+tx_j}\bigg(\frac{t(y_{\ell}-x_{\ell})}{1+tx_{\ell}}\bigg)\\
&=& \prod_{1\leq i\leq \ell-1 \atop m\leq j\leq \ell}\frac{1+ty_i}{1+tx_j} t(y_{\ell}-x_{\ell})
\end{eqnarray*}
Adding up the resulting identities for $m\leq \ell\leq n$ gives the identity
\[
\prod_{1\leq i\leq n \atop m\leq j\leq n}\frac{1+ty_i}{1+tx_j}\ \ - \ \ \prod_{1\leq i\leq m-1}(1+ty_i) \ \ =\ \  t\sum_{\ell=m}^n(y_{\ell}-x_{\ell})\prod_{1\leq i\leq \ell-1 \atop m\leq j\leq \ell}\frac{1+ty_i}{1+tx_j}
\]
We can rewrite this as
\[
\prod_{1\leq i\leq n \atop m\leq j\leq n}\frac{1+ty_i}{1+tx_j}\ \  =\ \  t\sum_{\ell=m}^n(y_{\ell}-x_{\ell})\prod_{1\leq i\leq \ell-1 \atop m-1\leq j\leq \ell-1}\frac{1+ty_i}{1+tx_{j+1}} + (\text{lower})
\]
where (lower) denotes terms of $t$-degree $<m$. Iterating this formula, we find that
\begin{eqnarray*}
&&\prod_{1\leq i\leq n\atop m\leq j\leq n}\frac{1+ty_i}{1+tx_j}  \\
&& \ \ \ \ \ \ \ \ = t^m\sum_{k_1=m}^n (y_{k_1}-x_{k_1})\sum_{k_2=m-1}^{k_1-1}(y_{k_2}-x_{k_2+1})\cdots\sum_{k_m=1}^{k_{m-1}-1} (y_{k_m}-x_{k_m+m-1})
\end{eqnarray*}
plus terms of $t$-degree $\neq m$.  Comparing the coefficients on $t^m$ on either side gives the identity in the statement.
\end{proof}
Recall the polynomials from Definition \ref{def-aij}:
\begin{equation}\label{eqn-aPolys}
a_{mj}(\xx,\yy):=\sum_{\gamma'} \prod_{i=1}^{m-1}(y_{\gamma_i} - x_{\gamma_i+m-i}),
\end{equation}
where $1\leq m\leq j\leq n$. The sum is over sequences of integers $\gamma' = (\gamma_1,\ldots,\gamma_{m-1})$ with $1\leq \gamma_1<\cdots <\gamma_{m-1}<j$.
\begin{aijRelation}
The $a_{ij}(\xx,\yy)$ satisfy $\sum_{j=k}^n a_{i,j}(\xx,\yy)(x_j-y_j) =0$ modulo $I_n$.  In particular $\sum_{j=k}^n a_{i,j}(\xx,\yy)(x_j-y_j) $ acts by zero on $B_{w_0}$.
\end{aijRelation}
\begin{proof}
Observe
\begin{eqnarray*}
\sum_{j=m}^n a_{mj} \cdot (y_j - x_j)  &=& \sum_{j=m}^n \sum_{\gamma} (y_j - x_j) \prod_{i=1}^{m=1}(y_{\gamma_i} - x_{\gamma_i+m-i})\\
& =& \sum_{\gamma} \prod_{i=1}^{m}(y_{\gamma_i} - x_{\gamma_i+m-i})\\
&=& z_{m,n}
\end{eqnarray*}
The first equality holds by definition.  In the second, we have reindexed by setting $\gamma=(\gamma_1,\ldots,\gamma_{m-1},j)$.  The third equality holds by Proposition \ref{prop-explicitRelations}.
\end{proof}

\bibliographystyle{plain}
\bibliography{bib}

\end{document}